\documentclass[thmsa,onecolumn,10pt]{article}
\usepackage{amssymb}
\usepackage{amsfonts}
\usepackage{graphicx}
\usepackage{amsmath}

\newtheorem{theorem}{Theorem}

\newtheorem{corollary}[theorem]{Corollary}

\newtheorem{lemma}[theorem]{Lemma}

\newtheorem{proposition}[theorem]{Proposition}

\newenvironment{proof}[1][Proof]{\textbf{#1.} }{\ \rule{0.5em}{0.5em}}

\begin{document}

\title{Stable convergence of generalized stochastic integrals and the
principle of conditioning: $L^{2}$ theory}
\author{Giovanni PECCATI\thanks{%
Laboratoire de Statistique Th\'{e}orique et Appliqu\'{e}e, Universit\'{e}
Paris VI, France. E-mail: \texttt{giovanni.peccati@gmail.com}}, and Murad S.
TAQQU\thanks{%
Boston University, Departement of Mathematics, 111 Cummington
Road, Boston (MA), USA. E-mail: murad@math.bu.edu. }.
\thanks{This research was partially supported by the NSF Grant
DNS-050547 at Boston University.}}
\date{April 25, 2006}
\maketitle

\begin{abstract}
Consider generalized adapted stochastic integrals with respect to
independently scattered random measures with second moments. We use a
decoupling technique, known as the \textquotedblleft principle of
conditioning\textquotedblright , to study their stable convergence towards
mixtures of infinitely divisible distributions. Our results apply, in
particular, to multiple integrals with respect to independently scattered
and square integrable random measures, as well as to Skorohod integrals on
abstract Wiener spaces. As a specific application, we establish a Central
Limit Theorem for sequences of double integrals with respect to a general
Poisson measure, thus extending the results contained in Nualart and Peccati
(2005) and Peccati and Tudor (2004) to a non-Gaussian context.


\end{abstract}

\section{Introduction}

In this paper we establish several criteria, ensuring the stable convergence
of sequences of \textquotedblleft generalized integrals\textquotedblright\
with respect to independently scattered random measures over abstract
Hilbert spaces. The notion of generalized integral is understood in a very
wide sense, and includes for example Skorohod integrals with respect to
isonormal Gaussian processes (see e.g. \cite{Nualart}), multiple Wiener-It%
\^{o} integrals associated to general Poisson measures (see \cite{NV}, or
\cite{Kab}), or the class of iterated integrals with respect to
orthogonalized Teugels martingales introduced in \cite{NuSch}. All these
random objects can be represented as appropriate generalized
\textquotedblleft adapted stochastic integrals\textquotedblright\ with
respect to a (possibly infinite) family of L\'{e}vy processes, constructed
by means of a well-chosen increasing family of orthogonal projections. These
adapted integrals are also the limit of sums of arrays of random variables
with a special dependence structure. We shall show, in particular, that
their asymptotic behavior can be naturally studied by means of a decoupling
technique, known as the \textquotedblleft principle of
conditioning\textquotedblright\ (see e.g. \cite{Jakubowki} and \cite{Xue}),
that we use in the framework of stable convergence (see \cite[Chapter 4]%
{JacSh}).

\bigskip

Our setup is roughly the following. We shall consider a centered and square
integrable random field $X=\left\{ X\left( h\right) :h\in \mathfrak{H}%
\right\} $, indexed by a separable Hilbert space $\mathfrak{H}$, and
verifying the isomorphic relation $\mathbb{E}\left[ X\left( h\right) X\left(
h^{\prime }\right) \right] $ $=$ $\left( h,h^{\prime }\right) _{\mathfrak{H}%
} $, where $\left( \cdot ,\cdot \right) _{\mathfrak{H}}$ is the inner
product on $\mathfrak{H}$. There is no time involved. To introduce time,
endow the space $\mathfrak{H}$ with an increasing family of orthogonal
projections, say $\pi _{t}$, $t\in \left[ 0,1\right] $, such that $\pi
_{0}=0 $ and $\pi _{1}=$ id.. Such projections operators induce the
(canonical) filtration $\mathcal{F}^{\pi }=\left\{ \mathcal{F}_{t}^{\pi
}:t\in \left[ 0,1\right] \right\} $, where each $\mathcal{F}_{t}^{\pi }$ is
generated by random variables of the type $X\left( \pi _{t}h\right) $, and
one can define (e.g., as in \cite{UZ} for Gaussian processes) a class of $%
\mathcal{F}^{\pi } $-adapted and $\mathfrak{H}$-valued random variables. If
for every $h\in \mathfrak{H}$ the application $t\mapsto X\left( \pi
_{t}h\right) $ is also a $\mathcal{F}^{\pi }$-L\'{e}vy process, then there
exists a natural It\^{o} type stochastic integral, of adapted and $\mathfrak{%
H}$-valued variables, with respect to the infinite dimensional process $%
t\mapsto \left\{ X\left( \pi _{t}h\right) :h\in \mathfrak{H}\right\} $.
Denote by $J_{X}\left( u\right) $ the integral of an adapted random variable
$u$ with respect to $X$. As will be made clear in the subsequent discussion,
several random objects appearing in stochastic analysis (such as Skorohod
integrals, or the multiple Poisson integrals quoted above) are in fact
generalized adapted integrals of the type $J_{X}\left( u\right) $, for some
well chosen random field $X$. Moreover, the definition of $J_{X}\left(
u\right) $ mimics in many instances the usual construction of adapted
stochastic integrals with respect to real-valued martingales. In particular:
(i) each stochastic integral $J_{X}\left( u\right) $ is associated to a $%
\mathcal{F}^{\pi }$-martingale, namely the process $t\mapsto J_{X}\left( \pi
_{t}u\right) $ and (ii) $J_{X}\left( u\right) $ is the limit (in $L^{2}$) of
finite \textquotedblleft adapted Riemann sums\textquotedblright\ of the kind
$S\left( u\right) =\sum_{j=1,...,n}F_{j}X\left( \left( \pi _{t_{j+1}}-\pi
_{t_{j}}\right) h_{j}\right) $, where $h_{j}\in \mathfrak{H}$, $%
t_{n}>t_{n-1} $ $>\cdot \cdot \cdot >t_{1}$ and $F_{j}\in \mathcal{F}%
_{t_{j}}^{\pi }$. We show that, by using a decoupling result known as
\textquotedblleft principle of conditioning\textquotedblright\ (see Theorem
1 in \cite{Xue}, and Section 2 below, for a very general form of such
principle), the stable and, in particular, the weak convergence of sequences
of sums such as $S\left( u\right) $ is completely determined by the
asymptotic behavior of random variables of the type%
\begin{equation*}
\widetilde{S}\left( u\right) =\sum_{j=1,...,n}F_{j}\widetilde{X}\left(
\left( \pi _{t_{j+1}}-\pi _{t_{j}}\right) h_{j}\right) \text{,}
\end{equation*}%
where $\widetilde{X}$ is an independent copy of $X$. Note that the vector
\begin{equation*}
\widetilde{V}=\left( F_{1}\widetilde{X}\left( \left( \pi _{t_{2}}-\pi
_{t_{1}}\right) h_{1}\right) ,...,F_{n}\widetilde{X}\left( \left( \pi
_{t_{n+1}}-\pi _{t_{n}}\right) h_{n}\right) \right) ,
\end{equation*}%
enjoys the specific property of being \textit{decoupled }(i.e.,
conditionally on the $F_{j}$'s, its components are independent) and \textit{%
tangent }to the \textquotedblleft original\textquotedblright\ vector%
\begin{equation*}
V=\left( F_{1}X\left( \left( \pi _{t_{2}}-\pi _{t_{1}}\right) h_{1}\right)
,...,F_{n}X\left( \left( \pi _{t_{n+1}}-\pi _{t_{n}}\right) h_{n}\right)
\right) \text{,}
\end{equation*}%
in the sense that for every $j$, and conditionally on the r.v.'s $F_{k}$, $%
k\leq j$, $F_{j}X\left( \left( \pi _{t_{j+1}}-\pi _{t_{j}}\right)
h_{j}\right) $ and $F_{j}\widetilde{X}\left( \left( \pi _{t_{j+1}}-\pi
_{t_{j}}\right) h_{j}\right) $ have the same law (the reader is referred to
\cite{JacSadi} or \cite{KW91} for a discussion of the general theory of
tangent processes). The convergence of sequences such as $J_{X}\left(
u_{n}\right) $, $n\geq 1$, where each $u_{n}$ is adapted, can therefore be
studied by means of simpler random variables $\widetilde{J}_{X}\left(
u_{n}\right) $, obtained from a decoupled and tangent version of the
martingale $t\mapsto J_{X}\left( \pi _{t}u_{n}\right) $. In particular (see
Theorem \ref{T : Main} below, as well as its consequences) we shall prove
that, since such decoupled processes can be shown to have conditionally
independent increments, the problem of the stable convergence of $%
J_{X}\left( u_{n}\right) $ can be reduced to the study of the convergence in
probability of sequences of random L\'{e}vy-Khinchine exponents. This
represents an extension of the techniques initiated in \cite{NuPe} and \cite%
{PT} where, in a purely Gaussian context, the CLTs for multiple Wiener-It%
\^{o} integrals are characterized by means of the convergence in probability
of the quadratic variation of Brownian martingales. We remark that the
extensions of \cite{NuPe} and \cite{PT} achieved in this paper go in two
directions: (a) we consider general (not necessarily Gaussian) square
integrable and independently scattered random measures, (b) we study stable
convergence, instead of weak convergence, so that, for instance, our results
can be used in the Gaussian case to obtain non-central limit theorems (see
e.g. Section 6 below, as well as \cite{PeTaq}).

\bigskip

When studying the stable convergence of random variables that are terminal
values of continuous-time martingales, one could alternatively use the
general criteria for the stable convergence of semimartingales, as developed
e.g. in \cite{LiSh}, \cite{Feigin} or \cite[Chapter 4]{JacSh}, instead of
the above decoupling techniques. However, the principle of conditioning
(which is in some sense the discrete-time skeleton of the general
semimartingale results), as formulated in the present paper, often requires
less stringent assumptions. For instance, conditions (\ref{emboite}) and (%
\ref{boxed}) below are weak versions of the \textit{nesting condition}
introduced by Feigin in the classic reference \cite{Feigin}.

\bigskip

The paper is organized as follows. In Section 2, we discuss a general
version of the principle of conditioning. In Section 3 we present a general
setup to which such decoupling techniques can be applied, and in Section 4
the above mentioned convergence results are established. In Section 5.1 and
5.2, we apply our techniques to sequences of multiple stochastic integrals
with respect to independently scattered random measures with second moments,
whereas in Section 3.3 we give a specific application to Central Limit
Theorems for double Poisson integrals. Finally, in Section 6 our results are
applied to study the stable convergence of Skorohod integrals with respect
to a general isonormal Gaussian process.

\section{The principle of conditioning}

We shall present a general version of the \textit{principle of conditioning}
(\textit{POC} in the sequel) for arrays of real valued random variables. Our
discussion is mainly inspired by a remarkable paper by X.-H. Xue \cite{Xue},
generalizing the classic results by Jakubowski \cite{Jakubowki} to the
framework of stable convergence. Note that the results discussed below refer
to a discrete time setting. However, thanks to some density arguments, we
will be able to apply most of the POC techniques to general stochastic
measures on abstract Hilbert spaces.

Instead of adopting the formalism of \cite{Xue} we choose, for the sake of
clarity, to rely in part on the slightly different language of \cite[Ch. 6
and 7]{GineDelaPena}. To this end, we shall recall some notions concerning
stable convergence, conditional independence and decoupled sequences of
random variables. From now on, all random objects are supposed to be defined
on an adequate probability space $\left( \Omega ,\mathcal{F},\mathbb{P}%
\right) $, and all $\sigma $-fields introduced below will be tacitly assumed
to be complete; $\overset{\mathbb{P}}{\rightarrow }$ means convergence in
probability; $\mathbf{%
\mathbb{R}
}$ stands for the set of real numbers; $\triangleq $ denotes a new
definition.

\bigskip

We start by defining the class $\mathbf{M}$ of random probability measures,
and the class $\widehat{\mathbf{M}}$ (resp. $\widehat{\mathbf{M}}_{0}$) of
random (resp. non-vanishing and random) characteristic functions.

\bigskip

\textbf{Definition A }(see e.g. \cite{Xue})\ -- Let $\mathcal{B}\left(
\mathbf{%
\mathbb{R}
}\right) $ denote the Borel $\sigma $-field on $\mathbf{%
\mathbb{R}
}$.

\begin{description}
\item[(A-i)] A map $\mu \left( \cdot ,\cdot \right) $, from $\mathcal{B}%
\left( \mathbf{%
\mathbb{R}
}\right) \times \Omega $ to $\mathbf{%
\mathbb{R}
}$ is called a \textit{random probability }(on $\mathbf{%
\mathbb{R}
}$) if, for every $C\in \mathcal{B}\left( \mathbf{%
\mathbb{R}
}\right) $, $\mu \left( C,\cdot \right) $ is a random variable and, for $%
\mathbb{P}$-a.e. $\omega $, the map $C\mapsto \mu \left( C,\omega \right) $,
$C\in \mathcal{B}\left( \mathbf{%
\mathbb{R}
}\right) $, defines a probability measure on $\mathbf{%
\mathbb{R}
}$. The class of all random probabilities is noted $\mathbf{M}$, and, for $%
\mu \in \mathbf{M}$, we write $\mathbb{E}\mu \left( \cdot \right) $ to
indicate the (deterministic) probability measure
\begin{equation}
\mathbb{E}\mu \left( C\right) \triangleq \mathbb{E}\left[ \mu \left( C,\cdot
\right) \right] \text{, \ \ }C\in \mathcal{B}\left( \mathbf{%
\mathbb{R}
}\right) .  \label{mixt}
\end{equation}

\item[(A-ii)] For a measurable map $\phi \left( \cdot ,\cdot \right) $, from
$\mathbf{%
\mathbb{R}
}\times \Omega $ to $\mathbb{C}$, we write $\phi \in \widehat{\mathbf{M}}$
whenever there exists $\mu \in \mathbf{M}$ such that
\begin{equation}
\phi \left( \lambda ,\omega \right) =\widehat{\mu }\left( \lambda \right)
\left( \omega \right) \text{, \ \ }\forall \lambda \in \mathbf{%
\mathbb{R}
}\text{, for }\mathbb{P}\text{-a.e. }\omega \text{,}  \label{randomfou}
\end{equation}%
where $\widehat{\mu }\left( \cdot \right) $ is defined as
\begin{equation}
\widehat{\mu }\left( \lambda \right) \left( \omega \right) =\left\{
\begin{array}{ll}
\int \exp \left( i\lambda x\right) \mu \left( dx,\omega \right) & \text{if }%
\mu \left( \cdot ,\omega \right) \text{ is a probability measure} \\
1 & \text{otherwise.}%
\end{array}%
\right. \text{, \ \ }\lambda \in \mathbf{%
\mathbb{R}
}.  \label{ftr}
\end{equation}

\item[(A-iii)] For a given $\phi \in \widehat{\mathbf{M}}$, we write $\phi
\in \widehat{\mathbf{M}}_{0}$ whenever $\mathbb{P}\left\{ \omega :\phi
\left( \lambda ,\omega \right) \neq 0\text{ \ }\forall \lambda \in \mathbf{%
\mathbb{R}
}\right\} =1$.
\end{description}

\bigskip

Observe that, for every $\omega \in \Omega $, $\widehat{\mu }\left( \lambda
\right) \left( \omega \right) $ is a continuous function of $\lambda $. The
probability $\mathbb{E}\mu \left( \cdot \right) =\int_{\Omega }\mu \left(
\cdot ,\omega \right) d\mathbb{P}\left( \omega \right) $ defined in (\ref%
{mixt}) is often called a \textit{mixture }of probability measures. The
following definition of \textit{stable convergence }extends the usual notion
of convergence in law.

\bigskip

\textbf{Definition B }(see e.g. \cite[Chapter 4]{JacSh} or \cite{Xue}) --
Let $\mathcal{F}^{\ast }\subseteq \mathcal{F}$ be a $\sigma $-field, and let
$\mu \in \mathbf{M}$. A sequence of real valued r.v.'s $\left\{ X_{n}:n\geq
1\right\} $ is said to \textit{converge }$\mathcal{F}^{\ast }$-\textit{%
stably }to $\mathbb{E}\mu \left( \cdot \right) $, written $X_{n}\rightarrow
_{\left( s,\mathcal{F}^{\ast }\right) }\mathbb{E}\mu \left( \cdot \right) $,
if, for every $\lambda \in \mathbf{%
\mathbb{R}
}$ and every bounded $\mathcal{F}^{\ast }$-measurable r.v. $Z$,%
\begin{equation}
\lim_{n\rightarrow +\infty }\mathbb{E}\left[ Z\times \exp \left( i\lambda
X_{n}\right) \right] =\mathbb{E}\left[ Z\times \widehat{\mu }\left( \lambda
\right) \right] \text{,}  \label{stable cond}
\end{equation}%
where the notation is the same as in (\ref{ftr}).

\bigskip

If $X_{n}$ converges $\mathcal{F}^{\ast }$-stably, then the conditional
distributions $\mathcal{L}\left( X_{n}\mid A\right) $ converge for any $A\in
\mathcal{F}^{\ast }$ such that $\mathbb{P}\left( A\right) >0$. (see e.g.
\cite[Section 5, \S 5c]{JacSh} for further characterizations of stable
convergence). Note that, by setting $Z=1$, we obtain that if $%
X_{n}\rightarrow _{\left( s,\mathcal{F}^{\ast }\right) }\mathbb{E}\mu \left(
\cdot \right) $, then the law of the $X_{n}$'s converges weakly to $\mathbb{E%
}\mu \left( \cdot \right) $. Moreover, by a monotone class argument, $%
X_{n}\rightarrow _{\left( s,\mathcal{F}^{\ast }\right) }\mathbb{E}\mu \left(
\cdot \right) $ if, and only if, (\ref{stable cond}) holds for random
variables with the form $Z=\exp \left( i\gamma Y\right) $, where $\gamma \in
\mathbf{%
\mathbb{R}
}$ and $Y$ is $\mathcal{F}^{\ast }$-measurable. Eventually, we note that, if
a sequence of random variables $\left\{ U_{n}:n\geq 0\right\} $ is such that
$\left( U_{n}-X_{n}\right) \rightarrow 0$ in $L^{1}\left( \mathbb{P}\right) $
and $X_{n}\rightarrow _{\left( s,\mathcal{F}^{\ast }\right) }\mathbb{E}\mu
\left( \cdot \right) $, then $U_{n}\rightarrow _{\left( s,\mathcal{F}^{\ast
}\right) }\mathbb{E}\mu \left( \cdot \right) $. The following definition
shows how to replace an array $X^{\left( 1\right) }$ of real-valued random
variables by a simpler, \textit{decoupled} array $X^{\left( 2\right) }$.

\bigskip

\textbf{Definition C }(see \cite[Chapter 7]{GineDelaPena}) \textbf{-- }Let $%
\left\{ N_{n}:n\geq 1\right\} $ be a sequence of positive natural numbers,
and let
\begin{equation*}
X^{\left( i\right) }\triangleq \left\{ X_{n,j}^{\left( i\right) }:0\leq
j\leq N_{n},\text{ }n\geq 1\right\} \text{, \ with\ }X_{n,0}^{\left(
i\right) }=0\text{,}
\end{equation*}%
$i=1,2$, be two arrays of real valued r.v.'s, such that, for $i=1,2$ and for
each $n$, the sequence%
\begin{equation*}
X_{n}^{\left( i\right) }\triangleq \left\{ X_{n,j}^{\left( i\right) }:0\leq
j\leq N_{n}\right\}
\end{equation*}%
is adapted to a discrete filtration $\left\{ \mathcal{F}_{n,j}:0\leq j\leq
N_{n}\right\} $ $\left( \text{of course, }\mathcal{F}_{n,j}\subseteq
\mathcal{F}\right) $. For a given $n\geq 1$, we say that $X_{n}^{\left(
2\right) }$ is a \textit{decoupled tangent sequence} to $X_{n}^{\left(
1\right) }$ if the following two conditions are verified:

\begin{description}
\item[$\bigstar $] (\textit{Tangency}) for each $j=1,...,N_{n}$%
\begin{equation}
\mathbb{E}\left[ \exp \left( i\lambda X_{n,j}^{\left( 1\right) }\right) \mid
\mathcal{F}_{n,j-1}\right] =\mathbb{E}\left[ \exp \left( i\lambda
X_{n,j}^{\left( 2\right) }\right) \mid \mathcal{F}_{n,j-1}\right]
\label{samedis1}
\end{equation}%
for each $\lambda \in \mathbf{%
\mathbb{R}
}$, a.s.-$\mathbb{P}$;

\item[$\bigstar $] (\textit{Conditional independence})\textit{\ }there
exists a $\sigma $-field $\mathcal{G}_{n}\subseteq \mathcal{F}$ such that,
for each $j=1,...,N_{n}$,%
\begin{equation}
\mathbb{E}\left[ \exp \left( i\lambda X_{n,j}^{\left( 2\right) }\right) \mid
\mathcal{F}_{n,j-1}\right] =\mathbb{E}\left[ \exp \left( i\lambda
X_{n,j}^{\left( 2\right) }\right) \mid \mathcal{G}_{n}\right]
\label{samedis2}
\end{equation}%
for each $\lambda \in \mathbf{%
\mathbb{R}
}$, a.s.-$\mathbb{P}$, and the random variables $X_{n,1}^{\left( 2\right)
},...,X_{n,N_{n}}^{\left( 2\right) }$ are conditionally independent given $%
\mathcal{G}_{n}$.
\end{description}

Observe that, in (\ref{samedis1}), $\mathcal{F}_{n,j-1}$ depends on $j$, but
$\mathcal{G}_{n}$ does not. The array $X^{\left( 2\right) }$ is said to be a
\textit{decoupled tangent array} to $X^{\left( 1\right) }$ if $X_{n}^{\left(
2\right) }$ is a decoupled tangent sequence to $X_{n}^{\left( 1\right) }$
for each $n\geq 1$.

\bigskip

\textbf{Remark -- }In general, given $X^{\left( 1\right) }$ as above, there
exists a canonical way to construct an array $X^{\left( 2\right) }$, which
is decoupled and tangent to $X^{\left( 1\right) }$. The reader is referred
to \cite[Section 2 and 3]{KW91} for a detailed discussion of this point, as
well as other relevant properties of decoupled tangent sequences.

\bigskip

The following result is essentially a translation of Theorem 2.1 in \cite%
{Xue} into the language of this section. It is a \textquotedblleft stable
convergence generalization\textquotedblright\ of the results obtained by
Jakubowski in \cite{Jakubowki}.

\bigskip

\begin{theorem}[Xue, 1991]
\label{T : Xue}Let $X^{\left( 2\right) }$ be a decoupled tangent array to $%
X^{\left( 1\right) }$, and let the notation of Definition C prevail (in
particular, the collection of $\sigma $-fields $\left\{ \mathcal{F}_{n,j},%
\mathcal{G}_{n}:0\leq j\leq N_{n},n\geq 1\right\} $ satisfies (\ref{samedis1}%
) and (\ref{samedis2})). We write, for every $n$ and every $k=0,...,N_{n}$, $%
S_{n,k}^{\left( i\right) }\triangleq \sum_{j=0,...,k}X_{n,j}^{\left(
i\right) }$, $i=1,2$. Suppose that there exists a sequence $\left\{
r_{n}:n\geq 1\right\} \subset \mathbb{N}$, and a sequence of $\sigma $%
-fields $\left\{ \mathcal{V}_{n}:n\geq 1\right\} $ such that%
\begin{equation}
\mathcal{V}_{n}\subseteq \mathcal{F}\text{ \ \ and \ }\mathcal{V}%
_{n}\subseteq \mathcal{V}_{n+1}\cap \mathcal{F}_{n,r_{n}}\text{, \ }n\geq 1%
\text{,}  \label{emboite}
\end{equation}%
and, as $n\rightarrow +\infty $,%
\begin{equation}
S_{n,r_{n}\wedge N_{n}}^{\left( 1\right) }\overset{\mathbb{P}}{\rightarrow }0%
\text{, \ \ }\mathbb{E}\left[ \exp \left( i\lambda S_{n,r_{n}\wedge
N_{n}}^{\left( 2\right) }\right) \mid \mathcal{G}_{n}\right] \overset{%
\mathbb{P}}{\rightarrow }1.  \label{CP1}
\end{equation}%
If moreover
\begin{equation}
\mathbb{E}\left[ \exp \left( i\lambda S_{n,N_{n}}^{\left( 2\right) }\right)
\mid \mathcal{G}_{n}\right] \overset{\mathbb{P}}{\rightarrow }\phi \left(
\lambda \right) =\phi \left( \lambda ,\omega \right) ,\text{ \ \ }\forall
\lambda \in \mathbf{%
\mathbb{R}
}\text{,}  \label{CP2}
\end{equation}%
where $\phi \in \widehat{\mathbf{M}}_{0}$ and, $\forall \lambda \in \mathbf{%
\mathbb{R}
}$, $\phi \left( \lambda \right) \in \vee _{n}\mathcal{V}_{n}$, then, as $%
n\rightarrow +\infty $,
\begin{equation}
\mathbb{E}\left[ \exp \left( i\lambda S_{n,N_{n}}^{\left( 1\right) }\right)
\mid \mathcal{F}_{n,r_{n}}\right] \overset{\mathbb{P}}{\rightarrow }\phi
\left( \lambda \right) ,\text{ \ \ }\forall \lambda \in \mathbf{%
\mathbb{R}
}\text{,}  \label{fourconv}
\end{equation}%
and%
\begin{equation}
S_{n,N_{n}}^{\left( 1\right) }\rightarrow _{\left( s,\mathcal{V}\right) }%
\mathbb{E}\mu \left( \cdot \right) ,  \label{stable}
\end{equation}%
where $\mathcal{V}\triangleq \vee _{n}\mathcal{V}_{n}$, and $\mu \in \mathbf{%
M}$ verifies (\ref{randomfou}).
\end{theorem}

\bigskip

\textbf{Remarks }-- (a) Condition (\ref{emboite}) says that $\mathcal{V}_{n}$%
, $n\geq 1$, must be an increasing sequence of $\sigma $-fields, whose $n$th
term is contained in $\mathcal{F}_{n,r_{n}}$, for every $n\geq 1$. Condition
(\ref{CP1}) ensures that, for $i=1,2$, the sum of the first $r_{n}$ terms of
the vector $X_{n}^{\left( i\right) }$ is asymptotically negligeable.

(b)\textbf{\ }There are some differences between the statement of Theorem %
\ref{T : Xue} above, and the original result presented in \cite{Xue}. On the
one hand, in \cite{Xue} the sequence $\left\{ N_{n}:n\geq 1\right\} $ is
such that each $N_{n}$ is a $\mathcal{F}_{n,\cdot }$-stopping time (but we
do not need such a generality). On the other hand, in \cite{Xue} one
considers only the case of the family of $\sigma $-fields $\mathcal{V}%
_{n}^{\ast }=\cap _{j\geq n}\mathcal{F}_{j,r_{n}}$, $n\geq 1$, where $r_{n}$
is non decreasing (note that, due to the monotonicity of $r_{n}$, the $%
\mathcal{V}_{n}^{\ast }$'s satisfy automatically (\ref{emboite})). However,
by inspection of the proof of \cite[Theorem 2.1 and Lemma 2.1]{Xue}, one
sees immediately that all is needed to prove Theorem \ref{T : Xue} is that
the $\mathcal{V}_{n}$'s verify condition (\ref{emboite}). For instance, if $%
r_{n}$ is a general sequence of natural numbers such that $\mathcal{F}%
_{n,r_{n}}\subseteq \mathcal{F}_{n+1,r_{n+1}}$ for each $n\geq 1$, then the
sequence $\mathcal{V}_{n}=\mathcal{F}_{n,r_{n}}$, $n\geq 1$, trivially
satisfies (\ref{emboite}), even if it does not fit Xue's original
assumptions.

(c) The main theorem in the paper by Jakubowski \cite[Theorem 1.1]{Jakubowki}
(which, to our knowledge, is the first systematic account of the POC)
corresponds to the special case $\mathcal{F}_{n,0}=\left\{ \varnothing
,\Omega \right\} $ and $r_{n}=0$, $n\geq 1$. Under such assumptions,
necessarily $\mathcal{V}_{n}=\mathcal{F}_{n,0},$ $S_{n,r_{n}\wedge
N_{n}}^{\left( i\right) }=0$, $i=1,2$, and $\phi \left( \lambda \right) $,
which is $\vee _{n}\mathcal{V}_{n}=\left\{ \varnothing ,\Omega \right\} $ --
measurable, is deterministic for every $\lambda $. In particular, relations (%
\ref{emboite}) and (\ref{CP1}) become immaterial. See also \cite[Theorem
5.8.3]{KW} and \cite[ Theorem 7.1.4]{GineDelaPena} for some detailed
discussions of the POC in this setting.

(d) For the case $r_{n}=0$ and $\mathcal{F}_{n,0}=\mathcal{A}$ ($n\geq 1$),
where $\mathcal{A}$ is not trivial, see also \cite[Section (1.c)]{JKM}.

\bigskip

The next proposition will be used in Section 5.

\bigskip

\begin{proposition}
\label{P : mom Xue}Let the notation of Theorem \ref{T : Xue} prevail,
suppose that the sequence $S_{n,N_{n}}^{\left( 1\right) }$ verifies (\ref%
{fourconv}) for some $\phi \in \widehat{\mathbf{M}}_{0}$, and assume
moreover that there exists a finite random variable $C\left( \omega \right)
>0$ such that, for some $\eta >0$,%
\begin{equation}
\mathbb{E}\left[ \left\vert S_{n,N_{n}}^{\left( 1\right) }\right\vert ^{\eta
}\mid \mathcal{F}_{n,r_{n}}\right] <C\left( \omega \right) \text{, \ }%
\forall n\geq 1\text{,\ \ a.s.-}\mathbb{P}.  \label{uniform bound}
\end{equation}%
Then, there exists a subsequence $\left\{ n\left( k\right) :k\geq 1\right\} $
such that, a.s. - $\mathbb{P}$,
\begin{equation}
\mathbb{E}\left[ \exp \left( i\lambda S_{n\left( k\right) ,N_{n\left(
k\right) }}^{\left( 1\right) }\right) \mid \mathcal{F}_{n\left( k\right)
,r_{n\left( k\right) }}\right] \underset{k\rightarrow +\infty }{\rightarrow }%
\phi \left( \lambda \right)  \label{asfconv}
\end{equation}%
for every real $\lambda $.
\end{proposition}

\begin{proof}
Combining (\ref{fourconv}) and (\ref{uniform bound}), we deduce the
existence of a set $\Omega ^{\ast }$ of probability one, as well as of a
subsequence $n\left( k\right) $, such that, for every $\omega \in \Omega
^{\ast }$, relation (\ref{uniform bound}) is satisfied and (\ref{asfconv})
holds for every rational $\lambda $. We now fix $\omega \in \Omega ^{\ast }$%
, and show that (\ref{asfconv}) holds for all real $\lambda $. Relations (%
\ref{fourconv}) and (\ref{uniform bound}) imply that%
\begin{equation*}
\mathbb{P}_{k}^{\omega }\left[ \cdot \right] =\mathbb{P}\left[ S_{n\left(
k\right) ,N_{n\left( k\right) }}^{\left( 1\right) }\in \cdot \mid \mathcal{F}%
_{n\left( k\right) ,r_{n\left( k\right) }}\right] \left( \omega \right) ,\ \
k\geq 1,
\end{equation*}%
is tight and hence relatively compact: every sequence of $n\left( k\right) $
has a further subsequence $\left\{ n\left( k_{r}\right) :r\geq 1\right\} $
such that $\mathbb{P}_{k_{r}}^{\omega }\left[ \cdot \right] $ is weakly
convergent, so that the corresponding characteristic function converges. In
view of (\ref{asfconv}), such characteristic function must also satisfy the
asymptotic relation
\begin{equation*}
\mathbb{E}\left[ \exp \left( i\lambda S_{n\left( k_{r}\right) ,N_{n\left(
k_{r}\right) }}^{\left( 1\right) }\right) \mid \mathcal{F}_{n\left(
k_{r}\right) ,r_{n\left( k_{r}\right) }}\right] \left( \omega \right)
\underset{r\rightarrow +\infty }{\rightarrow }\phi \left( \lambda \right)
\left( \omega \right)
\end{equation*}%
for every rational $\lambda $, hence for every real $\lambda $, because $%
\phi \left( \lambda \right) \left( \omega \right) $ is continuous in $%
\lambda $.
\end{proof}

\section{General framework for applications of the POC}

We now present a general framework in which the POC techniques discussed in
the previous paragraph can be applied. The main result of this section turns
out to be the key tool to obtain stable convergence results for multiple
stochastic integrals with respect to independently scattered random measures.

\bigskip

Our first goal is to define an It\^{o} type stochastic integral with respect
to a real valued and square integrable stochastic process $X$ (not
necessarily Gaussian) verifying the following three conditions: (i) $X$ is
indexed by the elements $f$ of a real separable Hilbert space $\mathfrak{H}$%
, (ii) $X$ satisfies the isomorphic relation
\begin{equation}
\mathbb{E}\left[ X\left( f\right) X\left( g\right) \right] =\left(
f,g\right) _{\mathfrak{H}}\text{, \ \ }\forall f,g\in \mathfrak{H}\text{,}
\label{basicISO}
\end{equation}%
and (iii) $X$ has independent increments (the notion of \textquotedblleft
increment\textquotedblright ,\ in this context, is defined through
orthogonal projections--see below). We shall then show that the asymptotic
behavior of such integrals can be studied by means of arrays of random
variables, to which the POC applies quite naturally. Note that the elements
of $\mathfrak{H}$ need not be functions -- they may be e.g. distributions on
$\mathbf{%
\mathbb{R}
}^{d}$, $d\geq 1$. Our construction is inspired by the theory initiated by
L. Wu (see \cite{Wu}) and A.S.\ \"{U}st\"{u}nel and M.\ Zakai (see \cite{UZ}%
), concerning Skorohod integrals and filtrations on abstract Wiener spaces.
These author have introduced the notion of time in the context of abstract
Wiener spaces by using resolutions of the identity.

\bigskip

\textbf{Definition D }(see e.g. \cite{Brod}, \cite{Yoshida} and \cite{UZ})%
\textbf{\ -- }Let $\mathfrak{H}$ be a separable real Hilbert space, endowed
with an inner product $\left( \cdot ,\cdot \right) _{\mathfrak{H}}$ ($%
\left\Vert \cdot \right\Vert _{\mathfrak{H}}$ is the corresponding norm). A
(continuous) \textit{resolution of the identity}, is a family $\pi =\left\{
\pi _{t}:t\in \left[ 0,1\right] \right\} $ of orthogonal projections
satisfying:

\begin{description}
\item[(D-i)] $\pi _{0}=0$, and $\pi _{1}=id.;$

\item[(D-ii)] $\forall 0\leq s<t\leq 1$, $\pi _{s}\mathfrak{H}\subseteq \pi
_{t}\mathfrak{H};$

\item[(D-iii)] $\forall 0\leq t_{0}\leq 1$, $\forall h\in \mathfrak{H}$, $%
\lim_{t\rightarrow t_{0}}\left\Vert \pi _{t}h-\pi _{t_{0}}h\right\Vert _{%
\mathfrak{H}}=0.$
\end{description}

A subset $F$ (not necessarily closed, nor linear) of $\mathfrak{H}$ is said
to be $\pi $-\textit{reproducing}, and is denoted $F_{\pi }$, if the linear
span of the set $\left\{ \pi _{t}f:f\in F_{\pi }\text{, }t\in \left[ 0,1%
\right] \right\} $ is dense in $\mathfrak{H}$ (in which case we say that
such a set is \textit{total} in $\mathfrak{H}$). The \textit{rank} of $\pi $
is the smallest of the dimensions of all the closed subspaces generated by
the $\pi $-reproducing subsets of $\mathfrak{H}$. A $\pi $-reproducing
subset $F_{\pi }$ of $\mathfrak{H}$ is called \textit{fully orthogonal }if $%
\left( \pi _{t}f,g\right) _{\mathfrak{H}}=0$ for every $t\in \left[ 0,1%
\right] $ and every $f,g\in F_{\pi }$. The class of all resolutions of the
identity satisfying conditions (\textbf{D-i})--(\textbf{D-iii}) is denoted $%
\mathcal{R}\left( \mathfrak{H}\right) $.

\bigskip

\textbf{Remarks -- }(a)\textbf{\ }Since $\mathfrak{H}$ is separable, for
every resolution of the identity $\pi $ there always exists a countable $\pi
$-reproducing subset of $\mathfrak{H}$.

(b) Let $\pi $ be a resolution of the identity, and note $\overline{\text{%
v.s.}}\left( A\right) $ the closure of the vector space generated by some $%
A\subseteq \mathfrak{H}$. By a standard Gram-Schmidt orthogonalization
procedure, it is easy to prove that for every $\pi $-reproducing subset $%
F_{\pi }$ of $\mathfrak{H}$ such that $\dim \left( \overline{\text{v.s.}}%
\left( F_{\pi }\right) \right) =rank\left( \pi \right) $, there exists a $%
\pi $-reproducing \textit{and fully orthogonal} subset $F_{\pi }^{\prime }$
of $\mathfrak{H}$, such that $\dim \left( \overline{\text{v.s.}}\left(
F_{\pi }^{\prime }\right) \right) =\dim \left( \overline{\text{v.s.}}\left(
F_{\pi }\right) \right) $ (see e.g. \cite[Lemma 23.2]{Brod}, or \cite[p. 27]%
{UZ}).

\bigskip

\textbf{Examples --} The following examples are related to the content of
Section 5 and Section 6.

(a) Take $\mathfrak{H}=L^{2}\left( \left[ 0,1\right] ,dx\right) $, i.e. the
space of square integrable functions on $\left[ 0,1\right] $. Then, a family
of projection operators naturally associated to $\mathfrak{H}$ can be as
follows: for every $t\in \left[ 0,1\right] $ and every $f\in \mathfrak{H}$,
\begin{equation}
\pi _{t}f\left( x\right) =f\left( x\right) \mathbf{1}_{\left[ 0,t\right]
}\left( x\right) .  \label{exPROJ}
\end{equation}%
It is easily seen that this family $\pi =\left\{ \pi _{t}:t\in \left[ 0,1%
\right] \right\} $ is a resolution of the identity verifying conditions
\textbf{(Di)}--\textbf{(Diii)} in Definition D. Also, $rank\left( \pi
\right) =1$, since the linear span of the projections of the function $%
f\left( x\right) \equiv 1$ generates $\mathfrak{H}$.

(b) If $\mathfrak{H}=L^{2}\left( \left[ 0,1\right] ^{2},dxdy\right) $, we
define: for every $t\in \left[ 0,1\right] $ and every $f\in \mathfrak{H}$,
\begin{equation}
\pi _{t}f\left( x,y\right) =f\left( x,y\right) \mathbf{1}_{\left[ 0,t\right]
^{2}}\left( x,y\right) .  \label{exPROJ2}
\end{equation}%
The family $\pi =\left\{ \pi _{t}:t\in \left[ 0,1\right] \right\} $
appearing in (\ref{exPROJ2}) is a resolution of the identity as in
Definition D. However, in this case $rank\left( \pi \right) =+\infty $.
Other choices of $\pi _{t}$ are also possible, for instance
\begin{equation*}
\pi _{t}f\left( x,y\right) =f\left( x,y\right) \mathbf{1}_{\left[ \frac{1}{2}%
-\frac{t}{2},\frac{1}{2}+\frac{t}{2}\right] ^{2}}\left( x,y\right) \text{, }
\end{equation*}%
which expands from the center of the square $\left[ 0,1\right] ^{2}.$

\bigskip

Now fix a real separable Hilbert space $\mathfrak{H}$, as well as a
probability space $\left( \Omega ,\mathcal{F},\mathbb{P}\right) $. In what
follows, we will write
\begin{equation}
X=X\left( \mathfrak{H}\right) =\left\{ X\left( f\right) :f\in \mathfrak{H}%
\right\}  \label{general isometric}
\end{equation}%
to denote a collection of centered random variables defined on $\left(
\Omega ,\mathcal{F},\mathbb{P}\right) $, indexed by the elements of $%
\mathfrak{H}$ and satisfying the isomorphic relation (\ref{basicISO}) (we
use the notation $X\left( \mathfrak{H}\right) $ when the role of the space $%
\mathfrak{H}$ is relevant to the discussion). Note that relation (\ref%
{basicISO}) implies that, for every $f,g\in \mathfrak{H}$, $X\left(
f+g\right) =X\left( f\right) +X\left( g\right) $, a.s.-$\mathbb{P}$.

\bigskip

Let $X\left( \mathfrak{H}\right) $ be defined as in (\ref{general isometric}%
). Then, for every resolution $\pi =\left\{ \pi _{t}:t\in \left[ 0,1\right]
\right\} \in \mathcal{R}\left( \mathfrak{H}\right) $, the following property
is verified: $\forall m\geq 2$, $\forall h_{1},...,h_{m}\in \mathfrak{H}$
and $\forall 0\leq t_{0}<t_{1}<...<t_{m}\leq 1$, the vector
\begin{equation}
\left( X\left( \left( \pi _{t_{1}}-\pi _{t_{0}}\right) h_{1}\right) ,X\left(
\left( \pi _{t_{2}}-\pi _{t_{1}}\right) h_{2}\right) ...,X\left( \left( \pi
_{t_{m}}-\pi _{t_{m-1}}\right) h_{m}\right) \right)  \label{indincrements}
\end{equation}%
is composed of uncorrelated random variables, because the $\pi _{t}$'s are
orthogonal projections. We stress that the class $\mathcal{R}\left(
\mathfrak{H}\right) $ depends only on the Hilbert space $\mathfrak{H}$, and
not on $X$. Now define $\mathcal{R}_{X}\left( \mathfrak{H}\right) $ to be
the subset of $\mathcal{R}\left( \mathfrak{H}\right) $ containing those $\pi
$ such that the vector (\ref{indincrements}) is composed of jointly
independent random variables, for any choice of $m\geq 2$, $%
h_{1},...,h_{m}\in \mathfrak{H}$ and $0\leq t_{0}<t_{1}<...<t_{m}\leq 1$.
The set $\mathcal{R}_{X}\left( \mathfrak{H}\right) $ depends in general of $%
X $. Note that, if $X\left( \mathfrak{H}\right) $ is a Gaussian family, then
$\mathcal{R}_{X}\left( \mathfrak{H}\right) =\mathcal{R}\left( \mathfrak{H}%
\right) $ (see Section 3 below). To every $\pi \in \mathcal{R}_{X}\left(
\mathfrak{H}\right) $ we associate the filtration%
\begin{equation}
\mathcal{F}_{t}^{\pi }\left( X\right) =\sigma \left\{ X\left( \pi
_{t}f\right) :f\in \mathfrak{H}\right\} \text{, \ \ }t\in \left[ 0,1\right]
\text{,}  \label{resfiltration}
\end{equation}%
so that, for instance, $\mathcal{F}_{1}^{\pi }\left( X\right) =\sigma \left(
X\right) .$

\bigskip

\textbf{Remark -- }Note that, for every $h\in \mathfrak{H}$ and every $\pi
\in \mathcal{R}_{X}\left( \mathfrak{H}\right) $, the stochastic process $%
t\mapsto X\left( \pi _{t}h\right) $ is a centered, square integrable $%
\mathcal{F}_{t}^{\pi }\left( X\right) $-martingale with independent
increments. Moreover, since $\pi $ is continuous and (\ref{basicISO}) holds,
$X\left( \pi _{s}h\right) \overset{\mathbb{P}}{\rightarrow }X\left( \pi
_{t}h\right) $ whenever $s\rightarrow t$. In the terminology of \cite[p. 3]%
{Sato}, this implies that $\left\{ X\left( \pi _{t}h\right) :t\in \left[ 0,1%
\right] \right\} $ is an \textit{additive process in law}. In particular, if
$\mathcal{R}_{X}\left( \mathfrak{H}\right) \ $is not empty, for every $h\in
\mathfrak{H}$ the law of $X\left( h\right) $ is infinitely divisible (see
e.g. \cite[Theorem 9.1]{Sato}). As a consequence (see \cite[Theorem 8.1 and
formula (8.8), p. 39]{Sato}), for every $h\in \mathfrak{H}$ there exists a
unique pair $\left( c^{2}\left( h\right) ,\nu _{h}\right) $ such that $%
c^{2}\left( h\right) \in \left[ 0,+\infty \right) $ and $\nu _{h}$ is a
measure on $\mathbb{R}$ satisfying
\begin{equation}
\nu _{h}\left( \left\{ 0\right\} \right) =0\text{, \ \ }\int_{\mathbb{R}%
}\left( x^{2}\wedge 1\right) \nu _{h}\left( dx\right) <+\infty \text{ \ \
and \ \ }\int_{\left\vert x\right\vert \geq 1}x^{2}\nu _{h}\left( dx\right)
<+\infty  \label{LK0}
\end{equation}%
(the last relation follows from the fact that $X\left( h\right) $ is square
integrable (see \cite[Section 5.25]{Sato})), and moreover, for every $%
\lambda \in \mathbb{R}$,
\begin{equation}
\mathbb{E}\left[ \exp \left( i\lambda X\left( h\right) \right) \right] =\exp %
\left[ -\frac{c^{2}\left( h\right) \lambda ^{2}}{2}+\int_{\mathbb{R}}\left(
\exp \left( i\lambda x\right) -1-i\lambda x\right) \nu _{h}\left( dx\right) %
\right] .  \label{L-K}
\end{equation}

\bigskip

Observe that, since the L\'{e}vy-Khintchine representation of an infinitely
divisible distribution is unique, the pair $\left( c^{2}\left( h\right) ,\nu
_{h}\right) $ does not depend on the choice of $\pi \in \mathcal{R}%
_{X}\left( \mathfrak{H}\right) $. In what follows, when $\mathcal{R}%
_{X}\left( \mathfrak{H}\right) \neq \varnothing $, we will use the notation:
for every $\lambda \in \mathbb{R}$ and every $h\in \mathfrak{H}$,
\begin{equation}
\psi _{\mathfrak{H}}\left( h;\lambda \right) \triangleq -\frac{c^{2}\left(
h\right) \lambda ^{2}}{2}+\int_{\mathbb{R}}\left( \exp \left( i\lambda
x\right) -1-i\lambda x\right) \nu _{h}\left( dx\right) ,  \label{exponent}
\end{equation}%
where the pair $\left( c^{2}\left( h\right) ,\nu _{h}\right) $,
characterizing the law of the random variable $X\left( h\right) $, is given
by (\ref{L-K}). Note that, if $h_{n}\rightarrow h$ in $\mathfrak{H}$, then $%
X\left( h_{n}\right) \rightarrow X\left( h\right) $ in $L^{2}\left( \mathbb{P%
}\right) $, and therefore $\psi _{\mathfrak{H}}\left( h_{n};\lambda \right)
\rightarrow \psi _{\mathfrak{H}}\left( h;\lambda \right) $ for every $%
\lambda \in \mathbb{R}$ (uniformly on compacts). We shall always endow $%
\mathfrak{H}$ with the $\sigma $-field $\mathcal{B}\left( \mathfrak{H}%
\right) $, generated by the open sets with respect to the distance induced
by the norm $\left\Vert \cdot \right\Vert _{\mathfrak{H}}$. Since, for every
real $\lambda $, the complex-valued application $h\mapsto \psi _{\mathfrak{H}%
}\left( h;\lambda \right) $ is continuous, it is also $\mathcal{B}\left(
\mathfrak{H}\right) $-measurable.

\bigskip

\textbf{Examples }-- (a) Take $\mathfrak{H}=L^{2}\left( \left[ 0,1\right]
,dx\right) $, suppose that $X\left( \mathfrak{H}\right) =\left\{ X\left(
h\right) :h\in \mathfrak{H}\right\} $ is a centered Gaussian family
verifying (\ref{basicISO}), and define the resolution of the identity $\pi
=\left\{ \pi _{t}:t\in \left[ 0,1\right] \right\} $ according to (\ref%
{exPROJ}). Then, if $\mathfrak{1}$ indicates the function which is
constantly equal to one, the process
\begin{equation}
W_{t}\triangleq X\left( \pi _{t}\mathfrak{1}\right) \text{, \ \ }t\in \left[
0,1\right] \text{,}  \label{BMexo}
\end{equation}%
is a standard Brownian motion started from zero,%
\begin{equation*}
\mathcal{F}_{t}^{\pi }\left( X\right) =\sigma \left\{ W_{s}:s\leq t\right\}
\text{, \ \ }\forall t\in \left[ 0,1\right] \text{,}
\end{equation*}%
and, for every $f\in \mathfrak{H}$,
\begin{equation*}
X\left( \pi _{t}f\right) =\int_{0}^{t}f\left( s\right) dW_{s},
\end{equation*}%
where the stochastic integration is in the usual Wiener-It\^{o} sense. Of
course, $X\left( \pi _{t}f\right) $ is a Gaussian $\mathcal{F}_{t}^{\pi
}\left( X\right) $ - martingale with independent increments, and also, by
using the notation (\ref{exponent}), for every $f\in L^{2}\left( \left[ 0,1%
\right] ,dx\right) $ and $\lambda \in \mathbb{R}$, $\psi _{\mathfrak{H}%
}\left( f;\lambda \right) =-\left( \lambda ^{2}/2\right) \int_{0}^{1}f\left(
x\right) ^{2}dx$.

(b) Take $\mathfrak{H}=L^{2}\left( \left[ 0,1\right] ^{2},dxdy\right) $ and
define the resolution $\pi =\left\{ \pi _{t}:t\in \left[ 0,1\right] \right\}
$ as in (\ref{exPROJ2}). We consider a compensated Poisson measure $\widehat{%
N}=\left\{ \widehat{N}\left( C\right) :C\in \mathcal{B}\left( \left[ 0,1%
\right] ^{2}\right) \right\} $ over $\left[ 0,1\right] ^{2}$. This means
that (1) for every $C\in \mathcal{B}\left( \left[ 0,1\right] ^{2}\right) $,
\begin{equation*}
\widehat{N}\left( C\right) \overset{\text{law}}{=}N\left( C\right) -\mathbb{E%
}\left( N\left( C\right) \right)
\end{equation*}%
where $N\left( C\right) $ is a Poisson random variable with parameter $%
Leb\left( C\right) $ (i.e., the Lebesgue measure of $C$), and (2) $\widehat{N%
}\left( C_{1}\right) $ and $\widehat{N}\left( C_{2}\right) $ are
stochastically independent whenever $C_{1}\cap C_{2}=\varnothing $. Then,
the family $X\left( \mathfrak{H}\right) =\left\{ X\left( h\right) :h\in
\mathfrak{H}\right\} $, defined by
\begin{equation*}
X\left( h\right) =\int_{\left[ 0,1\right] ^{2}}h\left( x,y\right) \widehat{N}%
\left( dx,dy\right) \text{, \ \ }h\in \mathfrak{H}\text{,}
\end{equation*}%
satisfies the isomorphic relation (\ref{basicISO}). Moreover%
\begin{equation*}
\mathcal{F}_{t}^{\pi }\left( X\right) =\sigma \left\{ \widehat{N}\left( %
\left[ 0,s\right] \times \left[ 0,u\right] \right) :s\vee u\leq t\right\}
\text{, \ \ }\forall t\in \left[ 0,1\right] \text{,}
\end{equation*}%
and for every $h\in \mathfrak{H}$, the process%
\begin{equation*}
X\left( \pi _{t}h\right) =\int_{\left[ 0,t\right] ^{2}}h\left( x,y\right)
\widehat{N}\left( dx,dy\right) \text{, \ \ }t\in \left[ 0,1\right] \text{,}
\end{equation*}%
is a $\mathcal{F}_{t}^{\pi }\left( X\right) $ -- martingale with independent
increments, and hence $\pi \in \mathcal{R}_{X}\left( \mathfrak{H}\right) $.
Moreover, for every $h\in L^{2}\left( \left[ 0,1\right] ^{2},dxdy\right) $
and $\lambda \in \mathbb{R}$ the exponent $\psi _{\mathfrak{H}}\left(
h;\lambda \right) $ in (\ref{exponent}) verifies the relation (see e.g. \cite%
[Proposition 19.5]{Sato})
\begin{equation*}
\psi _{\mathfrak{H}}\left( h;\lambda \right) =\int_{0}^{1}\int_{0}^{1}\left[
\exp \left( i\lambda h\left( x,y\right) \right) -1-i\lambda h\left(
x,y\right) \right] dxdy.
\end{equation*}

\bigskip

We now want to consider random variables with values in $\mathfrak{H}$, and
define an It\^{o} type stochastic integral with respect to $X$. To do so, we
let $L^{2}\left( \mathbb{P},\mathfrak{H},X\right) =L^{2}\left( \mathfrak{H}%
,X\right) $ be the space of $\sigma \left( X\right) $-measurable and $%
\mathfrak{H}$-valued random variables $Y$ satisfying $\mathbb{E}\left[
\left\Vert Y\right\Vert _{\mathfrak{H}}^{2}\right] $ $<$ $+\infty $ (note
that $L^{2}\left( \mathfrak{H},X\right) $ is a Hilbert space, with inner
product $\left( Y,Z\right) _{L^{2}\left( \mathfrak{H},X\right) }$ $=\mathbb{E%
}\left[ \left( Y,Z\right) _{\mathfrak{H}}\right] $). Following for instance
\cite{UZ} (which concerns uniquely the Gaussian case), we associate to every
$\pi \in \mathcal{R}_{X}\left( \mathfrak{H}\right) $ the subspace $L_{\pi
}^{2}\left( \mathfrak{H},X\right) $ of the $\pi $-\textit{adapted} elements
of $L^{2}\left( \mathfrak{H},X\right) $, that is: $Y\in L_{\pi }^{2}\left(
\mathfrak{H},X\right) $ if, and only if, $Y\in L^{2}\left( \mathfrak{H}%
,X\right) $ and, for every $t\in \left[ 0,1\right] $ and every $h\in
\mathfrak{H}$,%
\begin{equation}
\left( Y,\pi _{t}h\right) _{\mathfrak{H}}\in \mathcal{F}_{t}^{\pi }\left(
X\right) \text{.}  \label{adaptation}
\end{equation}

For any resolution $\pi \in \mathcal{R}_{X}\left( \mathfrak{H}\right) $, $%
L_{\pi }^{2}\left( \mathfrak{H},X\right) $ is a closed subspace of $%
L^{2}\left( \mathfrak{H},X\right) $. Indeed, if $Y_{n}\in L_{\pi }^{2}\left(
\mathfrak{H},X\right) $ and $Y_{n}\rightarrow Y$ in $L^{2}\left( \mathfrak{H}%
,X\right) $, then necessarily $\left( Y_{n},\pi _{t}h\right) _{\mathfrak{H}}%
\overset{\mathbb{P}}{\rightarrow }\left( Y,\pi _{t}h\right) _{\mathfrak{H}}$
$\forall t\in \left[ 0,1\right] $ and every $h\in \mathfrak{H}$, thus
yielding $Y\in L_{\pi }^{2}\left( \mathfrak{H},X\right) $. We will
occasionally write $\left( u,z\right) _{L_{\pi }^{2}\left( \mathfrak{H}%
\right) }$ instead of $\left( u,z\right) _{L^{2}\left( \mathfrak{H}\right) }$%
, when both $u$ and $z$ are in $L_{\pi }^{2}\left( \mathfrak{H},X\right) $.
Now define, for $\pi \in \mathcal{R}_{X}\left( \mathfrak{H}\right) $, $%
\mathcal{E}_{\pi }\left( \mathfrak{H},X\right) $ to be the space of ($\pi $%
-adapted) \textit{elementary} elements of $L_{\pi }^{2}\left( \mathfrak{H}%
,X\right) $, that is, $\mathcal{E}_{\pi }\left( \mathfrak{H},X\right) $ is
the collection of those elements of $L_{\pi }^{2}\left( \mathfrak{H}%
,X\right) $ that are linear combinations of $\mathfrak{H}$-valued random
variables of the type%
\begin{equation}
h=\Phi \left( t_{1}\right) \left( \pi _{t_{2}}-\pi _{t_{1}}\right) f\text{,}
\label{elementaryad}
\end{equation}%
where $t_{2}>t_{1}$, $f\in \mathfrak{H}$ and $\Phi \left( t_{1}\right) $ is
a random variable which is square-integrable and $\mathcal{F}_{t_{1}}^{\pi
}\left( X\right) $ - measurable.

\begin{lemma}
\label{L : adapt}For every $\pi \in \mathcal{R}_{X}\left( \mathfrak{H}%
\right) $, the set $\mathcal{E}_{\pi }\left( \mathfrak{H},X\right) $, of
adapted elementary elements, is total (i.e., its span is dense) in $L_{\pi
}^{2}\left( \mathfrak{H},X\right) $.
\end{lemma}

\begin{proof}
The proof is similar to \cite[Lemma 2.2]{UZ}. Suppose $u\in L_{\pi
}^{2}\left( \mathfrak{H},X\right) $ and $\left( u,g\right) _{L^{2}\left(
\mathfrak{H},X\right) }=0$ for every $g\in \mathcal{E}_{\pi }\left(
\mathfrak{H},X\right) $. We shall show that $u=0$, a.s. - $\mathbb{P}$. For
every $t_{i+1}>t_{i}$, every bounded and $\mathcal{F}_{t_{i}}^{\pi }\left(
X\right) $-measurable r.v. $\Phi \left( t_{i}\right) $, and every $f\in
\mathfrak{H}$%
\begin{equation*}
\mathbb{E}\left[ \left( \Phi \left( t_{i}\right) \left( \pi _{t_{i+1}}-\pi
_{t_{i}}\right) f,u\right) _{\mathfrak{H}}\right] =0,
\end{equation*}%
and therefore $t\mapsto \left( \pi _{t}f,u\right) _{\mathfrak{H}}$ is a
continuous (since $\pi $ is continuous) $\mathcal{F}_{t}^{\pi }\left(
X\right) $ - martingale starting from zero. Moreover, for every $%
0=t_{0}<\cdots <t_{n}=1$%
\begin{equation*}
\sum_{i=0}^{n-1}\left\vert \left( f,\left( \pi _{t_{i+1}}-\pi
_{t_{i}}\right) u\right) _{\mathfrak{H}}\right\vert \leq \left\Vert
u\right\Vert _{\mathfrak{H}}\left\Vert f\right\Vert _{\mathfrak{H}}\underset{%
\text{a.s.-}\mathbb{P}}{<}+\infty \text{,}
\end{equation*}%
which implies that the continuous martingale $t\mapsto \left( \pi
_{t}f,u\right) _{\mathfrak{H}}$ has also (a.s.-$\mathbb{P}$) bounded
variation. It is therefore constant and hence equal to zero (see e.g. \cite[%
Proposition 1.2]{RY}). It follows that, a.s.-$\mathbb{P}$, $\left(
f,u\right) _{\mathfrak{H}}=\left( \pi _{1}f,u\right) _{\mathfrak{H}}=0$ for
every $f\in \mathfrak{H}$, and consequently $u=0$, a.s.-$\mathbb{P}$.
\end{proof}

\bigskip

We now want to introduce, for every $\pi \in \mathcal{R}_{X}\left( \mathfrak{%
H}\right) $, an It\^{o} type stochastic integral with respect to $X$. To
this end, we fix $\pi \in \mathcal{R}_{X}\left( \mathfrak{H}\right) $ and
first consider simple integrands of the form $h=\sum_{i=1}^{n}\lambda
_{i}h_{i}\in \mathcal{E}_{\pi }\left( \mathfrak{H},X\right) $, where $%
\lambda _{i}\in
\mathbb{R}
$, $n\geq 1$, and $h_{i}$ is as in (\ref{elementaryad}), i.e.%
\begin{equation}
h_{i}=\Phi _{i}\left( t_{1}^{\left( i\right) }\right) \left( \pi
_{t_{2}^{\left( i\right) }}-\pi _{t_{1}^{\left( i\right) }}\right) f_{i}%
\text{, \ \ }f_{i}\in \mathfrak{H}\text{, \ \ }i=1,...,n,
\label{elementaryad2}
\end{equation}%
with $t_{2}^{\left( i\right) }>t_{1}^{\left( i\right) }$, and $\Phi
_{i}\left( t_{1}^{\left( i\right) }\right) \in \mathcal{F}_{t_{1}^{\left(
i\right) }}^{\pi }\left( X\right) $ and square integrable. Then, the
stochastic integral of such a $h$ with respect to $X$ and $\pi $, is defined
as%
\begin{equation}
J_{X}^{\pi }\left( h\right) =\sum_{i=1}^{n}\lambda _{i}J_{X}^{\pi }\left(
h_{i}\right) =\sum_{i=1}^{n}\lambda _{i}\Phi _{i}\left( t_{1}^{\left(
i\right) }\right) X\left( \left( \pi _{t_{2}^{\left( i\right) }}-\pi
_{t_{1}^{\left( i\right) }}\right) f_{i}\right) .  \label{elSkorohod}
\end{equation}

Observe that the $\left( \pi _{t_{2}^{\left( i\right) }}-\pi _{t_{1}^{\left(
i\right) }}\right) f_{i}$ in (\ref{elementaryad2}) becomes the argument of $%
X $ in (\ref{elSkorohod}). Note also that, although $X$ has $\pi $%
-independent increments, there may be a very complex dependence structure
between the random variables
\begin{equation*}
J_{X}^{\pi }\left( h_{i}\right) =\Phi _{i}\left( t_{1}^{\left( i\right)
}\right) X\left( \left( \pi _{t_{2}^{\left( i\right) }}-\pi _{t_{1}^{\left(
i\right) }}\right) f_{i}\right) \text{, \ \ }i=1,...,n,
\end{equation*}%
since the $\Phi _{i}$'s are non-trivial functionals of $X$. We therefore
introduce a \textquotedblleft decoupled\textquotedblright\ version of the
integral $J_{X}^{\pi }\left( h\right) $, by considering an independent copy
of $X$, noted $\widetilde{X}$, and by substituting $X$ with $\widetilde{X}$
in formula (\ref{elSkorohod}). That is, for every $h$ $\in \mathcal{E}_{\pi
}\left( \mathfrak{H},X\right) $ as in (\ref{elementaryad2}) we define%
\begin{equation}
J_{\widetilde{X}}^{\pi }\left( h\right) =\sum_{i=1}^{n}\lambda _{i}\Phi
_{i}\left( t_{1}^{\left( i\right) }\right) \widetilde{X}\left( \left( \pi
_{t_{2}^{\left( i\right) }}-\pi _{t_{1}^{\left( i\right) }}\right)
f_{i}\right) .  \label{decoupSko}
\end{equation}%
\textit{\ }

Note that if $h\in \mathcal{E}_{\pi }\left( \mathfrak{H},X\right) $ is non
random, i.e. $h\left( \omega \right) =h^{\ast }\in \mathfrak{H}$, a.s.-$%
\mathbb{P}\left( d\omega \right) $, then the integrals $J_{X}^{\pi }\left(
h\right) =X\left( h^{\ast }\right) $ and $J_{\widetilde{X}}^{\pi }\left(
h\right) =\widetilde{X}\left( h^{\ast }\right) $ are independent copies of
each other.

\begin{proposition}
\label{P : ItoInt}Fix $\pi \in \mathcal{R}_{X}\left( \mathfrak{H}\right) $.
Then, for every $h,h^{\prime }\in \mathcal{E}_{\pi }\left( \mathfrak{H}%
,X\right) $,%
\begin{eqnarray}
\mathbb{E}\left( J_{X}^{\pi }\left( h\right) J_{X}^{\pi }\left( h^{\prime
}\right) \right) &=&\left( h,h^{\prime }\right) _{L_{\pi }^{2}\left(
\mathfrak{H}\right) }  \label{higherisometyry} \\
\mathbb{E}\left( J_{\widetilde{X}}^{\pi }\left( h\right) J_{\widetilde{X}%
}^{\pi }\left( h^{\prime }\right) \right) &=&\left( h,h^{\prime }\right)
_{L_{\pi }^{2}\left( \mathfrak{H}\right) }.  \notag
\end{eqnarray}%
As a consequence, there exist two linear extensions of $J_{X}^{\pi }$ and $%
J_{\widetilde{X}}^{\pi }$ to $L_{\pi }^{2}\left( \mathfrak{H},X\right) $
satisfying the following two conditions:

\begin{enumerate}
\item if $h_{n}$ converges to $h$ in $L_{\pi }^{2}\left( \mathfrak{H}%
,X\right) $, then
\begin{equation*}
\lim_{n\rightarrow +\infty }\mathbb{E}\left[ \left( J_{X}^{\pi }\left(
h_{n}\right) -J_{X}^{\pi }\left( h\right) \right) ^{2}\right]
=\lim_{n\rightarrow +\infty }\mathbb{E}\left[ \left( J_{\widetilde{X}}^{\pi
}\left( h_{n}\right) -J_{\widetilde{X}}^{\pi }\left( h\right) \right) ^{2}%
\right] =0;
\end{equation*}

\item for every $h,h^{\prime }\in L_{\pi }^{2}\left( \mathfrak{H},X\right) $%
\begin{equation}
\mathbb{E}\left( J_{X}^{\pi }\left( h\right) J_{X}^{\pi }\left( h^{\prime
}\right) \right) =\mathbb{E}\left( J_{\widetilde{X}}^{\pi }\left( h\right)
J_{\widetilde{X}}^{\pi }\left( h^{\prime }\right) \right) =\left(
h,h^{\prime }\right) _{L_{\pi }^{2}\left( \mathfrak{H}\right) }.
\label{Hiso}
\end{equation}
\end{enumerate}

The two extensions $J_{X}^{\pi }$ and $J_{\widetilde{X}}^{\pi }$ are unique,
in the sense that if $\widehat{J}_{X}^{\pi }$ and $\widehat{J}_{\widetilde{X}%
}^{\pi }$ are two other extensions satisfying properties 1 and 2 above, then
necessarily, a.s.-$\mathbb{P}$,%
\begin{equation*}
J_{X}^{\pi }\left( h\right) =\widehat{J}_{X}^{\pi }\left( h\right) \text{ \
\ and \ \ }J_{\widetilde{X}}^{\pi }\left( h\right) =\widehat{J}_{\widetilde{X%
}}^{\pi }\left( h\right)
\end{equation*}%
for every $h\in L_{\pi }^{2}\left( \mathfrak{H},X\right) .$
\end{proposition}

\begin{proof}
It is sufficient to prove (\ref{higherisometyry}) when $h$ and $h^{\prime }$
are simple adapted elements of the kind (\ref{elementaryad}), and in this
case the result follows from elementary computations. Since, according to
Lemma \ref{L : adapt}, $\mathcal{E}_{\pi }\left( \mathfrak{H},X\right) $ is
dense in $L_{\pi }^{2}\left( \mathfrak{H},X\right) $, the result is obtained
from a standard density argument.
\end{proof}

\bigskip

The following property, which is a consequence of the above discussion,
follows immediately.

\begin{corollary}
\label{C : mart}For every $f\in L_{\pi }^{2}\left( \mathfrak{H},X\right) $,
the process
\begin{equation*}
t\mapsto J_{X}^{\pi }\left( \pi _{t}f\right) \text{, \ \ }t\in \left[ 0,1%
\right]
\end{equation*}%
is a real valued $\mathcal{F}_{t}^{\pi }$ - martingale initialized at zero.
\end{corollary}

Observe that the process $t\mapsto J_{X}^{\pi }\left( \pi _{t}f\right) $, $%
t\in \left[ 0,1\right] $, need not have independent (nor conditionally
independent) increments. On the other hand, due to the independence between $%
X$ and $\widetilde{X}$, and to (\ref{indincrements}), \textit{conditionally }%
on the $\sigma $-field $\sigma \left( X\right) $, the increments of the
process $t\mapsto J_{\widetilde{X}}^{\pi }\left( \pi _{t}f\right) $ are
independent (to see this, just consider the process $J_{\widetilde{X}}^{\pi
}\left( \pi _{t}f\right) $ for an elementary $f$ as in (\ref{decoupSko}),
and observe that, in this case, conditioning on $\sigma \left( X\right) $ is
equivalent to conditioning on the $\Phi _{i}$'s; the general case is
obtained once again by a density argument). It follows that the random
process $J_{\widetilde{X}}^{\pi }\left( \pi _{\cdot }f\right) $ can be
regarded as being \textit{decoupled and} \textit{tangent }to $J_{X}^{\pi
}\left( \pi _{\cdot }f\right) $, in a spirit similar to \cite[Definition 4.1]%
{KW91}, \cite{Jacod} or \cite{Jacod84}. We stress, however, that $J_{%
\widetilde{X}}^{\pi }\left( \pi _{\cdot }f\right) $ need not meet the
definition of a tangent process given in such references, which is based on
a notion of convergence in the Skorohod topology, rather than on the $L^{2}$%
-convergence adopted in the present paper. The reader is referred to \cite%
{Jacod} for an exhaustive characterization of processes with conditionally
independent increments.

Now, for $h\in \mathfrak{H}$ and $\lambda \in \mathbb{R}$, define the
exponent $\psi _{\mathfrak{H}}\left( h;\lambda \right) $ according to (\ref%
{exponent}), and observe that every $f\in L_{\pi }^{2}\left( \mathfrak{H}%
,X\right) $ is a random element with values in $\mathfrak{H}$.\ It follows
that the quantity $\psi _{\mathfrak{H}}\left( f\left( \omega \right)
;\lambda \right) $ is well defined for every $\omega \in \Omega $ and every $%
\lambda \in \mathbb{R}$, and moreover, since $\psi _{\mathfrak{H}}\left(
\cdot ;\lambda \right) $ is $\mathcal{B}\left( \mathfrak{H}\right) $%
-measurable, for every $f\in L_{\pi }^{2}\left( \mathfrak{H},X\right) $ and
every $\lambda \in \mathbb{R}$, the complex-valued application $\omega
\mapsto \psi _{\mathfrak{H}}\left( f\left( \omega \right) ;\lambda \right) $%
\ is $\mathcal{F}$-measurable.

\bigskip

\begin{proposition}
For every $\lambda \in \mathbb{R}$ and every $f\in L_{\pi }^{2}\left(
\mathfrak{H},X\right) $,
\begin{equation}
\mathbb{E}\left[ \exp \left( i\lambda J_{\widetilde{X}}^{\pi }\left(
f\right) \right) \mid \sigma \left( X\right) \right] =\exp \left[ \psi _{%
\mathfrak{H}}\left( f;\lambda \right) \right] \text{, \ \ a.s.-}\mathbb{P}%
\text{.}  \label{expo2}
\end{equation}
\end{proposition}

\begin{proof}
For $f\in \mathcal{E}_{\pi }\left( \mathfrak{H},X\right) $, formula (\ref%
{expo2}) follows immediately from the independence of $X$ and $\widetilde{X}$%
. Now fix $f\in L_{\pi }^{2}\left( \mathfrak{H},X\right) $, and select a
sequence $\left( f_{n}\right) \subset \mathcal{E}_{\pi }\left( \mathfrak{H}%
,X\right) $ such that
\begin{equation}
\mathbb{E}\left[ \left\Vert f_{n}-f\right\Vert _{\mathfrak{H}}^{2}\right]
\rightarrow 0  \label{p}
\end{equation}%
(such a sequence $f_{n}$ always exists, due to Lemma \ref{L : adapt}). Since
(\ref{p}) implies that $\left\Vert f_{n}-f\right\Vert _{\mathfrak{H}}\overset%
{\mathbb{P}}{\rightarrow }0$, for every subsequence $n_{k}$ there exists a
further subsequence $n_{k\left( r\right) }$ such that $\left\Vert
f_{n_{k\left( r\right) }}-f\right\Vert _{\mathfrak{H}}\rightarrow 0$, a.s. -
$\mathbb{P}$, thus implying $\psi _{\mathfrak{H}}\left( f_{n_{k\left(
r\right) }};\lambda \right) \rightarrow \psi _{\mathfrak{H}}\left( f;\lambda
\right) $ for every $\lambda \in \mathbb{R}$, a.s. - $\mathbb{P}$. Then, for
every $\lambda \in \mathbb{R}$, $\psi _{\mathfrak{H}}\left( f_{n};\lambda
\right) \overset{\mathbb{P}}{\rightarrow }\psi _{\mathfrak{H}}\left(
f;\lambda \right) $, and therefore $\exp \left[ \psi _{\mathfrak{H}}\left(
f_{n};\lambda \right) \right] \overset{\mathbb{P}}{\rightarrow }\exp \left[
\psi _{\mathfrak{H}}\left( f;\lambda \right) \right] $. On the other hand,%
\begin{eqnarray*}
\mathbb{E}\left\vert \mathbb{E}\left[ \exp \left( i\lambda J_{\widetilde{X}%
}^{\pi }\left( f_{n}\right) \right) -\exp \left( i\lambda J_{\widetilde{X}%
}^{\pi }\left( f\right) \right) \mid \sigma \left( X\right) \right]
\right\vert &\leq &\left\vert \lambda \right\vert \mathbb{E}\left\vert J_{%
\widetilde{X}}^{\pi }\left( f_{n}\right) -J_{\widetilde{X}}^{\pi }\left(
f\right) \right\vert \\
&\leq &\left\vert \lambda \right\vert \mathbb{E}\left[ \left( J_{\widetilde{X%
}}^{\pi }\left( f_{n}\right) -J_{\widetilde{X}}^{\pi }\left( f\right)
\right) ^{2}\right] ^{\frac{1}{2}} \\
&=&\left\vert \lambda \right\vert \mathbb{E}\left[ \left\Vert
f_{n}-f\right\Vert _{\mathfrak{H}}^{2}\right] ^{\frac{1}{2}}\rightarrow 0%
\text{,}
\end{eqnarray*}%
where the equality follows from (\ref{Hiso}), thus yielding
\begin{equation*}
\exp \left[ \psi _{\mathfrak{H}}\left( f_{n};\lambda \right) \right] =%
\mathbb{E}\left[ \exp \left( i\lambda J_{\widetilde{X}}^{\pi }\left(
f_{n}\right) \right) \mid \sigma \left( X\right) \right] \overset{\mathbb{P}}%
{\rightarrow }\mathbb{E}\left[ \exp \left( i\lambda J_{\widetilde{X}}^{\pi
}\left( f\right) \right) \mid \sigma \left( X\right) \right] ,
\end{equation*}%
and the desired conclusion is therefore obtained.
\end{proof}

\bigskip

\textbf{Examples -- }(a) Take $\mathfrak{H}=L^{2}\left( \left[ 0,1\right]
,dx\right) $ and suppose that $X\left( \mathfrak{H}\right) =\left\{ X\left(
h\right) :h\in \mathfrak{H}\right\} $ is a centered Gaussian family
verifying (\ref{basicISO}). Define also $\pi =\left\{ \pi _{t}:t\in \left[
0,1\right] \right\} \in \mathcal{R}\left( \mathfrak{H}\right) $ according to
(\ref{exPROJ}), and write $W$ to denote the Brownian motion introduced in (%
\ref{BMexo}). The subsequent discussion will make clear that $L_{\pi
}^{2}\left( \mathfrak{H},X\right) $ is, in this case, the space of square
integrable processes that are adapted to the Brownian filtration $\sigma
\left\{ W_{u}:u\leq t\right\} $, $t\in \left[ 0,1\right] $. Moreover, for
every $t\in \left[ 0,1\right] $ and $u\in L_{\pi }^{2}\left( \mathfrak{H}%
,X\right) $%
\begin{equation*}
J_{X}^{\pi }\left( \pi _{t}u\right) =\int_{0}^{t}u\left( s\right) dW_{s}%
\text{ \ \ and \ \ }J_{\widetilde{X}}^{\pi }\left( \pi _{t}u\right)
=\int_{0}^{t}u\left( s\right) d\widetilde{W}_{s}\text{,}
\end{equation*}%
where the stochastic integration is in the It\^{o} sense, and $\widetilde{W}%
_{t}\triangleq \widetilde{X}\left( \mathbf{1}_{\left[ 0,t\right] }\right) $
is a standard Brownian motion independent of $X$.

(b) (\textit{Orthogonalized Teugels martingales}, see \cite{NuSch}) Let $%
Z=\left\{ Z_{t}:t\in \left[ 0,1\right] \right\} $ be a real-valued and
centered L\'{e}vy process, initialized at zero and endowed with a L\'{e}vy
measure $\nu $ satisfying the condition: for some $\varepsilon ,\lambda >0$%
\begin{equation*}
\int_{\left( -\varepsilon ,\varepsilon \right) ^{c}}\exp \left( \lambda
\left\vert x\right\vert \right) \nu \left( dx\right) <+\infty .
\end{equation*}%
Then, for every $i\geq 2$, $\int_{\mathbb{R}}\left\vert x\right\vert ^{i}\nu
\left( dx\right) <+\infty $, and $Z_{t}$ has moments of all orders. Starting
from $Z$, for every $i\geq 1$ one can therefore define the \textit{%
compensated power jump process} (or \textit{Teugel martingale}) of order $i$%
, noted $Y^{\left( i\right) }$, as $Y_{t}^{\left( 1\right) }=Z_{t}$ for $%
t\in \left[ 0,1\right] $,\ and, for $i\geq 2$ and $t\in \left[ 0,1\right] $,%
\begin{equation*}
Y_{t}^{\left( i\right) }=\sum_{0<s\leq t}\left( \Delta Z_{t}\right) ^{i}-%
\mathbb{E}\sum_{0<s\leq t}\left( \Delta Z_{t}\right) ^{i}=\sum_{0<s\leq
t}\left( \Delta Z_{t}\right) ^{i}-t\int_{\mathbb{R}}x^{i}\nu \left(
dx\right) .
\end{equation*}%
Plainly, each $Y^{\left( i\right) }$ is a centered L\'{e}vy process.
Moreover, according to \cite[pp. 111-112]{NuSch}, for every $i\geq 1$ it is
possible to find (unique) real coefficients $a_{i,1},...,a_{i,i}$, such that
$a_{i,i}=1$ and the stochastic processes%
\begin{equation*}
H_{t}^{\left( i\right) }=Y_{t}^{\left( i\right) }+a_{i,i-1}Y_{t}^{\left(
i-1\right) }+\cdot \cdot \cdot +a_{i,1}Y_{t}^{\left( 1\right) },\text{ \ \ }%
t\in \left[ 0,1\right] \text{, \ \ }i\geq 1,
\end{equation*}%
are strongly orthogonal centered martingales (in the sense of \cite[p.148]%
{Protter}), also verifying $\mathbb{E}\left[ H_{t}^{\left( i\right)
}H_{s}^{\left( j\right) }\right] $ $=$ $\delta _{ij}\left( t\wedge s\right) $%
, where $\delta _{ij}$ is the Kronecker symbol. Observe that $H^{\left(
i\right) }$ is again a L\'{e}vy process, and that, for every deterministic $%
g,f\in L^{2}\left( \left[ 0,1\right] ,ds\right) $, the integrals $%
\int_{0}^{1}f\left( s\right) dH_{s}^{\left( i\right) }$ and $%
\int_{0}^{1}g\left( s\right) dH_{s}^{\left( j\right) }$ are well defined and
such that
\begin{equation}
\mathbb{E}\left[ \int_{0}^{1}f\left( s\right) dH_{s}^{\left( i\right)
}\int_{0}^{1}g\left( s\right) dH_{s}^{\left( j\right) }\right] =\delta
_{ij}\int_{0}^{1}g\left( s\right) f\left( s\right) ds.  \label{izo}
\end{equation}%
Now define $\mathfrak{H}=L^{2}\left( \mathbb{N}\times \left[ 0,1\right]
,\kappa \left( dm\right) \times ds\right) $, where $\kappa \left( dm\right) $
is the counting measure, and define, for $h\left( \cdot ,\cdot \right) \in
\mathfrak{H}$, $t\in \left[ 0,1\right] $, and $\left( m,s\right) \in \mathbb{%
N}\times \left[ 0,1\right] $,%
\begin{equation*}
\pi _{t}h\left( m,s\right) =h\left( m,s\right) \mathbf{1}_{\left[ 0,t\right]
}\left( s\right) \text{.}
\end{equation*}%
It is clear that $\pi =\left\{ \pi _{t}:t\in \left[ 0,1\right] \right\} \in
\mathcal{R}\left( \mathfrak{H}\right) $. Moreover, for every $h\left( \cdot
,\cdot \right) \in \mathfrak{H}$, we define%
\begin{equation*}
X\left( h\right) =\sum_{m=1}^{\infty }\int_{0}^{1}h\left( m,s\right)
dH_{s}^{\left( m\right) }\text{,}
\end{equation*}%
where the series is convergent in $L^{2}\left( \mathbb{P}\right) $, since $%
\mathbb{E}X\left( h\right) ^{2}=\sum \int_{0}^{1}h\left( m,s\right)
^{2}ds<+\infty $, due to (\ref{izo}) and the fact that $h\in \mathfrak{H}$.
Since the $H^{\left( m\right) }$ are strongly orthogonal and (\ref{izo})
holds, one sees immediately that, for every $h,h^{\prime }\in \mathfrak{H}$,
$\mathbb{E}\left[ X\left( h\right) X\left( h^{\prime }\right) \right]
=\left( h,h^{\prime }\right) _{\mathfrak{H}}$, and moreover, since for every
$m$ and every $h$ the process $t\mapsto \int_{0}^{1}\pi _{t}h\left(
m,s\right) dH_{s}^{\left( m\right) }=\int_{0}^{t}h\left( m,s\right)
dH_{s}^{\left( m\right) }$ has independent increments, $\pi \in \mathcal{R}%
_{X}\left( \mathfrak{H}\right) $. We can also consider random $h$, and, by
using \cite{NuSch}, give the following characterization of random variables $%
h\in L_{\pi }^{2}\left( \mathfrak{H},X\right) $, and the corresponding
integrals $J_{X}^{\pi }\left( h\right) $ and $J_{\widetilde{X}}^{\pi }\left(
h\right) $: (i) for every $h\in L_{\pi }^{2}\left( \mathfrak{H},X\right) $
there exists a family $\left\{ \phi _{m,t}^{\left( h\right) }:t\in \left[ 0,1%
\right] \text{, }m\geq 1\right\} $ of real-valued and $\mathcal{F}_{t}^{\pi
} $-predictable processes such that for every fixed $m$, the process $%
t\mapsto \phi _{m,t}^{\left( h\right) }$ is a modification of $t\mapsto
h\left( m,t\right) $; (ii) for every $h\in L_{\pi }^{2}\left( \mathfrak{H}%
,X\right) $,%
\begin{equation}
J_{X}^{\pi }\left( h\right) =\sum_{m=1}^{\infty }\int_{0}^{1}\phi
_{m,t}^{\left( h\right) }dH_{t}^{\left( m\right) }\text{,}  \label{NS1}
\end{equation}%
where the series is convergent in $L^{2}\left( \mathbb{P}\right) $; (iii)
for every $h\in L_{\pi }^{2}\left( \mathfrak{H},X\right) $,
\begin{equation}
J_{\widetilde{X}}^{\pi }\left( h\right) =\sum_{m=1}^{\infty
}\int_{0}^{1}\phi _{m,t}^{\left( h\right) }d\widetilde{H}_{t}^{\left(
m\right) },  \label{NS2}
\end{equation}%
where the series is convergent in $L^{2}\left( \mathbb{P}\right) $, and the
sequence $\left\{ \widetilde{H}^{\left( m\right) }:m\geq 1\right\} $ is an
independent copy of $\left\{ H^{\left( m\right) }:m\geq 1\right\} $. Note
that by using \cite[Theorem 1]{NuSch}, one would obtain an analogous
characterization in terms of iterated stochastic integrals of deterministic
kernels.

\section{Stable convergence}

We shall now apply Theorem \ref{T : Xue} to the setup outlined in the
previous paragraph. Let $\mathfrak{H}_{n}$, $n\geq 1$, be a sequence of real
separable Hilbert spaces, and, for each $n\geq 1$, let%
\begin{equation}
X_{n}=X_{n}\left( \mathfrak{H}_{n}\right) =\left\{ X_{n}\left( g\right)
:g\in \mathfrak{H}_{n}\right\} \text{,}  \label{IsoSequence}
\end{equation}%
be a centered, real-valued stochastic process, indexed by the elements of $%
\mathfrak{H}_{n}$ and such that $\mathbb{E}\left[ X_{n}\left( f\right)
X_{n}\left( g\right) \right] $ $=\left( f,g\right) _{\mathfrak{H}_{n}}$. The
processes $X_{n}$ are not necessarily Gaussian. As before, $\widetilde{X}%
_{n} $ indicates an independent copy of $X_{n}$, for every $n\geq 1$.

\begin{theorem}
\label{T : Main}Let the previous notation prevail, and suppose that the
processes $X_{n}$, $n\geq 1$, appearing in (\ref{IsoSequence}) (along with
the independent copies $\widetilde{X}_{n}$) are all defined on the same
probability space $\left( \Omega ,\mathcal{F},\mathbb{P}\right) $. For every
$n\geq 1$, let $\pi ^{\left( n\right) }\in \mathcal{R}_{X_{n}}\left(
\mathfrak{H}_{n}\right) $ and $u_{n}\in L_{\pi ^{\left( n\right)
}}^{2}\left( \mathfrak{H}_{n},X_{n}\right) $. Suppose also that there exists
a sequence $\left\{ t_{n}:n\geq 1\right\} $ $\subset \left[ 0,1\right] $ and
a collection of $\sigma $-fields $\left\{ \mathcal{U}_{n}:n\geq 1\right\} $,
such that
\begin{equation*}
\lim_{n\rightarrow +\infty }\mathbb{E}\left[ \left\Vert \pi _{t_{n}}^{\left(
n\right) }u_{n}\right\Vert _{\mathfrak{H}_{n}}^{2}\right] =0
\end{equation*}%
and
\begin{equation}
\mathcal{U}_{n}\subseteq \mathcal{U}_{n+1}\cap \mathcal{F}_{t_{n}}^{\pi
^{\left( n\right) }}\left( X_{n}\right) .  \label{boxed}
\end{equation}%
If%
\begin{equation}
\exp \left[ \psi _{\mathfrak{H}_{n}}\left( u_{n};\lambda \right) \right] =%
\mathbb{E}\left[ \exp \left( i\lambda J_{\widetilde{X}_{n}}^{\pi ^{\left(
n\right) }}\left( u_{n}\right) \right) \mid \sigma \left( X_{n}\right) %
\right] \overset{\mathbb{P}}{\rightarrow }\phi \left( \lambda \right) =\phi
\left( \lambda ,\omega \right) ,\text{ \ \ }\forall \lambda \in
\mathbb{R}
\text{,}  \label{CT7}
\end{equation}%
where $\psi _{\mathfrak{H}_{n}}\left( u_{n};\lambda \right) $ is defined
according to (\ref{exponent}), $\phi \in \widehat{\mathbf{M}}_{0}$ and, $%
\forall \lambda \in
\mathbb{R}
$,
\begin{equation*}
\phi \left( \lambda \right) \in \vee _{n}\mathcal{U}_{n}\triangleq \mathcal{U%
}^{\ast }\text{,}
\end{equation*}%
then, as $n\rightarrow +\infty $,
\begin{equation}
\mathbb{E}\left[ \exp \left( i\lambda J_{X_{n}}^{\pi ^{\left( n\right)
}}\left( u_{n}\right) \right) \mid \mathcal{F}_{t_{n}}^{\pi ^{\left(
n\right) }}\left( X_{n}\right) \right] \overset{\mathbb{P}}{\rightarrow }%
\phi \left( \lambda \right) \text{, \ \ }\forall \lambda \in
\mathbb{R}
\text{,}  \label{CV+}
\end{equation}%
and%
\begin{equation}
J_{X_{n}}^{\pi ^{\left( n\right) }}\left( u_{n}\right) \rightarrow _{\left(
s,\mathcal{U}^{\ast }\right) }\mathbb{E}\mu \left( \cdot \right) ,
\label{CV++}
\end{equation}%
where $\mu \in \mathbf{M}$ verifies (\ref{randomfou}).
\end{theorem}

\bigskip

\textbf{Remarks -- }(1) The proof of Theorem \ref{T : Main} uses Theorem \ref%
{T : Xue}, which assumes $\phi \in \widehat{\mathbf{M}}_{0}$, that is, $\phi
$ is non-vanishing. If $\phi \in \widehat{\mathbf{M}}$ (instead of $\widehat{%
\mathbf{M}}_{0}$) and if, for example, there exists a subsequence $n_{k}$
such that,
\begin{equation*}
\mathbb{P}\left\{ \omega :\exp \left[ \psi _{\mathfrak{H}_{n_{k}}}\left(
u_{n_{k}}\left( \omega \right) ;\lambda \right) \right] \rightarrow \phi
\left( \lambda ,\omega \right) \text{, \ \ }\forall \lambda \in \mathbb{R}%
\right\} =1\text{,}
\end{equation*}%
then, given the nature of $\psi _{\mathfrak{H}_{n_{k}}}$, $\phi \left(
\lambda ,\omega \right) $ is necessarily, for $\mathbb{P}$-a.e. $\omega $,
the Fourier transform of an infinitely divisible distribution (see e.g. \cite%
[Lemma 7.5]{Sato}), and therefore $\phi \in \widehat{\mathbf{M}}_{0}$. A
similar remark applies to Theorem \ref{T : Main 2} below.

(2)\textbf{\ }For $n\geq 1$, the process $t\mapsto J_{X_{n}}^{\pi ^{\left(
n\right) }}\left( \pi _{t}^{\left( n\right) }u_{n}\right) $ is a martingale
and hence admits a c\`{a}dl\`{a}g modification. Then, an alternative
approach to obtain results for stable convergence is to use the well-known
criteria for the stable convergence of continuous-time c\`{a}dl\`{a}g
semimartingales, as stated e.g. in \cite[Proposition 1 and Theorems 1 and 2 ]%
{Feigin} or \cite[Chapter 4]{JacSh}. However, the formulation in terms of
\textquotedblleft principle of conditioning\textquotedblright\ yields, in
our setting, more precise results, by using less stringent assumptions. For
instance, (\ref{boxed}) can be regarded as a weak version of the
\textquotedblleft nesting condition\textquotedblright\ used in \cite[p. 126 ]%
{Feigin}, whereas (\ref{CV+}) is a refinement of the conclusions that can be
obtained by means of \cite[Proposition 1]{Feigin}.

(3) Suppose that, under the assumptions of Theorem \ref{T : Main}, there
exists a c\`{a}dl\`{a}g process $Y=\left\{ Y_{t}:t\in \left[ 0,1\right]
\right\} $ such that, conditionally on $\mathcal{U}^{\ast }$, $Y$ has
independent increments and $\phi \left( \lambda \right) =\mathbb{E}\left[
\exp \left( i\lambda Y_{1}\right) \mid \mathcal{U}^{\ast }\right] $. In this
case, formula (\ref{CV++}) is equivalent to saying that $J_{X_{n}}^{\pi
^{\left( n\right) }}\left( u_{n}\right) $ converges $\mathcal{U}^{\ast }$%
-stably to $Y_{1}$. See \cite[Section 4]{Jacod} for several results
concerning the stable convergence (for instance, in the sense of finite
dimensional distributions) of semimartingales towards processes with
conditionally independent increments.

\bigskip

Before proving Theorem \ref{T : Main}, we consider the important case of a
\textit{nested }sequence of resolutions. More precisely, assume that $%
\mathfrak{H}_{n}=\mathfrak{H}$, $X_{n}=X$, for every $n\geq 1$, and that the
sequence $\pi ^{\left( n\right) }\in \mathcal{R}_{X}\left( \mathfrak{H}%
\right) $, $n\geq 1$, is nested in the following sense: for every $t\in %
\left[ 0,1\right] $ and every $n\geq 1$,
\begin{equation}
\pi _{t}^{\left( n\right) }\mathfrak{H}\subseteq \pi _{t}^{\left( n+1\right)
}\mathfrak{H}  \label{nestedres}
\end{equation}%
(note that if $\pi ^{\left( n\right) }=\pi $ for every $n$, then (\ref%
{nestedres}) is trivially satisfied); in this case, if $t_{n}$ is non
decreasing, the sequence $\mathcal{U}_{n}=\mathcal{F}_{t_{n}}^{\pi ^{\left(
n\right) }}\left( X\right) $, $n\geq 1$, automatically satisfies (\ref{boxed}%
). We therefore have the following consequence of Theorem \ref{T : Main}.

\begin{corollary}
\label{C : corNest}Under the above notation and assumptions, suppose that
the sequence $\pi ^{\left( n\right) }\in \mathcal{R}_{X}\left( \mathfrak{H}%
\right) $, $n\geq 1$, is nested in the sense of (\ref{nestedres}), and let $%
u_{n}\in L_{\pi ^{\left( n\right) }}^{2}\left( \mathfrak{H},X\right) $, $%
n\geq 1$. Suppose also that there exists a non-decreasing sequence $\left\{
t_{n}:n\geq 1\right\} \subset \left[ 0,1\right] $ s.t.%
\begin{equation}
\lim_{n\rightarrow +\infty }\mathbb{E}\left[ \left\Vert \pi _{t_{n}}^{\left(
n\right) }u_{n}\right\Vert _{\mathfrak{H}}^{2}\right] =0.  \label{QuadNeg}
\end{equation}%
If%
\begin{equation*}
\exp \left[ \psi _{\mathfrak{H}}\left( u_{n};\lambda \right) \right] \overset%
{\mathbb{P}}{\rightarrow }\phi \left( \lambda \right) =\phi \left( \lambda
,\omega \right) ,\text{ \ \ }\forall \lambda \in
\mathbb{R}
\text{,}
\end{equation*}%
where $\phi \in \widehat{\mathbf{M}}_{0}$ and, $\forall \lambda \in
\mathbb{R}
$, $\phi \left( \lambda \right) \in \vee _{n}\mathcal{F}_{t_{n}}^{\pi
^{\left( n\right) }}\left( X\right) \triangleq \mathcal{F}_{\ast }$, then,
as $n\rightarrow +\infty $,
\begin{equation*}
\mathbb{E}\left[ \exp \left( i\lambda J_{X}\left( u_{n}\right) \right) \mid
\mathcal{F}_{t_{n}}^{\pi ^{\left( n\right) }}\left( X\right) \right] \overset%
{\mathbb{P}}{\rightarrow }\phi \left( \lambda \right) ,\text{ \ \ }\forall
\lambda \in
\mathbb{R}
\text{,}
\end{equation*}%
and%
\begin{equation*}
J_{X}\left( u_{n}\right) \rightarrow _{\left( s,\mathcal{F}_{\ast }\right) }%
\mathbb{E}\mu \left( \cdot \right) ,
\end{equation*}%
where $\mu \in \mathbf{M}$ verifies (\ref{randomfou}).
\end{corollary}

In the next result $\left\{ u_{n}\right\} $ may still be random, but $\phi
\left( \lambda \right) $ is non-random. It follows from Corollary \ref{C :
corNest} by taking $t_{n}=0$ for every $n$, so that (\ref{QuadNeg}) is
immaterial, and $\mathcal{F}_{\ast }$ becomes the trivial $\sigma $-field.

\begin{corollary}
\label{C : corDet}Keep the notation of Corollary \ref{C : corNest}, and
consider a (not necessarily nested) sequence $\pi ^{\left( n\right) }\in
\mathcal{R}_{X}\left( \mathfrak{H}\right) $, $n\geq 1$. If%
\begin{equation*}
\exp \left[ \psi _{\mathfrak{H}}\left( u_{n};\lambda \right) \right] \overset%
{\mathbb{P}}{\rightarrow }\phi \left( \lambda \right) \text{, \ \ }\forall
\lambda \in
\mathbb{R}
\text{,}
\end{equation*}%
where $\phi $ is the Fourier transform of some non-random measure $\mu $
such that $\phi \left( \lambda \right) \neq 0$ for every $\lambda \in
\mathbb{R}
$, then, as $n\rightarrow +\infty $,
\begin{equation*}
\mathbb{E}\left[ \exp \left( i\lambda J_{X}^{\pi ^{\left( n\right) }}\left(
u_{n}\right) \right) \right] \rightarrow \phi \left( \lambda \right) ,\text{
\ \ }\forall \lambda \in
\mathbb{R}
\text{,}
\end{equation*}%
that is, the law of $J_{X}^{\pi ^{\left( n\right) }}\left( u_{n}\right) $
converges weakly to $\mu $.
\end{corollary}

\bigskip

\textbf{Proof of Theorem \ref{T : Main} -- }Since $u_{n}\in L_{\pi ^{\left(
n\right) }}^{2}\left( \mathfrak{H}_{n},X_{n}\right) $, there exists, thanks
to Lemma \ref{L : adapt} a sequence $u_{n}^{e}\in \mathcal{E}_{\pi ^{\left(
n\right) }}\left( \mathfrak{H}_{n},X_{n}\right) $, $n\geq 1$, such that (by
using the isometry properties of $J_{\widetilde{X}_{n}}^{\pi ^{\left(
n\right) }}$ and $J_{X_{n}}^{\pi ^{\left( n\right) }}$, as stated in
Proposition \ref{P : ItoInt})%
\begin{eqnarray}
0 &=&\lim_{n\rightarrow +\infty }\mathbb{E}\left[ \left\Vert
u_{n}-u_{n}^{e}\right\Vert _{\mathfrak{H}_{n}}^{2}\right] =\lim_{n%
\rightarrow +\infty }\mathbb{E}\left[ \left( J_{\widetilde{X}_{n}}^{\pi
^{\left( n\right) }}\left( u_{n}\right) -J_{\widetilde{X}_{n}}^{\pi ^{\left(
n\right) }}\left( u_{n}^{e}\right) \right) ^{2}\right]  \label{zero} \\
&=&\lim_{n\rightarrow +\infty }\mathbb{E}\left[ \left( J_{X_{n}}^{\pi
^{\left( n\right) }}\left( u_{n}\right) -J_{X_{n}}^{\pi ^{\left( n\right)
}}\left( u_{n}^{e}\right) \right) ^{2}\right]  \notag
\end{eqnarray}%
and
\begin{eqnarray}
0 &=&\lim_{n\rightarrow +\infty }\mathbb{E}\left[ \left\Vert \pi
_{t_{n}}^{\left( n\right) }u_{n}^{e}\right\Vert _{\mathfrak{H}_{n}}^{2}%
\right] =\lim_{n\rightarrow +\infty }\mathbb{E}\left[ \left( J_{\widetilde{X}%
_{n}}^{\pi ^{\left( n\right) }}\left( \pi _{t_{n}}^{\left( n\right)
}u_{n}^{e}\right) \right) ^{2}\right]  \label{zero2} \\
&=&\lim_{n\rightarrow +\infty }\mathbb{E}\left[ \left( J_{X_{n}}^{\pi
^{\left( n\right) }}\left( \pi _{t_{n}}^{\left( n\right) }u_{n}^{e}\right)
\right) ^{2}\right] .  \notag
\end{eqnarray}%
Without loss of generality, we can always suppose that $u_{n}^{e}$ has the
form
\begin{equation*}
u_{n}^{e}=\sum_{i=1}^{N_{n}}\left[ \sum_{j=1}^{M_{n}\left( i\right) }\Phi
_{j}^{\left( n\right) }\left( t_{i-1}^{\left( n\right) }\right) \left( \pi
_{t_{i}^{\left( n\right) }}^{\left( n\right) }-\pi _{t_{i-1}^{\left(
n\right) }}^{\left( n\right) }\right) f_{j}^{\left( n\right) }\right]
\end{equation*}%
where $0=t_{0}^{\left( n\right) }<...<t_{N_{n}}^{\left( n\right) }=1$, $%
f_{j}^{\left( n\right) }\in \mathfrak{H}_{n}$, $N_{n},M_{n}\left( i\right)
\geq 1$, $\Phi _{j}^{\left( n\right) }\left( t_{i-1}^{\left( n\right)
}\right) $ is square integrable and measurable with respect to $\mathcal{F}%
_{t_{i-1}^{\left( n\right) }}^{\pi ^{\left( n\right) }}\left( X_{n}\right) $
where one of the $t_{0}^{\left( n\right) },...,t_{N_{n}}^{\left( n\right) }$
equals $t_{n}$. Moreover, we have
\begin{eqnarray*}
J_{X_{n}}^{\pi ^{\left( n\right) }}\left( u_{n}^{e}\right)
&=&\sum_{i=1}^{N_{n}}\left[ \sum_{j=1}^{M_{n}\left( i\right) }\Phi
_{j}^{\left( n\right) }\left( t_{i-1}^{\left( n\right) }\right) X_{n}\left(
(\pi _{t_{i}^{\left( n\right) }}-\pi _{t_{i-1}^{\left( n\right)
}})f_{j}^{\left( n\right) }\right) \right] \\
J_{\widetilde{X}_{n}}^{\pi ^{\left( n\right) }}\left( u_{n}^{e}\right)
&=&\sum_{i=1}^{N_{n}}\left[ \sum_{j=1}^{M_{n}\left( i\right) }\Phi
_{j}^{\left( n\right) }\left( t_{i-1}^{\left( n\right) }\right) \widetilde{X}%
_{n}\left( (\pi _{t_{i}^{\left( n\right) }}-\pi _{t_{i-1}^{\left( n\right)
}})f_{j}^{\left( n\right) }\right) \right] .
\end{eqnarray*}%
Now define for $n\geq 1$ and $i=1,...,N_{n}$%
\begin{eqnarray*}
X_{n,i}^{\left( 1\right) } &=&\sum_{j=1}^{M_{n}\left( i\right) }\Phi
_{j}^{\left( n\right) }\left( t_{i-1}^{\left( n\right) }\right) X_{n}\left(
(\pi _{t_{i}^{\left( n\right) }}-\pi _{t_{i-1}^{\left( n\right)
}})f_{j}^{\left( n\right) }\right) \\
X_{n,i}^{\left( 2\right) } &=&\sum_{j=1}^{M_{n}\left( i\right) }\Phi
_{j}^{\left( n\right) }\left( t_{i-1}^{\left( n\right) }\right) \widetilde{X}%
_{n}\left( (\pi _{t_{i}^{\left( n\right) }}-\pi _{t_{i-1}^{\left( n\right)
}})f_{j}^{\left( n\right) }\right)
\end{eqnarray*}%
as well as $X_{n,0}^{\left( \ell \right) }=0$, $\ell =1,2$; introduce
moreover the filtration
\begin{equation}
\widehat{\mathcal{F}}_{t}^{\left( \pi ^{\left( n\right) },\mathfrak{H}%
_{n}\right) }=\mathcal{F}_{t}^{\pi ^{\left( n\right) }}\left( X_{n}\right)
\vee \sigma \left\{ \widetilde{X}\left( \pi _{t}^{\left( n\right) }f\right)
:f\in \mathfrak{H}_{n}\right\} \text{, \ }t\in \left[ 0,1\right] \text{,}
\label{2fil}
\end{equation}%
and let $\mathcal{G}_{n}=\sigma \left( X_{n}\right) $, $n\geq 1$. We shall
verify that the array $X^{\left( 2\right) }=\left\{ X_{n,i}^{\left( 2\right)
}:0\leq i\leq N_{n},\text{ }n\geq 1\right\} $ is decoupled and tangent to $%
X^{\left( 1\right) }=\left\{ X_{n,i}^{\left( 1\right) }:0\leq i\leq N_{n},%
\text{ }n\geq 1\right\} $, in the sense of Definition C of Section 2.
Indeed, for $\ell =1,2$, the sequence $\left\{ X_{n,i}^{\left( \ell \right)
}:0\leq i\leq N_{n}\right\} $ is adapted to the discrete filtration
\begin{equation}
\mathcal{F}_{n,i}\triangleq \widehat{\mathcal{F}}_{t_{i}^{\left( n\right)
}}^{\left( \pi ^{\left( n\right) },\mathfrak{H}_{n}\right) }\text{, \ \ }%
i=1,...,N_{n};  \label{DefSigmaF}
\end{equation}%
also (\ref{samedis1}) is satisfied, since, for every $j$ and every $%
i=1,...,N_{n},$%
\begin{equation*}
\Phi _{j}^{\left( n\right) }\left( t_{i-1}^{\left( n\right) }\right) \in
\mathcal{F}_{t_{i-1}^{\left( n\right) }}^{\pi ^{\left( n\right) }}\left(
X_{n}\right) \subset \mathcal{F}_{n,i-1},
\end{equation*}%
and
\begin{eqnarray*}
\mathbb{E}\left[ \exp \left( i\lambda X_{n}\left( (\pi _{t_{i}^{\left(
n\right) }}-\pi _{t_{i-1}^{\left( n\right) }})f_{j}^{\left( n\right)
}\right) \right) \mid \mathcal{F}_{n,i-1}\right] &=&\mathbb{E}\left[ \exp
\left( i\lambda X_{n}\left( (\pi _{t_{i}^{\left( n\right) }}-\pi
_{t_{i-1}^{\left( n\right) }})f_{j}^{\left( n\right) }\right) \right) \right]
\\
&=&\mathbb{E}\left[ \exp \left( i\lambda \widetilde{X}_{n}\left( (\pi
_{t_{i}^{\left( n\right) }}-\pi _{t_{i-1}^{\left( n\right) }})f_{j}^{\left(
n\right) }\right) \right) \mid \mathcal{F}_{n,i-1}\right] .
\end{eqnarray*}%
Since $\mathcal{G}_{n}=\sigma \left( X_{n}\right) $, we obtain immediately\ (%
\ref{samedis2}), because $\widetilde{X}_{n}$ is an independent copy of $%
X_{n} $. We now want to apply Theorem \ref{T : Xue} with%
\begin{eqnarray}
J_{X_{n}}^{\pi }\left( \pi _{t_{n}}^{\left( n\right) }u_{n}^{e}\right)
&=&\sum_{i=1}^{r_{n}}\left[ \sum_{l=1}^{M_{n}\left( i\right) }\Phi
_{l}^{\left( n\right) }\left( t_{i-1}^{\left( n\right) }\right) X_{n}\left(
(\pi _{t_{i}^{\left( n\right) }}-\pi _{t_{i-1}^{\left( n\right)
}})f_{l}^{\left( n\right) }\right) \right] =\sum_{i=1}^{r_{n}}X_{n,i}^{%
\left( 1\right) }=S_{n,r_{n}}^{\left( 1\right) }  \label{identification} \\
J_{\widetilde{X}_{n}}^{\pi }\left( \pi _{t_{n}}u_{n}^{e}\right)
&=&\sum_{i=1}^{r_{n}}\left[ \sum_{l=1}^{M_{n}\left( i\right) }\Phi
_{l}^{\left( n\right) }\left( t_{i-1}^{\left( n\right) }\right) \widetilde{X}%
_{n}\left( (\pi _{t_{i}^{\left( n\right) }}-\pi _{t_{i-1}^{\left( n\right)
}})f_{l}^{\left( n\right) }\right) \right] =\sum_{i=1}^{r_{n}}X_{n,i}^{%
\left( 1\right) }=S_{n,r_{n}}^{\left( 2\right) },  \notag
\end{eqnarray}%
where $r_{n}$ is the element of $\left\{ 1,...,N_{n}\right\} $ such that $%
t_{r_{n}}^{\left( n\right) }=t_{n}$. To do so, we need to verify the
remaining conditions of that theorem. To prove (\ref{emboite}), use (\ref%
{2fil}), (\ref{DefSigmaF}) and (\ref{boxed}), to obtain%
\begin{equation*}
\mathcal{F}_{n,r_{n}}=\widehat{\mathcal{F}}_{t_{r_{n}}^{\left( n\right)
}}^{\left( \pi ^{\left( n\right) },\mathfrak{H}_{n}\right) }\supset \mathcal{%
F}_{t_{n}}^{\pi ^{\left( n\right) }}\left( X_{n}\right) \supseteq \mathcal{U}%
_{n}\text{,}
\end{equation*}%
and hence (\ref{emboite}) holds with $\mathcal{V}_{n}=\mathcal{U}_{n}$. To
prove (\ref{CP1}), observe that the asymptotic relation in (\ref{zero2}) can
be rewritten as
\begin{equation}
\lim_{n\rightarrow +\infty }\mathbb{E}\left[ \left( S_{n,r_{n}}^{\left( \ell
\right) }\right) ^{2}\right] =0\text{, \ \ }\ell =1,2\text{,}
\label{quadrZero}
\end{equation}%
which immediately yields, as $n\rightarrow +\infty $,%
\begin{equation*}
S_{n,r_{n}}^{\left( 1\right) }\overset{\mathbb{P}}{\rightarrow }0\text{ \ \
and \ }\mathbb{E}\left[ \exp \left( i\lambda S_{n,r_{n}}^{\left( 2\right)
}\right) \mid \mathcal{G}_{n}\right] \overset{\mathbb{P}}{\rightarrow }1
\end{equation*}%
for every $\lambda \in
\mathbb{R}
.$ To justify the last relation, just observe that (\ref{quadrZero}) implies
that $\mathbb{E}\left[ \left( S_{n,r_{n}}^{\left( 2\right) }\right) ^{2}\mid
\mathcal{G}_{n}\right] \rightarrow 0$ in $L^{1}\left( \mathbb{P}\right) $,
and hence, for every diverging sequence $n_{k}$, there exists a subsequence $%
n_{k}^{\prime }$ such that, a.s.- $\mathbb{P}$,%
\begin{equation*}
\mathbb{E}\left[ \left( S_{n_{k}^{\prime },r_{n_{k}^{\prime }}}^{\left(
2\right) }\right) ^{2}\mid \mathcal{G}_{n_{k}^{\prime }}\right] \underset{%
k\rightarrow +\infty }{\rightarrow }0\text{,}
\end{equation*}%
which in turn yields that, a.s.-$\mathbb{P}$,
\begin{equation*}
\mathbb{E}\left[ \exp \left( i\lambda S_{n_{k}^{\prime },r_{n_{k}^{\prime
}}}^{\left( 2\right) }\right) \mid \mathcal{G}_{n_{k}^{\prime }}\right]
\underset{k\rightarrow +\infty }{\rightarrow }1\text{.}
\end{equation*}%
To prove (\ref{CP2}), observe that%
\begin{equation*}
\mathbb{E}\left\vert \exp \left( i\lambda J_{\widetilde{X}_{n}}^{\pi
^{\left( n\right) }}\left( u_{n}^{e}\right) \right) -\exp \left( i\lambda J_{%
\widetilde{X}_{n}}^{\pi ^{\left( n\right) }}\left( u_{n}\right) \right)
\right\vert \leq \left\vert \lambda \right\vert \mathbb{E}\left\vert J_{%
\widetilde{X}_{n}}^{\pi ^{\left( n\right) }}\left( u_{n}^{e}\right) -J_{%
\widetilde{X}_{n}}^{\pi ^{\left( n\right) }}\left( u_{n}\right) \right\vert
\underset{n\rightarrow +\infty }{\rightarrow }0,
\end{equation*}%
by (\ref{zero}). Hence, since (\ref{CT7}) holds for $u_{n}$, it also holds
when $u_{n}$ is replaced by the elementary sequence $u_{n}^{e}$. Since $J_{%
\widetilde{X}_{n}}^{\pi ^{\left( n\right) }}\left( u_{n}^{e}\right) =J_{%
\widetilde{X}_{n}}^{\pi ^{\left( n\right) }}\left( \pi _{1}^{\left( n\right)
}u_{n}^{e}\right) =S_{n,N_{n}}^{\left( 2\right) }$ and $\mathcal{G}%
_{n}=\sigma \left( X_{n}\right) $, relation (\ref{CP2}) holds. It follows
that the assumptions of Theorem \ref{T : Xue} are satisfied, and we deduce
that necessarily, as $n\rightarrow +\infty $,%
\begin{eqnarray*}
&&\mathbb{E}\left[ \exp \left( i\lambda J_{X_{n}}^{\pi ^{\left( n\right)
}}\left( u_{n}^{e}\right) \right) \mid \mathcal{F}_{t_{n}}^{\pi ^{\left(
n\right) }}\left( X_{n}\right) \right] \\
&=&\mathbb{E}\left[ \exp \left( i\lambda J_{X_{n}}^{\pi ^{\left( n\right)
}}\left( u_{n}^{e}\right) \right) \mid \widehat{\mathcal{F}}_{t_{n}}^{\left(
\pi ^{\left( n\right) },\mathfrak{H}_{n}\right) }\right] \overset{\mathbb{P}}%
{\rightarrow }\phi \left( \lambda \right) ,\text{ \ \ }\forall \lambda \in
\mathbb{R}
\text{,}
\end{eqnarray*}%
(the equality follows from the fact that $X_{n}$ and $\widetilde{X}_{n}$ are
independent). Theorem \ref{T : Xue} also yields%
\begin{equation}
J_{X_{n}}^{\pi ^{\left( n\right) }}\left( u_{n}^{e}\right) \rightarrow
_{\left( s,\mathcal{U}^{\ast }\right) }\mathbb{E}\mu \left( \cdot \right) .
\label{r1}
\end{equation}%
To go back from $u_{n}^{e}$ to $u_{n}$, we use%
\begin{equation}
\mathbb{E}\left\vert \exp \left( i\lambda J_{X_{n}}^{\pi ^{\left( n\right)
}}\left( u_{n}^{e}\right) \right) -\exp \left( i\lambda J_{X_{n}}^{\pi
^{\left( n\right) }}\left( u_{n}\right) \right) \right\vert \underset{%
n\rightarrow +\infty }{\rightarrow }0,  \label{r0}
\end{equation}%
which follows again from (\ref{zero}), and we deduce that
\begin{equation*}
\mathbb{E}\left[ \exp \left( i\lambda J_{X_{n}}^{\pi ^{\left( n\right)
}}\left( u_{n}^{e}\right) \right) -\exp \left( i\lambda J_{X_{n}}^{\pi
^{\left( n\right) }}\left( u_{n}\right) \right) \mid \mathcal{F}%
_{t_{n}}^{\pi ^{\left( n\right) }}\left( X_{n}\right) \right] \overset{L^{1}}%
{\underset{n\rightarrow +\infty }{\rightarrow }}0\text{,}
\end{equation*}%
and therefore%
\begin{equation*}
\mathbb{E}\left[ \exp \left( i\lambda J_{X_{n}}^{\pi ^{\left( n\right)
}}\left( u_{n}\right) \right) \mid \mathcal{F}_{t_{n}}^{\pi ^{\left(
n\right) }}\left( X_{n}\right) \right] \overset{\mathbb{P}}{\underset{%
n\rightarrow +\infty }{\rightarrow }}\phi \left( \lambda \right) ,\text{ \ \
}\forall \lambda \in
\mathbb{R}
\text{.}
\end{equation*}%
Finally, by combining (\ref{r1}) and (\ref{r0}), we obtain%
\begin{equation*}
J_{X_{n}}^{\pi ^{\left( n\right) }}\left( u_{n}\right) \rightarrow _{\left(
s,\mathcal{U}^{\ast }\right) }\mathbb{E}\mu \left( \cdot \right) .
\end{equation*}%
$\blacksquare $

\bigskip

By using the same approximation procedure as in the preceding proof, we may
use Proposition \ref{P : mom Xue} to prove the following refinement of
Theorem \ref{T : Main}.

\bigskip

\begin{proposition}
\label{P : mom Main}With the notation of Theorem \ref{T : Main}, suppose
that the sequence $J_{X_{n}}^{\pi ^{\left( n\right) }}\left( u_{n}\right) $
verifies (\ref{CV+}), and that there exists a finite random variable $%
C\left( \omega \right) >0$ such that, for some $\eta >0$,%
\begin{equation*}
\mathbb{E}\left[ \left\vert J_{X_{n}}^{\pi ^{\left( n\right) }}\left(
u_{n}\right) \right\vert ^{\eta }\mid \mathcal{F}_{t_{n}}^{\pi ^{\left(
n\right) }}\right] <C\left( \omega \right) \text{, \ }\forall n\geq 1\text{%
,\ \ a.s.-}\mathbb{P}.
\end{equation*}%
Then, there is a subsequence $\left\{ n\left( k\right) :k\geq 1\right\} $
such that, a.s. - $\mathbb{P}$,
\begin{equation*}
\mathbb{E}\left[ \exp \left( i\lambda J_{X_{n}}^{\pi ^{\left( n\right)
}}\left( u_{n}\right) \right) \mid \mathcal{F}_{t_{n\left( k\right) }}^{\pi
^{\left( n\left( k\right) \right) }}\right] \underset{k\rightarrow +\infty }{%
\rightarrow }\phi \left( \lambda \right) \text{, \ \ }\forall \lambda \in
\mathbb{R}.
\end{equation*}
\end{proposition}

\bigskip

Theorem \ref{T : Main} can also be extended to a slightly more general
framework. To this end, we introduce some further notation. Fix a closed
subspace $\mathfrak{H}^{\ast }\subseteq \mathfrak{H}$. For every $t\in \left[
0,1\right] $, we denote by $\pi _{s\leq t}\mathfrak{H}^{\ast }$ the closed
linear subspace of $\mathfrak{H}$, generated by the set $\left\{ \pi
_{s}f:f\in \mathfrak{H}^{\ast }\text{,\ }s\leq t\right\} $. Of course, $\pi
_{\leq t}\mathfrak{H}^{\ast }\subseteq \pi _{t}\mathfrak{H}=\pi _{\leq t}%
\mathfrak{H}$. For a fixed $\pi \in \mathcal{R}_{X}\left( \mathfrak{H}%
\right) $, we set $\mathcal{E}_{\pi }\left( \mathfrak{H},\mathfrak{H}^{\ast
},X\right) $ to be the subset of $\mathcal{E}_{\pi }\left( \mathfrak{H}%
,X\right) $ composed of $\mathfrak{H}$-valued random variables of the kind%
\begin{equation}
h=\Psi ^{\ast }\left( t_{1}\right) \left( \pi _{t_{2}}-\pi _{t_{1}}\right) g%
\text{,}  \label{SpElAd}
\end{equation}%
where $t_{2}>t_{1}$, $g\in \mathfrak{H}^{\ast }$ and $\Psi ^{\ast }\left(
t_{1}\right) $ is a square integrable random variable verifying the
measurability condition%
\begin{equation*}
\Psi ^{\ast }\left( t_{1}\right) \in \sigma \left\{ X\left( f\right) :f\in
\pi _{\leq t_{1}}\mathfrak{H}^{\ast }\right\} \text{,}
\end{equation*}%
whereas $L_{\pi }^{2}\left( \mathfrak{H},\mathfrak{H}^{\ast },X\right) $ is
defined as the closure of $\mathcal{E}_{\pi }\left( \mathfrak{H},\mathfrak{H}%
^{\ast },X\right) $ in $L_{\pi }^{2}\left( \mathfrak{H},X\right) $. Note
that, plainly, $\mathcal{E}_{\pi }\left( \mathfrak{H},X\right) $ $=$ $%
\mathcal{E}_{\pi }\left( \mathfrak{H},\mathfrak{H},X\right) $ and $L_{\pi
}^{2}\left( \mathfrak{H},X\right) =L_{\pi }^{2}\left( \mathfrak{H},\mathfrak{%
H},X\right) $. Moreover, for every $Y\in L_{\pi }^{2}\left( \mathfrak{H},%
\mathfrak{H}^{\ast },X\right) $ and every $t\in \left[ 0,1\right] $, the
following two poperties are verified: (i) the random element $\pi _{t}Y$
takes values in $\pi _{\leq t}\mathfrak{H}^{\ast }$, a.s.-$\mathbb{P}$, and
(ii) the random variable $J_{X}^{\pi }\left( \pi _{t}h\right) $ is
measurable with respect to the $\sigma $-field $\sigma \left\{ X\left(
f\right) :f\in \pi _{\leq t}\mathfrak{H}^{\ast }\right\} $ (such claims are
easily verified for $h$ as in (\ref{SpElAd}), and the general results follow
once again by standard density arguments).

\bigskip

\textbf{Remark -- }Note that, in general, even when $rank\left( \pi \right)
=1$ as in (\ref{exPROJ}), and $\mathfrak{H}^{\ast }$ is non-trivial, for $%
0<t\leq 1$ the set $\pi _{\leq t}\mathfrak{H}^{\ast }$ may be strictly
contained in $\pi _{t}\mathfrak{H}$. It follows that the $\sigma $-field $%
\sigma \left\{ X\left( f\right) :f\in \pi _{\leq t}\mathfrak{H}^{\ast
}\right\} $ can be strictly contained in $\mathcal{F}_{t}^{\pi }\left(
X\right) $, as defined in (\ref{resfiltration}). To see this, just consider
the case $\mathfrak{H}=L^{2}\left( \left[ 0,1\right] ,dx\right) $, $%
\mathfrak{H}^{\ast }=\left\{ f\in L^{2}\left( \left[ 0,1\right] ,dx\right)
:f=f\mathbf{1}_{\left[ 0,1/2\right] }\right\} $, $\pi _{s}f=f\mathbf{1}_{%
\left[ 0,s\right] }$ ($s\in \left[ 0,1\right] $), and take $t\in \left( 1/2,1%
\right] $. Indeed, in this case $X\left( \mathbf{1}_{\left[ 0,t\right]
}\right) $ is $\mathcal{F}_{t}^{\pi }\left( X\right) $-measurable but is not
$\sigma \left\{ X\left( f\right) :f\in \pi _{\leq t_{1}}\mathfrak{H}^{\ast
}\right\} $-measurable.

\bigskip

The following result can be proved along the lines of Lemma \ref{L : adapt}.

\bigskip

\begin{lemma}
For every closed subspace $\mathfrak{H}^{\ast }$ of $\mathfrak{H}$, a random
element $Y\ $is in $L_{\pi }^{2}\left( \mathfrak{H},\mathfrak{H}^{\ast
},X\right) $ if, and only if, $Y\in L^{2}\left( \mathfrak{H},X\right) $ and,
for every $t\in \left[ 0,1\right] $,%
\begin{equation*}
\left( Y,\pi _{t}h\right) _{\mathfrak{H}}\in \sigma \left\{ X\left( f\right)
:f\in \pi _{\leq t}\mathfrak{H}^{\ast }\right\} \text{.}
\end{equation*}
\end{lemma}

\bigskip

The next theorem can be proved by using arguments analogous to the ones in
the proof of Theorem \ref{T : Main}. Here, $\mathfrak{H}_{n}=\mathfrak{H}$
and $X_{n}\left( \mathfrak{H}_{n}\right) =X\left( \mathfrak{H}\right) $ for
every $n$.

\begin{theorem}
\label{T : Main 2}Under the above notation and assumptions, for every $n\geq
1$ let $\mathfrak{H}^{\left( n\right) }$ be a closed subspace of $\mathfrak{H%
}$, $\pi ^{\left( n\right) }\in \mathcal{R}_{X}\left( \mathfrak{H}\right) $,
and $u_{n}\in L_{\pi ^{\left( n\right) }}^{2}\left( \mathfrak{H},\mathfrak{H}%
^{\left( n\right) },X\right) $. Suppose also that there exists a sequence $%
\left\{ t_{n}:n\geq 1\right\} $ $\subset \left[ 0,1\right] $ and a
collection of closed subspaces of $\mathfrak{H}$, noted $\left\{ \mathfrak{U}%
_{n}:n\geq 1\right\} $, such that
\begin{equation*}
\lim_{n\rightarrow +\infty }\mathbb{E}\left[ \left\Vert \pi _{t_{n}}^{\left(
n\right) }u_{n}\right\Vert _{\mathfrak{H}}^{2}\right] =0
\end{equation*}%
and
\begin{equation*}
\mathfrak{U}_{n}\subseteq \mathfrak{U}_{n+1}\cap \pi _{\leq t_{n}}^{\left(
n\right) }\mathfrak{H}^{\left( n\right) }.
\end{equation*}%
If%
\begin{equation*}
\exp \left[ \psi _{\mathfrak{H}}\left( u_{n};\lambda \right) \right] \overset%
{\mathbb{P}}{\rightarrow }\phi \left( \lambda \right) =\phi \left( \lambda
,\omega \right) ,\text{ \ \ }\forall \lambda \in
\mathbb{R}
\text{,}
\end{equation*}%
where $\phi \in \widehat{\mathbf{M}}_{0}$ and, $\forall \lambda \in
\mathbb{R}
$,
\begin{equation*}
\phi \left( \lambda \right) \in \vee _{n}\sigma \left\{ X\left( f\right)
:f\in \mathfrak{U}_{n}\right\} \triangleq \mathcal{U}^{\ast }\text{,}
\end{equation*}%
then, as $n\rightarrow +\infty $,
\begin{equation*}
\mathbb{E}\left[ \exp \left( i\lambda J_{X}^{\pi ^{\left( n\right) }}\left(
u_{n}\right) \right) \mid X\left( f\right) :f\in \pi _{\leq t_{n}}^{\left(
n\right) }\mathfrak{H}^{\left( n\right) }\right] \overset{\mathbb{P}}{%
\rightarrow }\phi \left( \lambda \right) \text{, \ \ }\forall \lambda \in
\mathbb{R}
\text{,}
\end{equation*}%
and%
\begin{equation*}
J_{X}^{\pi ^{\left( n\right) }}\left( u_{n}\right) \rightarrow _{\left( s,%
\mathcal{U}^{\ast }\right) }\mathbb{E}\mu \left( \cdot \right) ,
\end{equation*}%
where $\mu \in \mathbf{M}$ verifies (\ref{randomfou}).
\end{theorem}

\section{Stable limit theorems for multiple integrals with respect
independently scattered measures}

This section concerns multiple integrals with respect to independently
scattered random measures (not necessarily Gaussian) and corresponding limit
theorems. In particular, we will use Theorem \ref{T : Main} to obtain new
central and non-central limit theorems for these multiple integrals,
extending part of the results proved in \cite{NuPe} and \cite{PT} in the
framework of multiple Wiener-It\^{o} integrals with respect to Gaussian
processes. A specific application is described in Section 3.3, where we deal
with sequences of double integrals with respect to Poisson random measures.
For further applications of the theory developed in Section 2 to the
asymptotic analysis of Gaussian fields, the reader is referred to Section 6,
as well as to the companion paper \cite{PeTaq}. For a general discussion
concerning multiple integrals with respect to random measures, see \cite%
{Engel} and \cite{RoWa}. For limit theorems involving multiple stochastic
integrals (and other related classes of random variables), see the two
surveys by Surgailis \cite{Sur} and \cite{SUr2}, and the references therein.

\subsection{Independently scattered random measures and multiple integrals}

From now on $\left( Z,\mathcal{Z},\mu \right) $ stands for a standard Borel
space, with $\mu $ a positive, non-atomic and $\sigma $-finite measure on $%
\left( Z,\mathcal{Z}\right) $. We denote by $\mathcal{Z}_{\mu }$ the subset
of $\mathcal{Z}$ composed of sets of finite $\mu $-measure. Observe that the
$\sigma $-finiteness of $\mu $ implies that $\mathcal{Z}=\sigma \left(
\mathcal{Z}_{\mu }\right) $.

\bigskip

\textbf{Definition E }-- An \textit{independently scattered random measure }$%
M$\textit{\ }on $\left( Z,\mathcal{Z}\right) $, with \textit{control measure
}$\mu $, is a collection of random variables%
\begin{equation*}
M=\left\{ M\left( B\right) :B\in \mathcal{Z}_{\mu }\right\} ,
\end{equation*}%
indexed by the elements of $\mathcal{Z}_{\mu }$ and such that: (\textbf{E-i}%
) for every $B\in \mathcal{Z}_{\mu }$ $M\left( B\right) \in L^{2}\left(
\mathbb{P}\right) $, (\textbf{E-ii}) for every finite collection of disjoint
sets $B_{1},...,B_{m}\in \mathcal{Z}_{\mu }$, the vector $\left( M\left(
B_{1}\right) ,...,M\left( B_{d}\right) \right) $ si composed of mutually
independent random variables; (\textbf{E-iii}) for every $B,C\in \mathcal{Z}%
_{\mu }$,%
\begin{equation}
\mathbb{E}\left[ M\left( B\right) M\left( C\right) \right] =\mu \left( C\cap
B\right) .  \label{intersection}
\end{equation}

\bigskip

Let $\mathfrak{H}_{\mu }=L^{2}\left( Z,\mathcal{Z},\mu \right) $ be the
Hilbert space of real-valued and square-integrable functions on $\left( Z,%
\mathcal{Z}\right) $ (with respect to $\mu $). Since relation (\ref%
{intersection}) holds, it is easily seen that there exists a unique
collection of centered and square-integrable random variables
\begin{equation}
X_{M}=X_{M}\left( \mathfrak{H}_{\mu }\right) =\left\{ X_{M}\left( h\right)
:h\in \mathfrak{H}_{\mu }\right\} \text{,}  \label{isomeasure}
\end{equation}%
such that the following two properties are verified: (a) for every
elementary function $h\in \mathfrak{H}_{\mu }$ with the form $h\left(
z\right) =\sum_{i=1,...,n}c_{i}\mathbf{1}_{B_{i}}\left( z\right) $, where $%
n=1,2,...$, $c_{i}\in \mathbb{R}$ and $B_{i}\in \mathcal{Z}_{\mu }$ are
disjoint, $X_{M}\left( h\right) =\sum_{i=1,...,n}c_{i}M\left( B_{i}\right) $%
, and (b) for every $h,h^{\prime }\in \mathfrak{H}_{\mu }$%
\begin{equation}
\mathbb{E}\left[ X_{M}\left( h\right) X_{M}\left( h^{\prime }\right) \right]
=\left( h,h^{\prime }\right) _{\mathfrak{H}_{\mu }}\triangleq
\int_{Z}h\left( z\right) h^{\prime }\left( z\right) \mu \left( dz\right) .
\label{CVX}
\end{equation}

Property (a) implies in particular that, $\forall B\in \mathcal{Z}_{\mu }$, $%
M\left( B\right) =X_{M}\left( \mathbf{1}_{B}\right) $. Note that $X_{M}$ is
a collection of random variables of the kind defined in formula (\ref%
{general isometric}) of Section 3. Moreover, for every $h\in \mathfrak{H}%
_{\mu }$, the random variable $X_{M}\left( h\right) $ has an infinitely
divisible law. It follows that, for every $h\in \mathfrak{H}_{\mu }$, there
exists a unique pair $\left( c^{2}\left( h\right) ,\nu _{h}\right) $ such
that $c^{2}\left( h\right) \in \left[ 0,+\infty \right) $ and $\nu _{h}$ is
a (L\'{e}vy) measure on $\mathbb{R}$ satisfying the three properties in (\ref%
{LK0}), so that, for every $\lambda \in \mathbb{R}$,
\begin{equation}
\mathbb{E}\left[ \exp \left( i\lambda X_{M}\left( h\right) \right) \right]
=\exp \left[ \psi _{\mathfrak{H}_{\mu }}\left( h;\lambda \right) \right] ,
\label{LK3}
\end{equation}%
where the L\'{e}vy-Khinchine exponent $\psi _{\mathfrak{H}_{\mu }}\left(
h;\lambda \right) $ is defined by (\ref{exponent}).

\bigskip

We now give a characterization of $\psi _{\mathfrak{H}_{\mu }}\left(
h;\lambda \right) $, based on the techniques developed in \cite{Raj Ros}
(but see also \cite[Section 5]{KW91}).

\bigskip

\begin{proposition}
\label{P : LK}For every $B\in \mathcal{Z}_{\mu }$, let $\left( c^{2}\left(
B\right) ,\nu _{B}\right) $ denote the pair such that $c^{2}\left( B\right)
\in \left[ 0,+\infty \right) $, $\nu _{B}$ verifies (\ref{LK0}) and%
\begin{equation}
\psi _{\mathfrak{H}_{\mu }}\left( \mathbf{1}_{B};\lambda \right) =-\frac{%
\lambda ^{2}}{2}c^{2}\left( B\right) +\int_{\mathbb{R}}\left( \exp \left(
i\lambda x\right) -1-i\lambda x\right) \nu _{B}\left( dx\right) .
\label{LK4}
\end{equation}%
Then,

\begin{enumerate}
\item The application $B\mapsto c^{2}\left( B\right) $, from $\mathcal{Z}%
_{\mu }$ to $\left[ 0,+\infty \right) $, extends to a unique $\sigma $%
-finite measure $c^{2}\left( dz\right) $ on $\left( Z,\mathcal{Z}\right) $,
such that $c^{2}\left( dz\right) \ll \mu \left( dz\right) .$

\item There exists a unique measure $\nu $ on $\left( Z\times \mathbb{R},%
\mathcal{Z}\times \mathcal{B}\left( \mathbb{R}\right) \right) $ such that $%
\nu \left( B\times C\right) =\nu _{B}\left( C\right) $, for every $B\in
\mathcal{Z}_{\mu }$ and $C\in \mathcal{B}\left( \mathbb{R}\right) $.

\item There exists a function $\rho _{\mu }:Z\times \mathcal{B}\left(
\mathbb{R}\right) \mapsto \left[ 0,+\infty \right] $ such that (i) for every
$z\in Z$, $\rho _{\mu }\left( z,\cdot \right) $ is a L\'{e}vy measure%
\footnote{%
That is, $\rho _{\mu }\left( z,\left\{ 0\right\} \right) =0$ and $\int_{%
\mathbb{R}}\min \left( 1,x^{2}\right) \rho _{\mu }\left( z,dx\right)
<+\infty $} on $\left( \mathbb{R},\mathcal{B}\left( \mathbb{R}\right)
\right) $ satisfying $\int_{Z}x^{2}\rho _{\mu }\left( z,dx\right) <+\infty $%
, (ii) for every $C\in \mathcal{B}\left( \mathbb{R}\right) $, $\rho _{\mu
}\left( \cdot ,C\right) $ is a Borel measurable function, (iii) for every
positive function $g\left( z,x\right) \in \mathcal{Z}\otimes \mathcal{B}%
\left( \mathbb{R}\right) $,
\begin{equation}
\int_{Z}\int_{\mathbb{R}}g\left( z,x\right) \rho _{\mu }\left( z,dx\right)
\mu \left( dz\right) =\int_{Z}\int_{\mathbb{R}}g\left( z,x\right) \nu \left(
dz,dx\right) .  \label{LKdensity}
\end{equation}

\item For every $\left( \lambda ,z\right) \in \mathbb{R}\times Z$, define
\begin{equation}
K_{\mu }\left( \lambda ,z\right) =-\frac{\lambda ^{2}}{2}\sigma _{\mu
}^{2}\left( z\right) +\int_{\mathbb{R}}\left( e^{i\lambda x}-1-i\lambda
x\right) \rho _{\mu }\left( z,dx\right) \text{,}  \label{innerLK}
\end{equation}%
where $\sigma _{\mu }^{2}\left( z\right) =\frac{dc^{2}}{d\mu }\left(
z\right) $, then, for every $h\in \mathfrak{H}_{\mu }=L^{2}\left( Z,\mathcal{%
Z},\mu \right) $, $\int_{Z}\left\vert K_{\mu }\left( \lambda h\left(
z\right) ,z\right) \right\vert \mu \left( dz\right) <+\infty $ and the
exponent $\psi _{\mathfrak{H}_{\mu }}$ in (\ref{LK3}) is given by%
\begin{eqnarray}
\psi _{\mathfrak{H}_{\mu }}\left( h;\lambda \right) &=&\int_{Z}K_{\mu
}\left( \lambda h\left( z\right) ,z\right) \mu \left( dz\right)  \label{LK5*}
\\
&=&-\frac{\lambda ^{2}}{2}\int_{Z}\sigma _{\mu }^{2}\left( z\right) \mu
\left( dz\right) +\int_{Z}\int_{\mathbb{R}}\left( e^{i\lambda x}-1-i\lambda
x\right) \rho _{\mu }\left( z,dx\right) \mu \left( dz\right)  \notag
\end{eqnarray}
\end{enumerate}
\end{proposition}

\begin{proof}
The proof follows from results contained in \cite[Section II]{Raj Ros}.
Point 1 is indeed a direct consequence of \cite[Proposition 2.1 (a)]{Raj Ros}%
. In particular, whenever $B\in \mathcal{Z}$ is such that $\mu \left(
B\right) =0$, then $M\left( B\right) =0$, a.s.-$\mathbb{P}$ (by applying (%
\ref{intersection}) with $B=C$), and therefore $c^{2}\left( B\right) =0$,
thus implying $c^{2}\ll \mu $. Point 2 follows from the first part of the
statement of \cite[Lemma 2.3]{Raj Ros}. To establish Point 3 define, as in
\cite[p. 456]{Raj Ros},%
\begin{equation*}
\gamma \left( A\right) =c^{2}\left( A\right) +\int_{\mathbb{R}}\min \left(
1,x^{2}\right) \nu _{A}\left( dx\right) \text{,}
\end{equation*}%
whenever $A\in \mathcal{Z}_{\mu }$, and observe (see \cite[Definition 2.2]%
{Raj Ros}) that $\gamma \left( \cdot \right) $ can be canonically extended
to a $\sigma $-finite and positive measure on $\left( Z,\mathcal{Z}\right) $%
. Moreover, since $\mu \left( B\right) =0$ implies $M\left( B\right) =0$ \
a.s.-$\mathbb{P}$, the uniqueness of the L\'{e}vy-Khinchine characteristics
implies as before $\gamma \left( A\right) =0$, and therefore $\gamma \left(
dz\right) \ll \mu \left( dz\right) $. Observe also that, by standard
arguments, one can select a version of the density $\left( d\gamma /d\mu
\right) \left( z\right) $ such that $\left( d\gamma /d\mu \right) \left(
z\right) <+\infty $ for every $z\in Z$. According to \cite[Lemma 2.3]{Raj
Ros}, there exists a function $\rho :Z\times \mathcal{B}\left( \mathbb{R}%
\right) \mapsto \left[ 0,+\infty \right] $, such that: (a) $\rho \left(
z,\cdot \right) $ is a L\'{e}vy measure on $\mathcal{B}\left( \mathbb{R}%
\right) $ for every $z\in Z$, (b) $\rho \left( \cdot ,C\right) $ is a Borel
measurable function for every $C\in \mathcal{B}\left( \mathbb{R}\right) $,
(c) for every positive function $g\left( z,x\right) \in \mathcal{Z}\otimes
\mathcal{B}\left( \mathbb{R}\right) $,
\begin{equation}
\int_{Z}\int_{\mathbb{R}}g\left( z,x\right) \rho \left( z,dx\right) \gamma
\left( dz\right) =\int_{Z}\int_{\mathbb{R}}g\left( z,x\right) \nu \left(
dz,dx\right) .  \label{ffpp}
\end{equation}%
In particular, by using (\ref{ffpp}) in the case $g\left( z,x\right) =%
\mathbf{1}_{A}\left( z\right) x^{2}$ for $A\in \mathcal{Z}_{\mu }$,%
\begin{equation*}
\int_{A}\int_{\mathbb{R}}x^{2}\rho \left( z,dx\right) \gamma \left(
dz\right) =\int_{\mathbb{R}}x^{2}\nu _{A}\left( dx\right) <+\infty \text{,}
\end{equation*}%
since $M\left( A\right) \in L^{2}\left( \mathbb{P}\right) $, and we deduce
that $\rho $ can be chosen in such a way that, for every $z\in Z$, $\int_{%
\mathbb{R}}x^{2}\rho \left( z,dx\right) <+\infty $. Now define, for every $%
z\in Z$ and $C\in \mathcal{B}\left( \mathbb{R}\right) $,
\begin{equation*}
\rho _{\mu }\left( z,C\right) =\frac{d\gamma }{d\mu }\left( z\right) \rho
\left( z,C\right) \text{,}
\end{equation*}%
and observe that, due to the previous discussion, the application $\rho
_{\mu }:Z\times \mathcal{B}\left( \mathbb{R}\right) \mapsto \left[ 0,+\infty %
\right] $ trivially satisfies properties (i)-(iii) in the statement of Point
3, which is therefore proved. To Prove point 4, first define a function $%
h\in \mathfrak{H}_{\mu }$ to be \textit{simple }if $h\left( z\right)
=\sum_{i=1}^{n}a_{i}\mathbf{1}_{A_{i}}\left( z\right) $, where $a_{i}\in
\mathbb{R}$, and $\left( A_{1},...,A_{n}\right) $ is a finite collection of
disjoints elements of $\mathcal{Z}_{\mu }$. Of course, the class of simple
functions (which is a linear space) is dense in $\mathfrak{H}_{\mu }$, and
therefore for every $h\in \mathfrak{H}_{\mu }$ there exists a sequence $%
h_{n} $, $n\geq 1$, of simple functions such that $\int_{Z}\left(
h_{n}\left( z\right) -h\left( z\right) \right) ^{2}\mu \left( dz\right)
\rightarrow 0$. As a consequence, since $\mu $ is $\sigma $-finite there
exists a subsequence $n_{k}$ such that $h_{n_{k}}\left( z\right) \rightarrow
h\left( z\right) $ for $\mu $-a.e. $z\in Z$ (and therefore for $\gamma $%
-a.e. $z\in Z $) and moreover, for every $A\in \mathcal{Z}$, the random
sequence $X_{M}\left( \mathbf{1}_{A}h_{n}\right) $ (where we use the
notation (\ref{isomeasure})) is a Cauchy sequence in $L^{2}\left( \mathbb{P}%
\right) $, and hence it converges in probability. In the terminology of \cite%
[p. 460]{Raj Ros}, this implies that every $h\in \mathfrak{H}_{\mu }$ is $M$%
-integrable, and that, for every $A\in \mathcal{Z}$, the random variable $%
X_{M}\left( h\mathbf{1}_{A}\right) $, defined according to (\ref{isomeasure}%
), coincides with $\int_{A}h\left( z\right) M\left( dz\right) $, i.e. the
integral of $h$ with respect to the restriction of $M\left( \cdot \right) $
to $A$, as defined in \cite[p. 460]{Raj Ros}. As a consequence, by using a
slight modification of \cite[Proposition 2.6]{Raj Ros}\footnote{%
The difference lies in the choice of the truncation.}, the function $K_{0}$
on $\mathbb{R}\times Z$ given by
\begin{equation*}
K_{0}\left( \lambda ,z\right) =-\frac{\lambda ^{2}}{2}\sigma _{0}^{2}\left(
z\right) +\int_{\mathbb{R}}\left( e^{i\lambda x}-1-i\lambda x\right) \rho
\left( z,dx\right) \text{,}
\end{equation*}%
where $\sigma _{0}^{2}\left( z\right) =\left( dc^{2}/d\gamma \right) \left(
z\right) $, is such that $\int_{Z}\left\vert K_{0}\left( \lambda h\left(
z\right) ,z\right) \right\vert \gamma \left( dz\right) <+\infty $ for every $%
h\in \mathfrak{H}_{\mu }$, and also
\begin{equation*}
\mathbb{E}\left[ \exp \left( i\lambda X_{M}\left( h\right) \right) \right]
=\int_{Z}K_{0}\left( \lambda h\left( z\right) ,z\right) \gamma \left(
dz\right) .
\end{equation*}%
Relation (\ref{LK3}) and the fact that, by definition,
\begin{equation*}
K_{\mu }\left( \lambda h\left( z\right) ,z\right) =K_{0}\left( \lambda
h\left( z\right) ,z\right) \frac{d\gamma }{d\mu }\left( z\right) \text{, \ \
}\forall z\in Z\text{, }\forall h\in \mathfrak{H}_{\mu }\text{, }\forall
\lambda \in \mathbb{R}\text{,}
\end{equation*}%
yield (\ref{LK5*}).
\end{proof}

\bigskip

\textbf{Examples -- }(a) If $M$ is a centered Gaussian measure with control $%
\mu $, then $\nu =0$ and, for $h\in \mathfrak{H}_{\mu }$,%
\begin{equation*}
\psi _{\mathfrak{H}_{\mu }}\left( h;\lambda \right) =-\frac{\lambda ^{2}}{2}%
\int_{Z}h^{2}\left( z\right) \mu \left( dz\right) .
\end{equation*}

(b) If $M$ is a centered Poisson measure with control $\mu $, then $%
c^{2}\left( \cdot \right) =0$ and $\rho _{\mu }\left( z,dx\right) =\delta
_{1}\left( dx\right) $ for all $z\in Z$, where $\delta _{1}$ is the Dirac
mass at $x$, and therefore, for $h\in \mathfrak{H}_{\mu }$,%
\begin{equation*}
\psi _{\mathfrak{H}_{\mu }}\left( h;\lambda \right) =\int_{Z}\left(
e^{i\lambda h\left( z\right) }-1-i\lambda h\left( z\right) \right) \mu
\left( dz\right) .
\end{equation*}

For instance, one can take $Z=\left[ 0,+\infty \right) \times \mathbb{%
R\times R}$, and $\mu \left( dx,du,dw\right) =dxdu\nu \left( dw\right) $,
where $\nu \left( dw\right) =\mathbf{1}_{\left\vert w\right\vert
<1}\left\vert w\right\vert ^{-\left( 1+\alpha \right) }dw$ and $\alpha \in
\left( 0,2\right) $. In this case, the centered Poisson measure\ $M$
generates the (standard)\ \textit{Poissonized Telecom process} $\left\{
Y_{P,\alpha }\left( t\right) :t\geq 0\right\} $, defined in \cite[Section 4.1%
]{CohTaq} as%
\begin{eqnarray*}
Y_{P,\alpha }\left( t\right) &=&\int_{0}^{\infty }\int_{\mathbb{R}}\int_{%
\mathbb{R}}\left( \left( \left( t+u\right) \wedge 0+x\right) _{+}\right. \\
&&\left. -\left( u\wedge 0+x\right) _{+}\right) x^{-\left( 1-\kappa \right)
-1/\alpha }wM\left( dx,du,dw\right) \text{,}
\end{eqnarray*}%
with $\kappa \in \left( 0,1-1/\alpha \right) $.

\bigskip

We now want to define multiple integrals, of functions vanishing on
diagonals, with respect to the random measure $M$. To this end, fix $d\geq 2$
and set $\mu ^{d}$ to be the canonical product measure on $\left( Z^{d},%
\mathcal{Z}^{d}\right) $ induced by $\mu $. We introduce the following
standard notation: (i) $L^{2}\left( \mu ^{d}\right) \triangleq L^{2}\left(
Z^{d},\mathcal{Z}^{d},\mu ^{d}\right) $ is the class of real-valued and
square-integrable functions on $\left( Z^{d},\mathcal{Z}^{d}\right) $; (ii) $%
L_{s}^{2}\left( \mu ^{d}\right) $ is the subset of $L^{2}\left( \mu
^{d}\right) $ composed of square integrable and symmetric functions; (iii) $%
L_{s,0}^{2}\left( \mu ^{d}\right) $ is the subset of $L_{s}^{2}\left( \mu
^{d}\right) $ composed of square integrable and symmetric functions
vanishing on diagonals.

\bigskip

Now define $\mathcal{S}_{s,0}\left( \mu ^{d}\right) $ to be subset of $%
L_{s,0}^{2}\left( \mu ^{d}\right) $ composed of functions with the form%
\begin{equation}
f\left( z_{1},...,z_{d}\right) =\sum_{\sigma \in \mathfrak{S}_{d}}\mathbf{1}%
_{B_{1}}\left( z_{\sigma \left( 1\right) }\right) \cdot \cdot \cdot \mathbf{1%
}_{B_{d}}\left( z_{\sigma \left( d\right) }\right) \text{,}
\label{permsimple}
\end{equation}%
where $B_{1},...,B_{d}\in \mathcal{Z}_{\mu }$ are pairwise disjoint sets,
and $\mathfrak{S}_{d}$ is the group of all permutations of $\left\{
1,...,d\right\} $. Recall (see e.g. \cite[Proposition 3]{RoWa}) that $%
\mathcal{S}_{s,0}\left( \mu ^{d}\right) $ is total in $L_{s,0}^{2}\left( \mu
^{d}\right) $. For $f\in L_{s,0}^{2}\left( \mu ^{d}\right) $ as in (\ref%
{permsimple}), we set
\begin{equation}
I_{d}^{M}\left( f\right) =d!M\left( B_{1}\right) \times M\left( B_{2}\right)
\times \cdot \cdot \cdot \times M\left( B_{d}\right)  \label{MSINT}
\end{equation}%
to be the \textit{multiple integral},\textit{\ }of order\textit{\ }$d$, of $%
f $ with respect to $M$. It is well known (see for instance \cite[Theorem 5]%
{RoWa}) that there exists a unique linear extension of $I_{d}^{M}$, from $%
\mathcal{S}_{s,0}\left( \mu ^{d}\right) $ to $L_{s,0}^{2}\left( \mu
^{d}\right) $, satisfying the following: (a) for every $f\in
L_{s,0}^{2}\left( \mu ^{d}\right) $, $I_{d}^{M}\left( f\right) $ is a
centered and square-integrable random variable, and (b) for every $f,g\in
L_{s,0}^{2}\left( \mu ^{d}\right) $%
\begin{equation*}
\mathbb{E}\left[ I_{d}^{M}\left( f\right) I_{d}^{M}\left( g\right) \right]
=d!\left( f,g\right) _{L^{2}\left( \mu ^{d}\right) }\triangleq
d!\int_{Z^{d}}f\left( \mathbf{z}_{d}\right) g\left( \mathbf{z}_{d}\right)
\mu ^{d}\left( d\mathbf{z}_{d}\right) \text{,}
\end{equation*}%
where $\mathbf{z}_{d}=\left( z_{1},...,z_{d}\right) $ stands for a generic
element of $Z^{d}$. Note that, by construction, if $d\neq d^{\prime }$, $%
\mathbb{E}\left[ I_{d}^{M}\left( f\right) I_{d^{\prime }}^{M}\left( g\right) %
\right] $ $=$ $0$ for every $f\in L_{s,0}^{2}\left( \mu ^{d}\right) $ and
every $g\in L_{s,0}^{2}\left( \mu ^{d^{\prime }}\right) $. Again, for $f\in
L_{s,0}^{2}\left( \mu ^{d}\right) $, $I_{d}^{M}\left( f\right) $ is called
the multiple integral, of order\textit{\ }$d$, of $f$ with respect to $M$.
When $f\in L_{s}^{2}\left( \mu ^{d}\right) $ (hence, $f$ does not
necessarily vanish on diagonals) we define%
\begin{equation}
I_{d}^{M}\left( f\right) \triangleq I_{d}^{M}\left( f\mathbf{1}%
_{Z_{0}^{d}}\right) \text{,}  \label{conv}
\end{equation}%
where
\begin{equation}
Z_{0}^{d}\triangleq \left\{ \left( z_{1},...,z_{d}\right) \in Z^{d}:\text{\
the }z_{j}\text{'s are all different }\right\} ,  \label{zeta0}
\end{equation}%
so that (since $\mu $ is non atomic, and therefore the product measures do
not charge diagonals), for every $f,g\in L^{2}\left( \mu ^{d}\right) $, $%
\mathbb{E}\left[ I_{d}^{M}\left( f\right) I_{d}^{M}\left( g\right) \right] $
$=d!\int_{Z_{0}^{d}}f\left( \mathbf{z}_{d}\right) g\left( \mathbf{z}%
_{d}\right) \mu ^{d}\left( d\mathbf{z}_{d}\right) $ $=d!\left( f,g\right)
_{L^{2}\left( \mu ^{d}\right) }$. Note that, for $d=1$, one usually sets $%
L_{s,0}^{2}\left( \mu ^{1}\right) =L_{s}^{2}\left( \mu ^{1}\right)
=L^{2}\left( \mu ^{1}\right) =\mathfrak{H}_{\mu }$, and $I_{1}^{M}\left(
f\right) =X_{M}\left( f\right) $, $f\in \mathfrak{H}_{\mu }$.

\bigskip

In what follows, we shall show that, for some well chosen
resolutions $\pi \in \mathcal{R}_{X_{M}}\left( \mathfrak{H}_{\mu
}\right) $, every multiple integral of the type $I_{d}^{M}\left(
f\right) $, $f\in L_{s,0}^{2}\left( \mu ^{d}\right) $, can be
represented in the form of a generalized adapted integral of the
kind introduced in Section 3. As a consequence, the asymptotic
behavior of $I_{d}^{M}\left( f\right) $ can be studied by means of
Theorem \ref{T : Main}.

\subsection{Representation of multiple integrals and limit theorems}

Under the notation and assumptions of this section, consider a
\textquotedblleft continuous\textquotedblright\ increasing family $\{Z_{t}:$
$t\in \left[ 0,1\right] \}$ of elements of $\mathcal{Z}$, such that $%
Z_{0}=\varnothing $, $Z_{1}=Z$, $Z_{s}\subseteq Z_{t}$ for $s<t$, and, for
every $g\in L^{1}\left( \mu \right) $ and every $t\in \left[ 0,1\right] $,
\begin{equation}
\lim_{s\rightarrow t}\int_{Z_{s}}g\left( x\right) \mu \left( dx\right)
=\int_{Z_{t}}g\left( x\right) \mu \left( dx\right) .  \label{conti}
\end{equation}

For example, for $Z=\left[ 0,1\right] ^{2}$, one can take $Z_{t}=\left[ 0,t%
\right] ^{2}$ or $Z_{t}=\left[ \left( 1-t\right) /2,\left( 1+t\right) /2%
\right] ^{2}$. To each $t\in \left[ 0,1\right] $, we associate the following
projection operator $\pi _{t}:\mathfrak{H}_{\mu }\mapsto \mathfrak{H}_{\mu }$%
: $\forall f\in \mathfrak{H}_{\mu }$,
\begin{equation}
\pi _{t}f\left( z\right) =\mathbf{1}_{Z_{t}}\left( z\right) f\left( z\right)
\text{, \ \ }z\in Z\text{,}  \label{mmres}
\end{equation}%
so that, since $M$ is independently scattered, the continuous resolution of
the identity $\pi =\left\{ \pi _{t}:t\in \left[ 0,1\right] \right\} $ is
such that, $\pi \in \mathcal{R}_{X_{M}}\left( \mathfrak{H}_{\mu }\right) $.
Note also that, thanks to (\pageref{conti}) and by uniform continuity, for
every $f\in \mathfrak{H}_{\mu }$, every $t\in \left( 0,1\right] $ and every
sequence of partitions of $\left[ 0,t\right] $,%
\begin{equation}
t^{\left( n\right) }=\left\{ 0=t_{0}^{\left( n\right) }<t_{1}^{\left(
n\right) }<...<t_{r_{n}}^{\left( n\right) }=t\right\} \text{, \ \ }n\geq 1%
\text{,}  \label{part}
\end{equation}%
such that $mesh\left( t^{\left( n\right) }\right) \triangleq
\max_{i=0,...,r_{n}-1}\left( t_{i+1}^{\left( n\right) }-t_{i}^{\left(
n\right) }\right) \rightarrow 0$,%
\begin{equation}
\max_{i=0,...,r_{n}-1}\left\Vert \left( \pi _{t_{i+1}^{\left( n\right)
}}-\pi _{t_{i}^{\left( n\right) }}\right) f\right\Vert _{\mathfrak{H}_{\mu
}}^{2}\rightarrow 0,  \label{meshh}
\end{equation}%
and in particular, for every $B\in \mathcal{Z}_{\mu }$,
\begin{equation}
\max_{i=0,...,r_{n}-1}\mu \left( B\cap \left( Z_{t_{i}^{\left( n\right)
}}\backslash Z_{t_{i-1}^{\left( n\right) }}\right) \right) \rightarrow 0.
\label{meshh2}
\end{equation}

The following result contains the key of the subsequent discussion.

\begin{proposition}
\label{P : MIrep}For every $d\geq 2$, every random variable of the form $%
I_{d}^{M}\left( f\mathbf{1}_{Z_{t}^{d}}\right) $, for some $f\in
L_{s,0}^{2}\left( \mu ^{d}\right) $ and $t\in \left( 0,1\right] $, can be
approximated in $L^{2}\left( \mathbb{P}\right) $ by linear combinations of
random variables of the type%
\begin{equation}
M\left( B_{1}\cap Z_{t_{1}}\right) \times M\left( B_{2}\cap \left(
Z_{t_{2}}\backslash Z_{t_{1}}\right) \right) \times \cdot \cdot \cdot \times
M\left( B_{d}\cap \left( Z_{t_{d}}\backslash Z_{t_{d-1}}\right) \right)
\text{,}  \label{MI*}
\end{equation}%
where the $t_{1},...,t_{d}$ are rational, $0\leq t_{1}<t_{2}<\cdot \cdot
\cdot <t_{d}\leq t$ and $B_{1},...,B_{d}\in \mathcal{Z}_{\mu }$ are
disjoint. In particular, $I_{d}^{M}\left( f\mathbf{1}_{Z_{t}^{d}}\right) \in
\mathcal{F}_{t}^{\pi }$, where the filtration $\mathcal{F}_{t}^{\pi }$, $%
t\in \left[ 0,1\right] $, is defined as in (\ref{resfiltration}).
\end{proposition}

\bigskip

\textbf{Remark -- }Observe that, if $f\in \mathcal{S}_{s,0}\left( \mu
^{d}\right) $ is such that
\begin{equation}
f\left( z_{1},...,z_{d}\right) =\sum_{\sigma \in \mathfrak{S}_{d}}\mathbf{1}%
_{B_{1}\cap Z_{t_{1}}}\left( z_{\sigma \left( 1\right) }\right) \cdot \cdot
\cdot \mathbf{1}_{B_{d}\cap \left( Z_{t_{d}}\backslash Z_{t_{d-1}}\right)
}\left( z_{\sigma \left( d\right) }\right) \text{,}  \label{permsym2}
\end{equation}
then, by (\ref{MSINT}),
\begin{equation}
d!M\left( B_{1}\cap Z_{t_{1}}\right) \times M\left( B_{2}\cap \left(
Z_{t_{2}}\backslash Z_{t_{1}}\right) \right) \times \cdot \cdot \cdot \times
M\left( B_{d}\cap \left( Z_{t_{d}}\backslash Z_{t_{d-1}}\right) \right)
=I_{d}^{M}\left( f\right) \text{.}  \label{MI}
\end{equation}

\bigskip

\begin{proof}
Observe first that, for every $f\in L_{s,0}^{2}\left( \mu ^{d}\right) $,
every $t\in \left( 0,1\right] $ and every sequence of rational numbers $%
t_{n}\rightarrow t$, $I_{d}^{M}\left( f\mathbf{1}_{Z_{t_{n}}^{d}}\right)
\rightarrow I_{d}^{M}\left( f\mathbf{1}_{Z_{t}^{d}}\right) $ in $L^{2}\left(
\mathbb{P}\right) $. By density, it is therefore sufficient to prove the
statement for multiple integrals of the type $I_{d}^{M}\left( f\mathbf{1}%
_{Z_{t}^{d}}\right) $, where $t\in \mathbb{Q\cap }\left( 0,1\right] $ and $%
f\in \mathcal{S}_{s,0}\left( \mu ^{d}\right) $ is as in (\ref{permsimple}).
Start with $d=2$. In this case,
\begin{equation*}
\frac{1}{2}I_{2}^{M}\left( f\mathbf{1}_{Z_{t}^{2}}\right) =M\left( B_{1}\cap
Z_{t}\right) M\left( B_{2}\cap Z_{t}\right)
\end{equation*}%
with $B_{1},B_{2}$ disjoints, and also, for every partition $\left\{
0=t_{0}<t_{1}<...<t_{r}=t\right\} $ (with $r\geq 1$) of $\left[ 0,t\right] $,%
\begin{eqnarray*}
\frac{1}{2}I_{2}^{M}\left( f\right) &=&\sum_{i=1}^{r}M\left( B_{1}\cap
\left( Z_{t_{i}}\backslash Z_{t_{i-1}}\right) \right) \sum_{j=1}^{r}M\left(
B_{2}\cap \left( Z_{t_{j}}\backslash Z_{t_{j-1}}\right) \right) \\
&=&\sum_{1\leq i\neq j\leq r}M\left( B_{1}\cap \left( Z_{t_{i}}\backslash
Z_{t_{i-1}}\right) \right) M\left( B_{2}\cap \left( Z_{t_{j}}\backslash
Z_{t_{j-1}}\right) \right) + \\
&&+\sum_{i=1}^{r}M\left( B_{1}\cap \left( Z_{t_{i}}\backslash
Z_{t_{i-1}}\right) \right) M\left( B_{2}\cap \left( Z_{t_{i}}\backslash
Z_{t_{i-1}}\right) \right) \triangleq \Sigma _{1}+\Sigma _{2}.
\end{eqnarray*}%
The summands in the first sum $\Sigma _{1}$ have the desired form (\ref{MI*}%
). It is therefore sufficient to prove that for every sequence of partitions
$t^{\left( n\right) }$, $n\geq 1$, as in (\ref{part}) and such that $%
mesh\left( t^{\left( n\right) }\right) \rightarrow 0$ and the $t_{1}^{\left(
n\right) },...,t_{r_{n}}^{\left( n\right) }$ are rational,%
\begin{equation}
\lim_{n\rightarrow \infty }\mathbb{E}\left[ \left( \sum_{i=1}^{r_{n}}M\left(
B_{1}\cap \left( Z_{t_{i}^{\left( n\right) }}\backslash Z_{t_{i-1}^{\left(
n\right) }}\right) \right) M\left( B_{2}\cap \left( Z_{t_{i}^{\left(
n\right) }}\backslash Z_{t_{i-1}^{\left( n\right) }}\right) \right) \right)
^{2}\right] =0.  \label{MI**}
\end{equation}%
Since $B_{1}$ and $B_{2}$ are disjoint, and thanks to the isometric
properties of $M$,%
\begin{eqnarray*}
&&\mathbb{E}\left[ \left( \sum_{i=1}^{r_{n}}M\left( B_{1}\cap \left(
Z_{t_{i}^{\left( n\right) }}\backslash Z_{t_{i-1}^{\left( n\right) }}\right)
\right) M\left( B_{2}\cap \left( Z_{t_{i}^{\left( n\right) }}\backslash
Z_{t_{i-1}^{\left( n\right) }}\right) \right) \right) ^{2}\right] \\
&=&\sum_{i=1}^{r_{n}}\mathbb{E}\left[ M\left( B_{1}\cap \left(
Z_{t_{i}^{\left( n\right) }}\backslash Z_{t_{i-1}^{\left( n\right) }}\right)
\right) ^{2}M\left( B_{2}\cap \left( Z_{t_{i}^{\left( n\right) }}\backslash
Z_{t_{i-1}^{\left( n\right) }}\right) \right) ^{2}\right] \\
&=&\sum_{i=1}^{r_{n}}\mu \left( B_{1}\cap \left( Z_{t_{i}^{\left( n\right)
}}\backslash Z_{t_{i-1}^{\left( n\right) }}\right) \right) \mu \left(
B_{2}\cap \left( Z_{t_{i}^{\left( n\right) }}\backslash Z_{t_{i-1}^{\left(
n\right) }}\right) \right) \\
&\leq &\mu \left( B_{1}\right) \max_{i=1,...,r_{n}}\mu \left( B_{2}\cap
\left( Z_{t_{i}^{\left( n\right) }}\backslash Z_{t_{i-1}^{\left( n\right)
}}\right) \right) \rightarrow 0\text{,}
\end{eqnarray*}%
thanks to (\ref{meshh2}). Now fix $d\geq 3$, and consider a random variable
of the type
\begin{equation}
F=M\left( B_{1}\cap Z_{t}\right) \times \cdot \cdot \cdot \times M\left(
B_{d-1}\cap Z_{t}\right) \times M\left( B_{d}\cap Z_{t}\right) \text{,}
\label{MI***}
\end{equation}%
where $B_{1},...,B_{d}\in \mathcal{Z}_{\mu }$ are disjoint. The above
discussion yields that $F$ can be approximated by linear combinations of
random variables of the type%
\begin{eqnarray}
&&M\left( B_{1}\cap Z_{t}\right) \times \cdot \cdot \cdot \times M\left(
B_{d-3}\cap Z_{t}\right) \times  \label{MI****} \\
&&\times \left[ M\left( B_{d-2}\right) \times M\left( B_{d-1}\cap \left(
Z_{s}\backslash Z_{r}\right) \right) \times M\left( B_{d}\cap \left(
Z_{v}\backslash Z_{u}\right) \right) \right] \text{,}  \notag
\end{eqnarray}%
where $r<s<u<v\leq t$ are rational. We will proceed by induction focusing
first on the terms in the brackets in (\ref{MI****}). Express $Z_{t}$ as the
union of five disjoint sets $Z_{t}=\left( Z_{t}\backslash Z_{v}\right) $ $%
\cup \left( Z_{v}\backslash Z_{u}\right) $ $\cup \left( Z_{u}\backslash
Z_{s}\right) $ $\cup \left( Z_{s}\backslash Z_{r}\right) $ $\cup Z_{r}$, and
decompose $M\left( B_{d-2}\cap Z_{t}\right) $ accordingly. One gets%
\begin{eqnarray}
&&M\left( B_{d-2}\cap Z_{t}\right) M\left( B_{d-1}\cap \left(
Z_{t}\backslash Z_{r}\right) \right) M\left( B_{d}\cap \left(
Z_{v}\backslash Z_{u}\right) \right)  \label{f} \\
&=&M\left( B_{d-2}\cap \left( Z_{s}\backslash Z_{r}\right) \right) M\left(
B_{d-1}\cap Z_{s}\backslash Z_{r}\right) M\left( B_{d}\cap \left(
Z_{v}\backslash Z_{u}\right) \right)  \notag \\
&&+M\left( B_{d-2}\cap \left( Z_{v}\backslash Z_{u}\right) \right) M\left(
B_{d-1}\cap Z_{s}\backslash Z_{r}\right) M\left( B_{d}\cap \left(
Z_{v}\backslash Z_{u}\right) \right)  \notag \\
&&+M\left( B_{d-2}\cap \left( Z_{u}\backslash Z_{s}\right) \right) M\left(
B_{d-1}\cap Z_{s}\backslash Z_{r}\right) M\left( B_{d}\cap \left(
Z_{v}\backslash Z_{u}\right) \right)  \notag \\
&&+M\left( B_{d-2}\cap \left( Z_{t}\backslash Z_{v}\right) \right) M\left(
B_{d-1}\cap Z_{s}\backslash Z_{r}\right) M\left( B_{d}\cap \left(
Z_{v}\backslash Z_{u}\right) \right)  \notag \\
&&+M\left( B_{d-2}\cap \left( Z_{r}\backslash Z_{0}\right) \right) M\left(
B_{d-1}\cap Z_{s}\backslash Z_{r}\right) M\left( B_{d}\cap \left(
Z_{v}\backslash Z_{u}\right) \right) .  \notag
\end{eqnarray}%
Observe that the last three summands involve disjoint subsets of $Z$ and
hence are of the form (\ref{MI*}). Since each of the first two summands
involve two identical subsets of $Z$ (e.\.{g}. $\left( Z_{s}\backslash
Z_{r}\right) $) and a disjoint subset (e.\.{g}. $\left( Z_{s}\backslash
Z_{r}\right) $), they can be dealt with in the same way as (\ref{MI**})
above. Thus, linear combinations of the five summands on the RHS of (\ref{f}%
) can be approximated by linear combinations of random variables of the type%
\begin{equation*}
M\left( C_{1}\cap \left( Z_{t_{2}}\backslash Z_{t_{1}}\right) \right)
M\left( C_{2}\cap \left( Z_{t_{3}}\backslash Z_{t_{2}}\right) \right)
M\left( C_{3}\cap \left( Z_{t_{3}}\backslash Z_{t_{2}}\right) \right) ,
\end{equation*}%
where $C_{1},C_{2},C_{3}\in \mathcal{Z}_{\mu }$ are disjoints, and $%
t_{1}<t_{2}<t_{3}\leq t$ are rational. The general result is obtained by
recurrence.
\end{proof}

\bigskip

Proposition \ref{P : MIrep} will be used to prove that, whenever there
exists $\pi \in \mathcal{R}_{X_{M}}\left( \mathfrak{H}_{\mu }\right) $
defined as in formula (\ref{mmres}), multiple integrals can be represented
as generalized adapted integrals of the kind described in Section 3. To do
this, we introduce a partial ordering on $Z$ as follows: for every $%
z,z^{\prime }\in Z$,
\begin{equation}
z\prec _{\pi }z^{\prime }  \label{PO}
\end{equation}%
if, and only if, there exists $t\in \mathbb{Q\cap }\left( 0,1\right) $ such
that $z\in Z_{t}$ and $z^{\prime }\in Z_{t}^{c}$, where $Z_{t}^{c}$ stands
for the complement of $Z_{t}$. For a fixed $d\geq 2$, we define the $\pi $-%
\textit{purely non-diagonal }subset of $Z^{d}$ as
\begin{equation*}
Z_{\pi ,0}^{d}=\left\{ \left( z_{1},...,z_{d}\right) \in Z^{d}:z_{\sigma
\left( 1\right) }\prec _{\pi }z_{\sigma \left( 2\right) }\prec _{\pi }\cdot
\cdot \cdot \prec _{\pi }z_{\sigma \left( d\right) }\text{, \ for some }%
\sigma \in \mathfrak{S}_{d}\right\} \text{.}
\end{equation*}%
Note that $Z_{\pi ,0}^{d}\in \mathcal{Z}^{d}$, and also that not every pair
of distinct points of $Z$ can be ordered, that is, in general, $Z_{\pi
,0}^{d}\neq Z_{0}^{d}$, where $d\geq 2$ and $Z_{0}^{d}$ is defined in (\ref%
{zeta0}) (for illustration, think of $Z=\left[ 0,1\right] ^{2}$, $Z_{t}=%
\left[ 0,t\right] ^{2}$, $t\in \left[ 0,1\right] $; indeed $\left( \left(
1/8,1/4\right) ,\left( 1/4,1/4\right) \right) \in Z_{0}^{2}$, but $\left(
1/4,1/4\right) $ and $\left( 1/8,1/4\right) $ cannot be ordered). However,
because of the continuity condition (\ref{conti}) and for every $d\geq 2$,
the class of the elements of $Z_{0}^{d}$ whose components cannot be ordered
has measure $\mu ^{d}$ equal to zero, as indicated by the following
corollary.

\begin{corollary}
\label{C : diagSets}For every $d\geq 2$ and every $f\in L_{s,0}^{2}\left(
\mu ^{d}\right) $,
\begin{equation*}
I_{d}^{M}\left( f\right) =I_{d}^{M}\left( f\mathbf{1}_{Z_{\pi
,0}^{d}}\right) .
\end{equation*}%
As a consequence, $\mu ^{d}\left( Z_{0}^{d}\backslash Z_{\pi ,0}^{d}\right)
=0$, where $Z_{0}$ is defined in (\ref{zeta0}).
\end{corollary}

\begin{proof}
First observe that the class of r.v.'s of the type $I_{d}^{M}\left( f\mathbf{%
1}_{Z_{\pi ,0}^{d}}\right) $, $f\in L_{s,0}^{2}\left( \mu ^{d}\right) $ is a
closed vector space. Plainly, every $f\in \mathcal{S}_{s,0}\left( \mu
^{d}\right) $ with the form (\ref{permsym2}) is such that $f\left( \mathbf{z}%
_{d}\right) =f\left( \mathbf{z}_{d}\right) \mathbf{1}_{Z_{\pi ,0}^{d}}\left(
\mathbf{z}_{d}\right) $ for every $\mathbf{z}_{d}\in Z^{d}$. Since, by
Proposition \ref{P : MIrep} and relation (\ref{MI}), the class of functions
of the type (\ref{permsym2}) are total in $L_{s,0}^{2}\left( \mu ^{d}\right)
$, the result is obtained by a density argument. The last assertion follows
from the facts that $\forall f\in L_{s,0}^{2}\left( \mu ^{d}\right) $ one
has $f=f\mathbf{1}_{Z_{\pi ,0}^{d}}$, a.e.-$\mu ^{d}$, and $I_{d}^{M}\left(
f\right) =I_{d}^{M}\left( g\right) $ if and only if $f=g$, a.e.-$\mu ^{d}$.
\end{proof}

\bigskip

For $\pi \in \mathcal{R}_{X_{M}}\left( \mathfrak{H}_{\mu }\right) $ as in
formula (\ref{mmres}), the vector spaces $L_{\pi }^{2}\left( \mathfrak{H}%
_{\mu },X_{M}\right) $ and $\mathcal{E}_{\pi }\left( \mathfrak{H}_{\mu
},X_{M}\right) $, composed respectively of adapted and elementary adapted
elements of $L^{2}\left( \mathfrak{H}_{\mu },X_{M}\right) $, are defined as
in Section 3 (in particular, \textit{via }formulae (\ref{adaptation}) and (%
\ref{elementaryad})). Recall that, according to Lemma \ref{L : adapt}, the
closure of $\mathcal{E}_{\pi }\left( \mathfrak{H}_{\mu },X_{M}\right) $
coincides with $L_{\pi }^{2}\left( \mathfrak{H}_{\mu },X_{M}\right) $. For
every $h\in L_{\pi }^{2}\left( \mathfrak{H}_{\mu },X_{M}\right) $, the
random variable $J_{X_{M}}^{\pi }\left( h\right) $ is defined by means of
Proposition \ref{P : ItoInt} and formula (\ref{elSkorohod}). The following
result states that every multiple integral with respect to $M$ is indeed a
generalized adapted integral of the form $J_{X_{M}}^{\pi }\left( h\right) $,
for some $h\in L_{\pi }^{2}\left( \mathfrak{H}_{\mu },X_{M}\right) $. In
what follows, for every $d\geq 1$, every $f\in L_{s,0}^{2}\left( \mu
^{d}\right) $ and every fixed $z\in Z$, the symbol $f\left( z,\cdot \right)
\mathbf{1}\left( \cdot \prec _{\pi }z\right) $ stands for the element of $%
L_{s,0}^{2}\left( \mu ^{d-1}\right) $, given by
\begin{equation}
\left( z_{1},...,z_{d-1}\right) \mapsto f\left( z,z_{1},...,z_{d-1}\right)
\prod_{j=1}^{d-1}\mathbf{1}_{\left( z_{j}\prec _{\pi }z\right) }.
\label{projected K}
\end{equation}

\bigskip

\begin{proposition}
\label{P : measClOc}Fix $d\geq 2$, and let $f\in L_{s,0}^{2}\left( \mu
^{d}\right) $. Then,

\begin{enumerate}
\item the random function%
\begin{equation}
z\mapsto h_{\pi }\left( f\right) \left( z\right) =d\times I_{d-1}^{M}\left(
f\left( z,\cdot \right) \mathbf{1}\left( \cdot \prec _{\pi }z\right) \right)
\text{, \ \ }z\in Z\text{,}  \label{COF}
\end{equation}%
is an element of $L_{\pi }^{2}\left( \mathfrak{H}_{\mu },X_{M}\right) $;

\item $I_{M}^{d}\left( f\right) =J_{X_{M}}^{\pi }\left( h_{\pi }\left(
f\right) \right) $, where $h_{\pi }\left( f\right) $ is defined as in (\ref%
{COF}).
\end{enumerate}

Moreover, if a random variable $F\in L^{2}\left( \mathbb{P}\right) $ has the
form $F=\sum_{d=1}^{\infty }I_{M}^{d}\left( f^{\left( d\right) }\right) $,
where $f^{\left( d\right) }\in L_{s,0}^{2}\left( \mu ^{d}\right) $ for $%
d\geq 1$ and the series is convergent in $L^{2}\left( \mathbb{P}\right) $,
then
\begin{equation}
F=J_{X_{M}}^{\pi }\left( h_{\pi }\left( F\right) \right) \text{,}
\label{COF2}
\end{equation}%
where
\begin{equation}
h_{\pi }\left( F\right) \left( z\right) =\sum_{d=1}^{\infty }h_{\pi }\left(
f^{\left( d\right) }\right) \left( z\right) \text{, \ \ }z\in Z\text{,}
\label{P14}
\end{equation}%
and the series in (\ref{P14}) is convergent in $L_{\pi }^{2}\left( \mathfrak{%
H}_{\mu },X_{M}\right) $.
\end{proposition}

\begin{proof}
It is clear that $h_{\pi }\left( f\right) \in L^{2}\left( \mathfrak{H}_{\mu
},X_{M}\right) $ (the class of square integrable, but not necessarily
adapted processes). Now observe that, thanks to Proposition \ref{P : MIrep},
if $g\in L_{s,0}^{2}\left( \mu ^{d}\right) $ has support in $Z_{t}^{d}$ for
some $t\in \left( 0,1\right] $, then $I_{M}^{d}\left( g\right) \in \mathcal{F%
}_{t}^{\pi }$. As a consequence, since for any fixed $z\in Z_{t}$, $t\in
\left( 0,1\right] $, the symmetric function (on $Z^{d-1}$) $f\left( z,\cdot
\right) \mathbf{1}\left( \cdot \prec _{\pi }z\right) $ has support in $%
Z_{t}^{d}$, for every $b\in \mathfrak{H}_{\mu }$ and $t\in \left( 0,1\right]
$,
\begin{equation*}
\left( h_{\pi }\left( f\right) ,\pi _{t}b\right) _{\mathfrak{H}_{\mu
}}=\int_{Z_{t}}h_{\pi }\left( f\right) \left( z\right) b\left( z\right) \mu
\left( dz\right) =d\int_{Z_{t}}b\left( z\right) I_{d-1}^{M}\left( f\left(
z,\cdot \right) \mathbf{1}\left( \cdot \prec _{\pi }z\right) \right) \mu
\left( dz\right) \in \mathcal{F}_{t}^{\pi }\text{,}
\end{equation*}%
and therefore $h_{\pi }\left( f\right) \in L_{\pi }^{2}\left( \mathfrak{H}%
_{\mu },X_{M}\right) $. This proves Point 1. By density, it is sufficient to
prove Point 2 for random variables of the type $I_{M}^{d}\left( f\right) $,
where $f\in \mathcal{S}_{s,0}\left( \mu ^{d}\right) $ is as in (\ref%
{permsym2}). Indeed, for such an $f$ and for every $\left(
z,z_{1},...,z_{d-1}\right) \in Z^{d}$%
\begin{eqnarray*}
&&f\left( z,z_{1},...,z_{d-1}\right) \prod_{j=1}^{d-1}\mathbf{1}_{\left(
z_{j}\prec _{\pi }z\right) } \\
&=&\sum_{\sigma \in \mathfrak{S}_{d-1}}\mathbf{1}_{B_{1}\cap
Z_{t_{1}}}\left( z_{\sigma \left( 1\right) }\right) \cdot \cdot \cdot
\mathbf{1}_{B_{d-1}\cap \left( Z_{t_{d-1}}\backslash Z_{t_{d-2}}\right)
}\left( z_{\sigma \left( d-1\right) }\right) \mathbf{1}_{B_{d}\cap \left(
Z_{t_{d}}\backslash Z_{t_{d-1}}\right) }\left( z\right) ,
\end{eqnarray*}%
so that
\begin{equation*}
d\times h_{\pi }\left( f\right) \left( z\right) =d\left( d-1\right) !M\left(
B_{1}\cap Z_{t_{1}}\right) \times \cdot \cdot \cdot \times M\left(
B_{d-1}\cap \left( Z_{t_{d-1}}\backslash Z_{t_{d-2}}\right) \right) \mathbf{1%
}_{B_{d}\cap \left( Z_{t_{d}}\backslash Z_{t_{d-1}}\right) }\left( z\right) ,
\end{equation*}%
and finally, thanks to (\ref{elSkorohod}) and (\ref{MI}),
\begin{equation*}
J_{X_{M}}^{\pi }\left( h_{\pi }\left( f\right) \right) =d!M\left( B_{1}\cap
Z_{t_{1}}\right) \times \cdot \cdot \cdot \times M\left( B_{d}\cap \left(
Z_{t_{d}}\backslash Z_{t_{d-1}}\right) \right) =I_{M}^{d}\left( f\right)
\text{.}
\end{equation*}

The last assertion in the statement is an immediate consequence of the
orthogonality relations between multiple integrals of different orders.
\end{proof}

\bigskip

\textbf{Remarks -- }(1) Formula (\ref{COF}) implies that, for $t\in \left[
0,1\right] $ and $f\in L_{s,0}^{2}\left( \mu ^{d}\right) $,
\begin{equation*}
I_{M}^{d}\left( f\mathbf{1}_{Z_{t}^{d}}\right) =J_{X_{M}}^{\pi }\left( \pi
_{t}h_{\pi }\left( f\right) \right) ,
\end{equation*}%
and therefore, since $t\mapsto J_{X_{M}}^{\pi }\left( \pi _{t}h_{\pi }\left(
f\right) \right) $ is a $\mathcal{F}_{t}^{\pi }$-martingale (see Corollary %
\ref{C : mart}),%
\begin{equation}
\mathbb{E}\left[ I_{M}^{d}\left( f\right) \mid \mathcal{F}_{t}^{\pi }\right]
=I_{M}^{d}\left( f\mathbf{1}_{Z_{t}^{d}}\right) \text{, \ \ }t\in \left[ 0,1%
\right] \text{.}  \label{Mart}
\end{equation}

(2)\textbf{\ }The random process $z\mapsto dI_{d-1}^{M}\left( f\left(
z,\cdot \right) \right) \triangleq D_{z}I_{d}^{M}\left( f\right) $ is a
\textquotedblleft formal\textquotedblright\ Malliavin-Shikegawa derivative
of the random variable $I_{d}^{M}\left( f\right) $, whereas $z\mapsto
dI_{d-1}^{M}\left( f\left( z,\cdot \right) \mathbf{1}\left( \cdot \prec
_{\pi }z\right) \right) $ is the projection of $D_{z}I_{d}^{M}\left(
f\right) $ on the space of adapted integrands $L_{\pi }^{2}\left( \mathfrak{H%
}_{\mu },X_{M}\right) $. In this sense, formula (\ref{COF2}) can be
interpreted as a \textquotedblleft generalized Clark-Ocone
formula\textquotedblright , in the same spirit of the results proved by L.
Wu in \cite{Wu}. See also the discussion contained in Section 6.

\bigskip

We now state the announced convergence result, which is a consequence of
Proposition \ref{P : measClOc} and Theorem \ref{T : Main}. In what follows, $%
\left( Z_{n},\mathcal{Z}_{n},\mu _{n}\right) $, $n\geq 1$, is a sequence of
measurable spaces and, for each $n$, $M_{n}$ is an independently scattered
random measures on $\left( Z_{n},\mathcal{Z}_{n}\right) $ with control $\mu
_{n}$ (the $M_{n}$'s are defined on the same probability space); also $%
\mathfrak{H}_{\mu _{n}}=L^{2}\left( Z_{n},\mathcal{Z}_{n},\mu _{n}\right) .$
The collection of random variables $X_{M_{n}}=X_{M_{n}}\left( \mathfrak{H}%
_{\mu _{n}}\right) $ is defined through formula (\ref{isomeasure}), with L%
\'{e}vy-Khinchine exponent $\psi _{\mathfrak{H}_{\mu _{n}}}\left( h,\lambda
\right) $, $h\in \mathfrak{H}_{\mu _{n}}$, $\lambda \in \mathbb{R}$, given
by (\ref{LK3}). Moreover, for every $n\geq 1$, $\pi ^{\left( n\right)
}=\left\{ \pi _{t}^{\left( n\right) }:t\in \left[ 0,1\right] \right\} \in
\mathcal{R}_{X_{M_{n}}}\left( \mathfrak{H}_{\mu _{n}}\right) $ is a
continuous resolution of the identity defined as%
\begin{equation}
\pi _{t}^{\left( n\right) }h\left( z\right) =\mathbf{1}_{Z_{n,t}}\left(
z\right) h\left( z\right) \text{, \ \ }z\in Z\text{, \ \ }h\in \mathfrak{H}%
_{\mu _{n}}\text{,}  \label{RR15}
\end{equation}%
where $Z_{n,t}$, $t\in \left[ 0,1\right] $ is an increasing collection of
measurable sets such that $Z_{n,0}=\varnothing $, $Z_{n,1}=Z_{n}$ and
verifying the continuity condition (\ref{conti}).

\bigskip

\begin{theorem}
\label{T : Main 3}Under the previous notation and assumptions, let $d_{n}$, $%
n\geq 1$, be a sequence of natural numbers such that $d_{n}\geq 1$, and let $%
\pi ^{\left( n\right) }\in \mathcal{R}_{X_{M_{n}}}\left( \mathfrak{H}_{\mu
_{n}}\right) $ be as in (\ref{RR15}). Let moreover $f_{d_{n}}^{\left(
n\right) }\in L_{s,0}^{2}\left( \mu _{n}^{d_{n}}\right) $, $n\geq 1$, and
suppose there exists a sequence $\left\{ t_{n}:n\geq 1\right\} $ $\subset %
\left[ 0,1\right] $ and $\sigma $-fields $\left\{ \mathcal{U}_{n}:n\geq
1\right\} $, such that
\begin{equation}
\lim_{n\rightarrow +\infty }d_{n}!\left\Vert f_{d_{n}}^{\left( n\right) }%
\mathbf{1}_{Z_{n,t_{n}}^{d_{n}}}\right\Vert _{L^{2}\left( \mu
_{n}^{d_{n}}\right) }^{2}=0  \label{CO+O}
\end{equation}%
and
\begin{equation}
\mathcal{U}_{n}\subseteq \mathcal{U}_{n+1}\cap \mathcal{F}_{t_{n}}^{\pi
^{\left( n\right) }}\left( X_{M_{n}}\right) .  \label{CO+B}
\end{equation}%
Define also $h_{\pi ^{\left( n\right) }}\left( f_{d_{n}}^{\left( n\right)
}\right) \in L_{\pi ^{\left( n\right) }}^{2}\left( \mathfrak{H}_{\mu
_{n}},X_{M_{n}}\right) $ via formula (\ref{COF}) when $d_{n}\geq 2$, and set
$h_{\pi ^{\left( n\right) }}\left( f_{d_{n}}^{\left( n\right) }\right)
=f_{d_{n}}^{\left( n\right) }$ when $d_{n}=1$. If%
\begin{equation}
\exp \left[ \int_{Z_{n}}K_{\mu _{n}}\left( \lambda h_{\pi ^{\left( n\right)
}}\left( f_{d_{n}}^{\left( n\right) }\right) \left( z\right) ,z\right) \mu
_{n}\left( dz\right) \right] \overset{\mathbb{P}}{\rightarrow }\phi \left(
\lambda ,\omega \right) ,\text{ \ \ }\forall \lambda \in
\mathbb{R}
\text{,}  \label{CO+I}
\end{equation}%
where $K_{\mu _{n}}\left( t,z\right) $, $\left( t,z\right) \in \mathbb{R}%
\times Z$, is given by (\ref{innerLK}), $\phi \in \widehat{\mathbf{M}}_{0}$
and $\phi \left( \lambda \right) \in \vee _{n}\mathcal{U}_{n}\triangleq
\mathcal{U}^{\ast }$, then, as $n\rightarrow +\infty $,
\begin{equation}
\mathbb{E}\left[ \exp \left( i\lambda I_{d_{n}}^{M_{n}}\left(
f_{d_{n}}^{\left( n\right) }\right) \right) \mid \mathcal{F}_{t_{n}}^{\pi
^{\left( n\right) }}\left( X_{M_{n}}\right) \right] \overset{\mathbb{P}}{%
\rightarrow }\phi \left( \lambda \right) \text{, \ \ }\forall \lambda \in
\mathbb{R}
\text{,}  \label{CO+II}
\end{equation}%
and%
\begin{equation}
I_{d_{n}}^{M_{n}}\left( f_{d_{n}}^{\left( n\right) }\right) \rightarrow
_{\left( s,\mathcal{U}^{\ast }\right) }\mathbb{E}\mu \left( \cdot \right) ,
\label{CO+III}
\end{equation}%
where $\mu \in \mathbf{M}$ is as in (\ref{randomfou}).
\end{theorem}

\begin{proof}
It is sufficient to observe that, thanks to Proposition \ref{P : measClOc}, $%
I_{d_{n}}^{M_{n}}\left( f_{d_{n}}^{\left( n\right) }\right)
=J_{X_{M_{n}}}^{\pi ^{\left( n\right) }}\left( h_{\pi ^{\left( n\right)
}}\left( f_{d_{n}}^{\left( n\right) }\right) \right) $, $n\geq 1$. As a
consequence, by using (\ref{Mart}),%
\begin{eqnarray*}
d_{n}!\left\Vert f_{d_{n}}^{\left( n\right) }\mathbf{1}%
_{Z_{n,t_{n}}^{d_{n}}}\right\Vert _{L^{2}\left( \mu _{n}^{d_{n}}\right)
}^{2} &=&\mathbb{E}\left[ I_{d_{n}}^{M_{n}}\left( f_{d_{n}}^{\left( n\right)
}\right) ^{2}\right] \\
&=&\mathbb{E}\left[ J_{X_{M_{n}}}^{\pi ^{\left( n\right) }}\left( \pi
_{t_{n}}^{\left( n\right) }h_{\pi ^{\left( n\right) }}\left(
f_{d_{n}}^{\left( n\right) }\right) \right) ^{2}\right] \\
&=&\left\Vert \pi _{t_{n}}^{\left( n\right) }h_{\pi ^{\left( n\right)
}}\left( f_{d_{n}}^{\left( n\right) }\right) \right\Vert _{L_{\pi ^{\left(
n\right) }}^{2}\left( \mathfrak{H}_{\mu _{n}},X_{M_{n}}\right) }^{2}\text{.}
\end{eqnarray*}%
Moreover, according to Proposition \ref{P : LK},
\begin{equation*}
\int_{Z_{n}}K_{\mu _{n}}\left( \lambda h_{\pi ^{\left( n\right) }}\left(
f_{d_{n}}^{\left( n\right) }\right) \left( z\right) ,z\right) \mu _{n}\left(
dz\right) =\psi _{\mathfrak{H}_{\mu _{n}}}\left( h_{\pi ^{\left( n\right)
}}\left( f_{d_{n}}^{\left( n\right) }\right) ,\lambda \right) .
\end{equation*}%
The conclusion is now a direct consequence of Theorem \ref{T : Main}.
\end{proof}

\bigskip

\textbf{Remark -- }Starting from Theorem \ref{T : Main 3}, one can prove an
analogous of Corollary \ref{C : corNest} (for nested resolutions) and
Corollary \ref{C : corDet} (for non random $\phi \left( \lambda \right) $).
Moreover, Theorem \ref{T : Main 3} can be immediately extended to sequences
of random variables of the type $F_{n}=\sum_{d=1}^{\infty
}I_{M_{n}}^{d}\left( f_{n}^{\left( d\right) }\right) $, $n\geq 1$, by using
the last part of Proposition \ref{P : measClOc} (just replace $h_{\pi
^{\left( n\right) }}\left( f_{d_{n}}^{\left( n\right) }\right) $ with $%
h_{\pi ^{\left( n\right) }}\left( F_{n}\right) $).

\bigskip

Condition (\ref{CO+I}) can be difficult to verify, since it involves the
sequence of random integrands $h_{\pi ^{\left( n\right) }}\left(
f_{d_{n}}^{\left( n\right) }\right) $, which may be complex functionals of
the kernels $f_{d_{n}}^{\left( n\right) }$. In the next section, we will
show that, in the specific framework of double Poisson integrals, one can
establish neat sufficient conditions for (\ref{CO+I}), with a deterministic $%
\phi \left( \lambda \right) $, by using a version of the multiplication
formula for multiple stochastic integrals. The techniques developed below
can be extended to integrals of higher orders, to a random $\phi \left(
\lambda \right) $, and even to non-Poissonian random measures, as long as a
version of the multiplication formula is available (one might use, for
instance, the general theory of \textquotedblleft diagonal
measures\textquotedblright\ developed in \cite{Engel} and \cite{RoWa}).
These extensions will be discussed in a separate paper.

\subsection{Application: CLTs for double Poisson integrals}

In this section $\left( Z,\mathcal{Z},\mu \right) $ is a Borel measure
space, with $\mu $ non-atomic, $\sigma $-finite and positive. Also, $%
\widehat{N}$ stands for a compensated Poisson random measure on $\left( Z,%
\mathcal{Z}\right) $ with control $\mu .$ This means that $\widehat{N}%
=\left\{ \widehat{N}\left( B\right) :B\in \mathcal{Z}_{\mu }\right\} $ is an
independently scattered random measure as in Definition E, such that, for
every $B\in \mathcal{Z}_{\mu }$%
\begin{equation*}
\widehat{N}\left( B\right) \overset{\text{law}}{=}N\left( B\right) -\mu
\left( B\right) \text{,}
\end{equation*}%
where $N\left( B\right) $ is a Poisson random variable with parameter $\mu
\left( B\right) $. Note that, for every $h\in L^{2}\left( Z,\mathcal{Z},\mu
\right) =\mathfrak{H}_{\mu }$,
\begin{equation*}
X_{\widehat{N}}\left( h\right) =\int_{Z}h\left( z\right) \widehat{N}\left(
dz\right) \text{,}
\end{equation*}%
where $X_{\widehat{N}}$ is defined by (\ref{isomeasure}). Moreover, for
every $h\in \mathfrak{H}_{\mu }$ and $\lambda \in \mathbb{R}$, the L\'{e}%
vy-Khinchine exponent $\psi _{\mathfrak{H}_{\mu }}\left( h,\lambda \right) $
appearing in (\ref{LK3}), is such that (see again \cite[Proposition 19.5]%
{Sato})%
\begin{equation}
\psi _{\mathfrak{H}_{\mu }}\left( h,\lambda \right) =\int_{Z}\exp \left(
i\lambda h\left( z\right) -1-i\lambda h\left( z\right) \right) \mu \left(
dz\right)  \label{PLK}
\end{equation}%
(recall that this corresponds to the case $\rho _{\mu }\left( z,dx\right)
=\delta _{1}\left( dx\right) $ in Proposition \ref{P : LK}).

\bigskip

As an application of the previous theory, we shall study the asymptotic
behavior of a sequence of random variables of the type
\begin{equation}
F_{n}=I_{2}^{\widehat{N}}\left( f_{n}\right) \text{, \ \ }n\geq 1\text{,}
\label{seq}
\end{equation}%
where $f_{n}\in L_{s,0}^{2}\left( \mu ^{2}\right) $. In particular, we want
to use Theorem 15 to establish sufficient conditions, ensuring that $F_{n}$
converges in law to a standard Gaussian distribution. We will suppose the
following:

\bigskip

\textbf{Assumption N -- }(\textbf{N}$_{1}$) The sequence $f_{n}$, $n\geq 1$,
in (\ref{seq}) verifies:

\begin{description}
\item[(\textbf{N}$_{1}$\textbf{-i})] (\textit{Integrability condition}) $%
\forall n\geq 1$,
\begin{equation}
\int_{Z}f_{n}\left( z,\cdot \right) ^{2}\mu \left( dz\right) \in L^{2}\left(
\mu \right) ;  \label{N-1}
\end{equation}

\item[(\textbf{N}$_{1}$\textbf{-ii})] (\textit{Normalization condition}) As $%
n\rightarrow +\infty $,
\begin{equation}
2\int_{Z}\int_{Z}f_{n}\left( z,z^{\prime }\right) ^{2}\mu \left( dz\right)
\mu \left( dz^{\prime }\right) \rightarrow 1;  \label{N0}
\end{equation}

\item[(\textbf{N}$_{1}$\textbf{-iii})] (\textit{Fourth moment condition}) As
$n\rightarrow +\infty $,
\begin{equation}
\int_{Z}\int_{Z}f_{n}\left( z,z^{\prime }\right) ^{4}\mu \left( dz\right)
\mu \left( dz^{\prime }\right) \rightarrow 0  \label{N1}
\end{equation}%
(this implies, in particular, that $f_{n}\in L^{4}\left( \mu ^{2}\right) $).
\end{description}

(\textbf{N}$_{2}$) For every $n\geq 1$, there exists a collection $\left\{
Z_{n,t}:t\in \left[ 0,1\right] \right\} \subset \mathcal{Z}$, such that $%
Z_{n,0}=\varnothing $, $Z_{n,1}=Z$, $Z_{n,s}\subseteq Z_{n,t}$ for $s<t$,
and satisfying condition (\ref{conti}) (with $Z_{n,t}$ substituting $Z_{t}$%
). Note that one can take $Z_{n,t}=Z_{1,t}$ for every $n$ and $t$.

\bigskip

\textbf{Remarks -- }(1) Suppose there exists a set $B$, independent of $n$,
such that $\mu \left( B\right) <+\infty $, and, for each $n$, $f_{n}=f_{n}%
\mathbf{1}_{B}$, a.e.--$d\mu ^{2}$ (this is true, in particular, when $\mu $
is finite). Then, by the Cauchy-Schwarz inequality, if (\ref{N1}) is
verified $\left( f_{n}\right) $ must necessarily converge to zero in $%
L_{s,0}^{2}\left( \mu ^{2}\right) .$ To get more general sequences $\left(
f_{n}\right) $ we need to suppose $\mu \left( Z\right) =+\infty $.

(2) Assumption N is satisfied by a properly normalized sequence of uniformly
bounded functions, with supports \textquotedblleft slowly converging to $Z$%
\textquotedblright . For instance, consider a sequence $g_{n}\in
L_{s,0}^{2}\left( \mu ^{2}\right) $ such that, for $n\geq 1$, $\left\vert
g_{n}\left( \cdot ,\cdot \cdot \right) \right\vert \leq c<+\infty $ ($c$
independent of $n$) and the support of $g_{n}$ is contained in a set of the
type $B_{n}\times B_{n}$, where $0<\mu \left( B_{n}\right) <+\infty $ and $%
\mu \left( B_{n}\right) \rightarrow +\infty $. Then, if
\begin{equation*}
\mu \left( B_{n}\right) ^{-2}\int_{Z}\int_{Z}g_{n}\left( z,z^{\prime
}\right) ^{2}\mu \left( dz\right) \mu \left( dz^{\prime }\right) \rightarrow
1\text{,}
\end{equation*}%
the sequence $f_{n}\triangleq \mu \left( B_{n}\right) ^{-1}g_{n}$, $n\geq 1$%
, verifies Assumption N. Indeed, since $\left\vert f_{n}\right\vert \leq
c\mu \left( B_{n}\right) ^{-1}$,
\begin{eqnarray*}
\int_{Z}\int_{Z}f_{n}\left( z,z^{\prime }\right) ^{4}\mu \left( dz\right)
\mu \left( dz^{\prime }\right) &\leq &\frac{c^{4}}{\mu \left( B_{n}\right)
^{2}}\rightarrow 0 \\
\int_{Z}\left( \int_{Z}f_{n}\left( z,z^{\prime }\right) ^{2}\mu \left(
dz\right) \right) ^{2}\mu \left( dz^{\prime }\right) &\leq &\frac{c^{4}}{\mu
\left( B_{n}\right) }<+\infty .
\end{eqnarray*}

\bigskip

Before stating the main result of the section, we recall a useful version of
the \textit{multiplication formula }for multiple Poisson integrals. To this
end, we define, for $q,p\geq 1$, $f\in L_{s,0}^{2}\left( \mu ^{p}\right) $, $%
g\in L_{s,0}^{2}\left( \mu ^{q}\right) $, $r=0,...,q\wedge p$ and $l=1,...,r$%
, the (contraction) kernel on $Z^{p+q-r-l}$, which reduces the number of
variables in the product $fg$ from $p+q$ to $p+q-r-l$ as follows: $r$
variables are identified and, among these, $l$ are integrated out. This
contraction kernel is formally defined as follows:
\begin{eqnarray*}
&&f\star _{r}^{l}g(\gamma _{1},\ldots ,\gamma _{r-l},t_{1},\ldots
,t_{p-r},s_{1},\ldots ,s_{q-r}) \\
&=&\int_{Z^{l}}f(z_{1},\ldots ,z_{l},\gamma _{1},\ldots ,\gamma
_{r-l},t_{1},\ldots ,t_{p-r})g(z_{1},\ldots ,z_{l},\gamma _{1},\ldots
,\gamma _{r-l},s_{1},\ldots ,s_{q-r})\mu ^{l}\left( dz_{1}...dz_{l}\right)
\text{,}
\end{eqnarray*}%
and, for $l=0$,
\begin{equation}
f\star _{r}^{0}g(\gamma _{1},\ldots ,\gamma _{r},t_{1},\ldots
,t_{p-r},s_{1},\ldots ,s_{q-r})=f(\gamma _{1},\ldots ,\gamma
_{r},t_{1},\ldots ,t_{p-r})g(\gamma _{1},\ldots ,\gamma _{r},s_{1},\ldots
,s_{q-r}),  \label{Ocont}
\end{equation}%
so that $f\star _{0}^{0}g(t_{1},\ldots ,t_{p},s_{1},\ldots
,s_{q})=f(t_{1},\ldots ,t_{p})g(s_{1},\ldots ,s_{q})$. For example, if $%
p=q=2 $,
\begin{eqnarray}
f\star _{1}^{0}g\left( \gamma ,t,s\right) &=&f\left( \gamma ,t\right)
g\left( \gamma ,s\right) \text{, \ \ }f\star _{1}^{1}g\left( t,s\right)
=\int_{Z}f\left( z,t\right) g\left( z,s\right) \mu \left( dz\right)
\label{exCont1} \\
f\star _{2}^{1}g\left( \gamma \right) &=&\int_{Z}f\left( z,\gamma \right)
g\left( z,\gamma \right) \mu \left( dz\right) \text{, \ \ }f\star
_{2}^{2}g=\int_{Z}\int_{Z}f\left( z_{1},z_{2}\right) g\left(
z_{1},z_{2}\right) \mu \left( dz_{1}\right) \mu \left( dz_{2}\right) .
\label{exCont2}
\end{eqnarray}

The following product formula for two Poisson multiple integrals is proved
e.g. in \cite{Kab}: let $f\in L_{s,0}^{2}\left( \mu ^{p}\right) $ and $g\in
L_{s,0}^{2}\left( \mu ^{q}\right) $, $p,q\geq 1$, and suppose moreover that $%
f\star _{r}^{l}g\in L^{2}(\mu ^{p+q-r-l})$ for every $r=0,...,p\wedge q$ and
$l=1,...,r$, then
\begin{equation}
I_{p}^{\widehat{N}}(f)I_{q}^{\widehat{N}}(g)=\sum_{r=0}^{p\wedge q}r!\dbinom{%
p}{r}\dbinom{q}{r}\sum_{l=0}^{r}I_{q+p-r-l}^{\widehat{N}}(\widetilde{f\star
_{r}^{l}g}),  \label{pproduct}
\end{equation}%
where the tilde ( $\widetilde{}$ ) stands for symmetrization (note that $%
\widetilde{f\star _{r}^{l}g}$ need not vanish on diagonals, and that we use
convention (\ref{conv})).

\bigskip

The next result is the announced central limit theorem.

\bigskip

\begin{theorem}
\label{T : PoissCLT}Define the sequence $F_{n}$ and $f_{n}\in
L_{s,0}^{2}(\mu ^{2})$, $n\geq 1$, as in (\ref{seq}), and suppose Assumption
N holds. Then, $f_{n}\star _{1}^{0}f_{n}\in L^{2}(\mu ^{3})$ and $f_{n}\star
_{1}^{1}f_{n}\in L_{s}^{2}(\mu ^{2})$ for every $n\geq 1$, and moreover:

\begin{enumerate}
\item if
\begin{equation}
f_{n}\star _{1}^{1}f_{n}\rightarrow 0\text{ in }L_{s,0}^{2}(\mu ^{2})\text{
\ and \ }f_{n}\star _{2}^{1}f_{n}\rightarrow 0\text{ in }L^{2}(\mu )
\label{G*}
\end{equation}%
then
\begin{equation}
F_{n}\overset{\text{law}}{\rightarrow }N\left( 0,1\right) \text{,}
\label{GG}
\end{equation}%
where $N\left( 0,1\right) $ is a standard Gaussian random variable;

\item if $F_{n}\in L^{4}\left( \mathbb{P}\right) $ for every $n$, then a
sufficient condition to have (\ref{G*}) is that%
\begin{equation}
\mathbb{E}\left( F_{n}^{4}\right) \rightarrow 3;  \label{GGG}
\end{equation}

\item if the sequence $\left\{ F_{n}^{4}:n\geq 1\right\} $ is uniformly
integrable, then conditions (\ref{G*}), (\ref{GG}) and (\ref{GGG}) are
equivalent.
\end{enumerate}
\end{theorem}

\bigskip

\textbf{Remarks }-- (a) Note that the statement of Theorem \ref{T : PoissCLT}
does not involve any resolution of the identity. However, part (\textbf{N}$%
_{2}$) of Assumption N will play a crucial role in the proof.

(b) Observe that
\begin{eqnarray}
\left\Vert f_{n}\star _{1}^{1}f_{n}\right\Vert _{L^{2}\left( \mu ^{2}\right)
}^{2} &=&\int_{Z}\left( \int_{Z}f\left( a,z\right) f\left( b,z\right) \mu
\left( dz\right) \right) ^{2}\mu \left( da\right) \mu \left( db\right)
\label{norm1} \\
\left\Vert f_{n}\star _{2}^{1}f_{n}\right\Vert _{L^{2}\left( \mu \right)
}^{2} &=&\int_{Z}\left( \int_{Z}f\left( a,z\right) ^{2}\mu \left( da\right)
\right) ^{2}\mu \left( dz\right) .  \label{norm2}
\end{eqnarray}

(c) Let $G$ be a Gaussian measure on $\left( Z,\mathcal{Z}\right) $, with
control $\mu $, and, for $n\geq 1$, let $H_{n}=I_{2}^{G}\left( h_{n}\right) $
be the double Wiener-It\^{o} integral of a function $h_{n}\in
L_{s,0}^{2}(\mu ^{2})$. In \cite[Theorem 1]{NuPe} it is proved that, if $%
2\left\Vert h_{n}\right\Vert ^{2}\rightarrow 1$ and \textit{regardless of
Assumption N}, the following three conditions are equivalent: (i) $H_{n}%
\overset{\text{law}}{\rightarrow }N\left( 0,1\right) $, (ii) $\mathbb{E}%
\left( H_{n}^{4}\right) \rightarrow 3$, (iii) $h_{n}\star
_{1}^{1}h_{n}\rightarrow 0$. Note also that Theorem 1\ in \cite{NuPe}
applies to multiple integrals of arbitrary order.

(d) A sufficient condition for the uniform integrability of $\left(
F_{n}^{4}\right) $ is clearly that $\sup_{n}\mathbb{E}\left(
F_{n}^{4+\varepsilon }\right) <+\infty $ for some $\varepsilon >0$. Note
that in the Gaussian framework of \cite[Theorem 1]{NuPe} the uniform
integrability condition is always satisfied. Indeed, by noting $%
H_{n}=I_{2}^{G}\left( h_{n}\right) $ ($n\geq 1$) the sequence of double
integrals introduced in the previous remark, for every $p>2$ there exists a
finite constant $c_{p}$ such that $\sup_{n}\mathbb{E}\left( \left\vert
H_{n}\right\vert ^{p}\right) \leq c_{p}\sup_{n}\mathbb{E}\left(
H_{n}^{2}\right) <+\infty $, where the last relation follows from the
normalization condition $\mathbb{E}\left( H_{n}^{2}\right) =2\left\Vert
h_{n}\right\Vert ^{2}\rightarrow 1$.

\bigskip

\textbf{Proof of Theorem \ref{T : PoissCLT}. }Since
\begin{equation*}
f_{n}\star _{1}^{1}f_{n}\left( t,s\right) =\int_{Z}f_{n}\left( s,z\right)
f_{n}\left( t,z\right) \mu \left( dz\right) \text{,}
\end{equation*}%
and $f\in L^{2}\left( \mu ^{2}\right) $, the relation $f_{n}\star
_{1}^{1}f_{n}\in L_{s}^{2}(\mu ^{2})$ is a consequence of the Cauchy-Schwarz
inequality. On the other hand, by (\ref{exCont1}),
\begin{eqnarray*}
\int_{Z^{3}}\left( f_{n}\star _{1}^{0}f_{n}\left( \gamma ,t,s\right) \right)
^{2}\mu ^{3}\left( d\gamma ,dt,ds\right) &=&\int_{Z}\left(
\int_{Z}f_{n}\left( \gamma ,t\right) ^{2}\mu \left( dt\right) \times
\int_{Z}f_{n}\left( \gamma ,s\right) ^{2}\mu \left( ds\right) \right) \mu
\left( d\gamma \right) \\
&=&\int_{Z}\left( \int_{Z}f_{n}\left( \gamma ,s\right) ^{2}\mu \left(
ds\right) \right) ^{2}\mu \left( d\gamma \right) <+\infty ,
\end{eqnarray*}%
due to part (\textbf{N}$_{1}$\textbf{-i}) in Assumption N, so that $%
f_{n}\star _{1}^{0}f_{n}\in L^{2}(\mu ^{3})$.

(\textit{Proof of point 1})\textbf{\ }We shall apply Theorem \ref{T : Main 3}
in the following case: $\phi \left( \lambda \right) =\exp \left( -\lambda
^{2}/2\right) $, $\lambda \in \mathbb{R}$, and, for $n\geq 1$, $\left( Z_{n},%
\mathcal{Z}_{n}\right) =\left( Z,\mathcal{Z}\right) $, $d_{n}=2$, $M_{n}=%
\widehat{N}$, $f_{d_{n}}^{\left( n\right) }=f_{2}^{\left( n\right) }=f_{n}$,
$\mu _{n}=\mu $, $t_{n}=0$ and $\pi ^{\left( n\right) }h\left( z\right) =%
\mathbf{1}_{Z_{n,t}}\left( z\right) h\left( z\right) $, $\forall h\in
\mathfrak{H}_{\mu }$ ($=L^{2}\left( Z,\mathcal{Z},\mu \right) $), where the
sets $Z_{n,t}$ are defined as in part (\textbf{N}$_{2}$) of Assumption N. In
this case, for $n\geq 1$, $Z_{n,t_{n}}=Z_{n,0}=\varnothing $ by definition,
and therefore $\mathcal{F}_{t_{n}}^{\pi ^{\left( n\right) }}=\left\{
\varnothing ,\Omega \right\} $, so that assumptions (\ref{CO+O}) and (\ref%
{CO+B}) are immaterial. Moreover, by using (\ref{PLK}), for every $h\in
\mathfrak{H}_{\mu }$%
\begin{eqnarray}
\int_{Z_{n}}K_{\mu _{n}}\left( \lambda h\left( z\right) ,z\right) \mu
_{n}\left( dz\right) &=&\int_{Z}K_{\mu }\left( \lambda h\left( z\right)
,z\right) \mu \left( dz\right)  \label{expl} \\
&=&\int_{Z}\left( \exp \left( i\lambda h\left( z\right) \right) -1-i\lambda
h\left( z\right) \right) \mu \left( dz\right) \text{,}  \notag
\end{eqnarray}%
where $K_{\mu }$ is defined by (\ref{innerLK}). Now define, for $n\geq 1$
and $z\in Z$,
\begin{equation*}
h_{\pi ^{\left( n\right) }}\left( f_{n}\right) \left( z\right) =2I_{1}^{%
\widehat{N}}\left( f_{n}\left( z,\cdot \right) \mathbf{1}\left( \cdot \prec
_{\pi ^{\left( n\right) }}z\right) \right) \text{,}
\end{equation*}%
where the notation is the same as in (\ref{projected K}). According to
Theorem \ref{T : Main 3} and (\ref{expl}), to prove Point 1 it is sufficient
to show that, when Assumption N is verified, condition (\ref{G*}) implies
that, as $n\rightarrow +\infty $,%
\begin{equation*}
\int_{Z}\left( \exp \left( i\lambda h_{\pi ^{\left( n\right) }}\left(
f_{n}\right) \left( z\right) \right) -1-i\lambda h_{\pi ^{\left( n\right)
}}\left( f_{n}\right) \left( z\right) \right) \mu \left( dz\right) \overset{%
\mathbb{P}}{\rightarrow }-\frac{1}{2}.
\end{equation*}%
To this end, we write
\begin{eqnarray}
&&\int_{Z}\left( \exp \left( i\lambda h_{\pi ^{\left( n\right) }}\left(
f_{n}\right) \left( z\right) \right) -1-i\lambda h_{\pi ^{\left( n\right)
}}\left( f_{n}\right) \left( z\right) \right) \mu \left( dz\right)
\label{dff} \\
&=&-\frac{1}{2}\int_{Z}\left( h_{\pi ^{\left( n\right) }}\left( f_{n}\right)
\left( z\right) \right) ^{2}\mu \left( dz\right)  \notag \\
&&+\int_{Z}\left( \exp \left( i\lambda h_{\pi ^{\left( n\right) }}\left(
f_{n}\right) \left( z\right) \right) -1-i\lambda h_{\pi ^{\left( n\right)
}}\left( f_{n}\right) \left( z\right) +\frac{1}{2}\left( h_{\pi ^{\left(
n\right) }}\left( f_{n}\right) \left( z\right) \right) ^{2}\right) \mu
\left( dz\right)  \notag \\
&\triangleq &U_{n}+V_{n}\text{,}  \notag
\end{eqnarray}%
and we shall show that, under the assumptions of Theorem \ref{T : PoissCLT},
$U_{n}\overset{\mathbb{P}}{\rightarrow }-\frac{1}{2}$ and $V_{n}\overset{%
\mathbb{P}}{\rightarrow }0.$ We now apply (\ref{pproduct}) in the case $%
p=q=1 $, to have
\begin{eqnarray}
&&\left( h_{\pi ^{\left( n\right) }}\left( f_{n}\right) \left( z\right)
\right) ^{2}=4I_{1}^{\widehat{N}}\left( f_{n}\left( z,\cdot \right) \mathbf{1%
}\left( \cdot \prec _{\pi ^{\left( n\right) }}z\right) \right) ^{2}  \notag
\\
&\triangleq &4\int_{Z}f_{n}\left( z,x\right) ^{2}\mathbf{1}\left( x\prec
_{\pi ^{\left( n\right) }}z\right) \mu \left( dx\right) +4I_{1}^{\widehat{N}%
}\left( f_{n}\left( z,\cdot \right) ^{2}\mathbf{1}\left( \cdot \prec _{\pi
^{\left( n\right) }}z\right) \right) +4I_{2}^{\widehat{N}}\left( g_{n}\left(
z;\cdot ,\cdot \cdot \right) \right)  \label{cc}
\end{eqnarray}%
where $g_{n}\left( z;\cdot ,\cdot \cdot \right) \in L_{s,0}^{2}\left( \mu
^{2}\right) $ is given by%
\begin{eqnarray}
g_{n}\left( z;a,b\right) &=&f_{n}\star _{0}^{0}f_{n}\left( z,a;z,b\right)
\mathbf{1}_{\left( a\prec _{\pi ^{\left( n\right) }}z\right) }\mathbf{1}%
_{\left( b\prec _{\pi ^{\left( n\right) }}z\right) }  \label{c0} \\
&=&f_{n}\left( z,a\right) f_{n}\left( z,b\right) \mathbf{1}_{\left( a\prec
_{\pi ^{\left( n\right) }}z\right) }\mathbf{1}_{\left( b\prec _{\pi ^{\left(
n\right) }}z\right) }.  \notag
\end{eqnarray}%
We deal with each term in (\ref{cc}) in succession. For the first term
observe that, due to Corollary \ref{C : diagSets} and the symmetry of $f_{n}$%
,%
\begin{eqnarray*}
&&\int_{Z}\int_{Z}f_{n}\left( z,x\right) ^{2}\mu \left( dx\right) \mu \left(
dz\right) \\
&=&\int_{Z}\int_{Z}f_{n}\left( z,x\right) ^{2}\left[ \mathbf{1}\left( x\prec
_{\pi ^{\left( n\right) }}z\right) +\mathbf{1}\left( z\prec _{\pi ^{\left(
n\right) }}x\right) \right] \mu \left( dx\right) \mu \left( dz\right) \\
&=&2\int_{Z}\int_{Z}f_{n}\left( z,x\right) ^{2}\mathbf{1}\left( x\prec _{\pi
^{\left( n\right) }}z\right) \mu \left( dx\right) \mu \left( dz\right)
\end{eqnarray*}%
and therefore, thanks to Assumption N,%
\begin{equation}
-\frac{1}{2}\int_{Z}4\int_{Z}f_{n}\left( z,x\right) ^{2}\mathbf{1}\left(
x\prec _{\pi ^{\left( n\right) }}z\right) \mu \left( dx\right) \mu \left(
dz\right) =-\frac{1}{2}\left[ 2\int_{Z}\int_{Z}f_{n}\left( z,x\right)
^{2}\mu \left( dx\right) \mu \left( dz\right) \right] \rightarrow -\frac{1}{2%
}.  \label{c}
\end{equation}%
For the second term in (\ref{cc}) one has
\begin{eqnarray}
\mathbb{E}\left[ \int_{Z}I_{1}^{\widehat{N}}\left( f_{n}\left( z,\cdot
\right) ^{2}\mathbf{1}\left( \cdot \prec _{\pi ^{\left( n\right) }}z\right)
\right) \mu \left( dz\right) \right] ^{2} &=&\mathbb{E}\left[ I_{1}^{%
\widehat{N}}\left( \int_{Z}f_{n}\left( z,\cdot \right) ^{2}\mathbf{1}\left(
\cdot \prec _{\pi ^{\left( n\right) }}z\right) \mu \left( dz\right) \right) %
\right] ^{2}  \notag \\
&=&\int_{Z}\left[ \int_{Z}f_{n}\left( z,x\right) ^{2}\mathbf{1}\left( x\prec
_{\pi ^{\left( n\right) }}z\right) \mu \left( dz\right) \right] ^{2}\mu
\left( dx\right)  \notag \\
&\leq &\int_{Z}\left[ \int_{Z}f_{n}\left( z,x\right) ^{2}\mu \left(
dz\right) \right] ^{2}\mu \left( dx\right) \rightarrow 0,  \label{ccc}
\end{eqnarray}%
due to (\ref{norm1}) and (\ref{G*}). Now consider the third term in (\ref{cc}%
), and observe that, by a Fubini argument,
\begin{equation*}
\int_{Z}I_{2}^{\widehat{N}}\left( g_{n}\left( z;\cdot ,\cdot \cdot \right)
\right) \mu \left( dz\right) =I_{2}^{\widehat{N}}\left( h_{n}\right) ,
\end{equation*}%
where, thanks to (\ref{c0}), $h_{n}\in L_{s,0}^{2}\left( \mu ^{2}\right) $
is s.t.%
\begin{equation}
h_{n}\left( a,b\right) =\int_{Z}f_{n}\left( z,a\right) f_{n}\left(
z,b\right) \mathbf{1}_{\left( a\prec _{\pi ^{\left( n\right) }}z\right) }%
\mathbf{1}_{\left( b\prec _{\pi ^{\left( n\right) }}z\right) }\mu \left(
dz\right) .  \label{df}
\end{equation}%
We now want to show that $f_{n}\star _{1}^{1}f_{n}\left( a,b\right)
=\int_{Z}f_{n}\left( z,a\right) f_{n}\left( z,b\right) \mu \left( dz\right)
\rightarrow 0$ implies that $h_{n}\rightarrow 0$. To do this, we start by
observing that, a.e.-$\mu ^{2}\left( da,db\right) $ and thanks to Corollary %
\ref{C : diagSets},%
\begin{equation*}
f_{n}\left( a,b\right) =f_{n}\left( a,b\right) \mathbf{1}_{\left( a\prec
_{\pi ^{\left( n\right) }}b\right) \cup \left( b\prec _{\pi ^{\left(
n\right) }}a\right) }.
\end{equation*}%
As a consequence, by noting (for fixed $z$)
\begin{eqnarray*}
\left( z\prec _{\pi ^{\left( n\right) }}a\vee b\right) &=&\left[ \left(
z\prec _{\pi ^{\left( n\right) }}a\right) \cap \left( z\prec _{\pi ^{\left(
n\right) }}b\right) \right] \cup \left( a\prec _{\pi ^{\left( n\right)
}}z\prec _{\pi ^{\left( n\right) }}b\right) \cup \left( b\prec _{\pi
^{\left( n\right) }}z\prec _{\pi ^{\left( n\right) }}a\right) \\
\left( a\vee b\prec _{\pi ^{\left( n\right) }}z\right) &=&\left( a\prec
_{\pi ^{\left( n\right) }}z\right) \cap \left( b\prec _{\pi ^{\left(
n\right) }}z\right) \text{,}
\end{eqnarray*}%
we obtain that
\begin{eqnarray}
&&\int_{Z^{2}}\left( f_{n}\star _{1}^{1}f_{n}\left( a,b\right) \right)
^{2}\mu ^{2}\left( da,db\right)  \notag \\
&=&\int_{Z^{2}}\left( \int_{Z}f_{n}\left( z,a\right) \mathbf{1}_{\left(
a\prec _{\pi ^{\left( n\right) }}z\right) \cup \left( z\prec _{\pi ^{\left(
n\right) }}a\right) }\mathbf{1}_{\left( z\prec _{\pi ^{\left( n\right)
}}b\right) \cup \left( b\prec _{\pi ^{\left( n\right) }}z\right)
}f_{n}\left( z,b\right) \mu \left( dz\right) \right) ^{2}\mu ^{2}\left(
da,db\right)  \notag \\
&=&\int_{Z^{2}}\left( \int_{Z}f_{n}\left( z,a\right) f_{n}\left( z,b\right)
\left( \mathbf{1}_{\left( z\prec _{\pi ^{\left( n\right) }}a\vee b\right) }+%
\mathbf{1}_{\left( a\vee b\prec _{\pi ^{\left( n\right) }}z\right) }\right)
\mu \left( dz\right) \right) ^{2}\mu ^{2}\left( da,db\right)  \label{vis0} \\
&=&\int_{Z^{2}}\left( \int_{Z}f_{n}\left( z,a\right) f_{n}\left( z,b\right)
\mathbf{1}_{\left( z\prec _{\pi ^{\left( n\right) }}a\vee b\right) }\mu
\left( dz\right) \right) ^{2}\mu ^{2}\left( da,db\right)  \label{vis1} \\
&&+\int_{Z^{2}}\left( \int_{Z}f_{n}\left( z,a\right) f_{n}\left( z,b\right)
\mathbf{1}_{\left( a\vee b\prec _{\pi ^{\left( n\right) }}z\right) }\mu
\left( dz\right) \right) ^{2}\mu ^{2}\left( da,db\right)  \notag \\
&&+2\int_{Z^{2}}\left( \int_{Z}f_{n}\left( z,a\right) f_{n}\left( z,b\right)
\mathbf{1}_{\left( a\vee b\prec _{\pi ^{\left( n\right) }}z\right) }\mu
\left( dz\right) \right) \times  \notag \\
&&\times \left( \int_{Z}f_{n}\left( z^{\prime },a\right) f_{n}\left(
z^{\prime },b\right) \mathbf{1}_{\left( z\prec _{\pi ^{\left( n\right)
}}a\vee b\right) }\mu \left( dz^{\prime }\right) \right) \mu ^{2}\left(
da,db\right) .  \notag
\end{eqnarray}%
Now we note $\left( a\prec _{\pi ^{\left( n\right) }}z\wedge z^{\prime
}\right) =\left( a\prec _{\pi ^{\left( n\right) }}z^{\prime }\right) \cap
\left( a\prec _{\pi ^{\left( n\right) }}z\right) $, so that, by a Fubini
argument,
\begin{eqnarray}
&&\int_{Z^{2}}\left( \int_{Z}f_{n}\left( z,a\right) f_{n}\left( z,b\right)
\mathbf{1}_{\left( a\vee b\prec _{\pi ^{\left( n\right) }}z\right) }\mu
\left( dz\right) \right) ^{2}\mu ^{2}\left( da,db\right)  \notag \\
&=&\int_{Z^{2}}\left( \int_{Z}f_{n}\left( z,a\right) f_{n}\left( z^{\prime
},a\right) \mathbf{1}_{\left( a\prec _{\pi ^{\left( n\right) }}z\wedge
z^{\prime }\right) }\mu \left( da\right) \right) ^{2}\mu ^{2}\left(
dz,dz^{\prime }\right)  \label{vis2}
\end{eqnarray}%
and also, with obvious notation,
\begin{eqnarray*}
&&2\int_{Z^{2}}\left( \int_{Z}f_{n}\left( z,a\right) f_{n}\left( z,b\right)
\mathbf{1}_{\left( a\vee b\prec _{\pi ^{\left( n\right) }}z\right) }\mu
\left( dz\right) \right) \times \\
&&\times \left( \int_{Z}f_{n}\left( z^{\prime },a\right) f_{n}\left(
z^{\prime },b\right) \mathbf{1}_{\left( z^{\prime }\prec _{\pi ^{\left(
n\right) }}a\vee b\right) }\mu \left( dz^{\prime }\right) \right) \mu
^{2}\left( da,db\right) \\
&=&\int_{Z^{2}}\left( \int_{Z}f_{n}\left( z,a\right) f_{n}\left( z,b\right)
\mathbf{1}_{\left( a\vee b\prec _{\pi ^{\left( n\right) }}z\vee z^{\prime
}\right) }\mu \left( dz\right) \right) \times \\
&&\times \left( \int_{Z}f_{n}\left( z^{\prime },a\right) f_{n}\left(
z^{\prime },b\right) \mathbf{1}_{\left( z^{\prime }\wedge z\prec _{\pi
^{\left( n\right) }}a\vee b\right) }\mu \left( dz^{\prime }\right) \right)
\mu ^{2}\left( da,db\right) ,
\end{eqnarray*}%
so that the relation
\begin{eqnarray*}
\mathbf{1}_{\left( a\vee b\prec _{\pi ^{\left( n\right) }}z\vee z^{\prime
}\right) }\mathbf{1}_{\left( z^{\prime }\wedge z\prec _{\pi ^{\left(
n\right) }}a\vee b\right) } &=&\mathbf{1}_{\left( z^{\prime }\wedge z\prec
_{\pi ^{\left( n\right) }}a,b\prec _{\pi ^{\left( n\right) }}z\vee z^{\prime
}\right) }+\mathbf{1}_{\left( z^{\prime }\wedge z\prec _{\pi ^{\left(
n\right) }}a\prec _{\pi ^{\left( n\right) }}z\vee z^{\prime }\right) }%
\mathbf{1}_{\left( b\prec _{\pi ^{\left( n\right) }}z\vee z^{\prime }\right)
} \\
&&+\mathbf{1}_{\left( z^{\prime }\wedge z\prec _{\pi ^{\left( n\right)
}}b\prec _{\pi ^{\left( n\right) }}z\vee z^{\prime }\right) }\mathbf{1}%
_{\left( a\prec _{\pi ^{\left( n\right) }}z\vee z^{\prime }\right) }
\end{eqnarray*}%
gives
\begin{eqnarray}
&&2\int_{Z^{2}}\left( \int_{Z}f_{n}\left( z,a\right) f_{n}\left( z,b\right)
\mathbf{1}_{\left( a\vee b\prec _{\pi ^{\left( n\right) }}z\right) }\mu
\left( dz\right) \right) \times  \notag \\
&&\times \left( \int_{Z}f_{n}\left( z^{\prime },a\right) f_{n}\left(
z^{\prime },b\right) \mathbf{1}_{\left( z^{\prime }\prec _{\pi ^{\left(
n\right) }}a\vee b\right) }\mu \left( dz^{\prime }\right) \right) \mu
^{2}\left( da,db\right)  \notag \\
&=&\int_{Z^{2}}\left( \int_{Z}f_{n}\left( z,a\right) f_{n}\left( z^{\prime
},a\right) \mathbf{1}_{\left( z\wedge z^{\prime }\prec _{\pi ^{\left(
n\right) }}a\prec _{\pi ^{\left( n\right) }}z\vee z^{\prime }\right) }\mu
\left( da\right) \right) ^{2}\mu ^{2}\left( dz,dz^{\prime }\right)
\label{vis3} \\
&&+2\int_{Z^{2}}\left( \int_{Z}f_{n}\left( z,a\right) f_{n}\left( z^{\prime
},a\right) \mathbf{1}_{\left( z\wedge z^{\prime }\prec _{\pi ^{\left(
n\right) }}a\prec _{\pi ^{\left( n\right) }}z\vee z^{\prime }\right) }\mu
\left( da\right) \right) \times  \label{vis4} \\
&&\times \left( \int_{Z}f_{n}\left( z,a\right) f_{n}\left( z^{\prime
},a\right) \mathbf{1}_{\left( a\prec _{\pi ^{\left( n\right) }}z\wedge
z^{\prime }\right) }\mu \left( da\right) \right) \mu ^{2}\left(
dz,dz^{\prime }\right) .  \notag
\end{eqnarray}%
Observe that the terms (\ref{vis1}), (\ref{vis2}), (\ref{vis3}) and (\ref%
{vis4}) are integrals of terms respectively of the form $\left( A+B\right)
^{2}$, $A^{2}$, $B^{2}$ and $2AB$, whose sum therefore equals $2\left(
A+B\right) ^{2}$, yielding%
\begin{equation}
\int_{Z^{2}}\left( f_{n}\star _{1}^{1}f_{n}\left( a,b\right) \right) ^{2}\mu
^{2}\left( da,db\right) =2\int_{Z^{2}}\left( \int_{Z}f_{n}\left( z,a\right)
f_{n}\left( z,b\right) \mathbf{1}_{\left( z\prec _{\pi ^{\left( n\right)
}}a\vee b\right) }\mu \left( dz\right) \right) ^{2}\mu ^{2}\left(
da,db\right) .  \label{vis5}
\end{equation}%
Since $h_{n}$ (as defined in (\ref{df})) is such that
\begin{eqnarray*}
&&\int_{Z^{2}}h_{n}\left( a,b\right) ^{2}\mu ^{2}\left( da,db\right) \\
&=&\int_{Z^{2}}\left( \int_{Z}f_{n}\left( z,a\right) f_{n}\left( z,b\right)
\mathbf{1}_{\left( a\vee b\prec _{\pi ^{\left( n\right) }}z\right) }\mu
\left( dz\right) \right) ^{2}\mu ^{2}\left( da,db\right)
\end{eqnarray*}%
and
\begin{eqnarray*}
&&\int_{Z^{2}}\left( f_{n}\star _{1}^{1}f_{n}\left( a,b\right) \right)
^{2}\mu ^{2}\left( da,db\right) \\
&=&\int_{Z^{2}}\left[ \int_{Z}f_{n}\left( z,a\right) f_{n}\left( z,b\right)
\left( \mathbf{1}_{\left( a\vee b\prec _{\pi ^{\left( n\right) }}z\right) }+%
\mathbf{1}_{\left( z\prec _{\pi ^{\left( n\right) }}a\vee b\right) }\right)
\mu \left( dz\right) \right] ^{2}\mu ^{2}\left( da,db\right) \text{,}
\end{eqnarray*}%
relation (\ref{vis5}) gives the implication: if $f_{n}\star
_{1}^{1}f_{n}\rightarrow 0$ in $L^{2}\left( \mu ^{2}\right) $, then $%
h_{n}\rightarrow 0$ in $L^{2}\left( \mu ^{2}\right) $. This last result,
combined with (\ref{c}) and (\ref{ccc}) implies that the sequence $U_{n}$, $%
n\geq 1$, as defined in (\ref{dff}), converges to $-\frac{1}{2}$ in
probability.

To show that $V_{n}\overset{\mathbb{P}}{\rightarrow }0$, observe that $%
\left\vert \exp \left( i\lambda x\right) -1-i\lambda x+\frac{1}{2}%
x^{2}\right\vert \leq \left\vert x\right\vert ^{3}/6$, and consequently, by
Cauchy-Schwarz
\begin{eqnarray*}
\left\vert V_{n}\right\vert &\leq &\frac{1}{6}\int_{Z}\left\vert h_{\pi
^{\left( n\right) }}\left( f_{n}\right) \left( z\right) \right\vert ^{3}\mu
\left( dz\right) \\
&\leq &\frac{1}{6}\left( \int_{Z}\left\vert h_{\pi ^{\left( n\right)
}}\left( f_{n}\right) \left( z\right) \right\vert ^{4}\mu \left( dz\right)
\right) ^{\frac{1}{2}}\left( \int_{Z}\left\vert h_{\pi ^{\left( n\right)
}}\left( f_{n}\right) \left( z\right) \right\vert ^{2}\mu \left( dz\right)
\right) ^{\frac{1}{2}}.
\end{eqnarray*}%
Since the first part of the proof implies that, under (\ref{G*}), $\left(
\int_{Z}\left\vert h_{\pi ^{\left( n\right) }}\left( f_{n}\right) \left(
z\right) \right\vert ^{2}\mu \left( dz\right) \right) ^{\frac{1}{2}}\overset{%
\mathbb{P}}{\rightarrow }1$, to conclude the proof of point 1 it is
sufficient to show that, under Assumption N and (\ref{G*}), $%
\int_{Z}\left\vert h_{\pi ^{\left( n\right) }}\left( f_{n}\right) \left(
z\right) \right\vert ^{4}\mu \left( dz\right) \rightarrow 0$ in $L^{1}\left(
\mathbb{P}\right) $. To do this, one can use (\ref{cc}) and the
orthogonality of multiple integrals of different orders to obtain that, for
any fixed $z$,%
\begin{eqnarray*}
\mathbb{E}\left[ \left( h_{\pi ^{\left( n\right) }}\left( f_{n}\right)
\left( z\right) \right) ^{4}\right] &=&\mathbb{E}\left[ \left( \left( h_{\pi
^{\left( n\right) }}\left( f_{n}\right) \left( z\right) \right) ^{2}\right)
^{2}\right] \\
&=&16\left( \int_{Z}f_{n}\left( z,x\right) ^{2}\mathbf{1}\left( x\prec _{\pi
^{\left( n\right) }}z\right) \mu \left( dx\right) \right) ^{2} \\
&&+16\int_{Z}f_{n}\left( z,\cdot \right) ^{4}\mathbf{1}\left( \cdot \prec
_{\pi ^{\left( n\right) }}z\right) \mu \left( da\right) \\
&&+32\int_{Z^{2}}f_{n}\left( z,a\right) ^{2}f_{n}\left( z,b\right) ^{2}%
\mathbf{1}_{\left( a\prec _{\pi ^{\left( n\right) }}z\right) }\mathbf{1}%
_{\left( b\prec _{\pi ^{\left( n\right) }}z\right) }\mu ^{2}\left(
da,db\right) ,
\end{eqnarray*}%
and therefore
\begin{eqnarray*}
\mathbb{E}\int_{Z}\left\vert h_{\pi ^{\left( n\right) }}\left( f_{n}\right)
\left( z\right) \right\vert ^{4}\mu \left( dz\right) &=&\int_{Z}\mathbb{E}%
\left\vert h_{\pi ^{\left( n\right) }}\left( f_{n}\right) \left( z\right)
\right\vert ^{4}\mu \left( dz\right) \\
&\leq &16\int_{Z}\left( \int_{Z}f_{n}\left( z,x\right) ^{2}\mu \left(
dx\right) \right) ^{2}\mu \left( dz\right) \\
&&+16\int_{Z}\int_{Z}f_{n}\left( z,a\right) ^{4}\mu \left( da\right) \mu
\left( dz\right) \\
&&+32\int_{Z}\left( \int_{Z}f_{n}\left( z,x\right) ^{2}\mu \left( dx\right)
\right) ^{2}\mu \left( dz\right) \\
&\rightarrow &0\text{,}
\end{eqnarray*}%
since Assumption N and (\ref{G*}) are in order. This concludes the proof of
Point 1.

(\textit{Proof of Point 2}) To proof Point 2, use the product formula
expansion (\ref{pproduct}) (from the term with $r=0$ to the terms with $r=2$%
) to write
\begin{eqnarray*}
F_{n}^{2} &=&I_{3}^{\widehat{N}}\left( f_{n}\right) ^{2}=I_{4}^{\widehat{N}%
}\left( \widetilde{f_{n}\star _{0}^{0}f_{n}}\right) +4I_{3}^{\widehat{N}%
}\left( \widetilde{f_{n}\star _{1}^{0}f_{n}}\right) +4I_{2}^{\widehat{N}%
}\left( f_{n}\star _{1}^{1}f_{n}\right) \\
&&+2I_{2}^{\widehat{N}}\left( f_{n}\star _{2}^{0}f_{n}\right) +2I_{1}^{%
\widehat{N}}\left( f_{n}\star _{2}^{1}f_{n}\right) +2\left\Vert
f_{n}\right\Vert _{L_{s,0}^{2}\left( \mu ^{2}\right) }^{2}\text{,}
\end{eqnarray*}%
and observe that, since Assumption N holds and $f_{n}\star
_{2}^{0}f_{n}\left( a,b\right) =f_{n}\left( a,b\right) ^{2}$ (by (\ref{Ocont}%
)), $I_{2}^{\widehat{N}}\left( f_{n}\star _{2}^{0}f_{n}\right) \rightarrow 0$
in $L^{2}\left( \mathbb{P}\right) $ by (\ref{N1}), and therefore the
assumption $\mathbb{E}\left( F_{n}^{4}\right) \rightarrow 3$ implies that
\begin{eqnarray}
&&\mathbb{E}\left[ \left( F_{n}^{2}-2I_{2}^{\widehat{N}}\left( f_{n}\star
_{2}^{0}f_{n}\right) \right) ^{2}\right]  \notag \\
&=&\mathbb{E}\left[ \left( I_{4}^{\widehat{N}}\left( \widetilde{f_{n}\star
_{0}^{0}f_{n}}\right) +4I_{3}^{\widehat{N}}\left( \widetilde{f_{n}\star
_{1}^{0}f_{n}}\right) +4I_{2}^{\widehat{N}}\left( f_{n}\star
_{1}^{1}f_{n}\right) \right. \right.  \label{develop} \\
&&\left. \left. +2I_{1}^{\widehat{N}}\left( f_{n}\star _{2}^{1}f_{n}\right)
+2\left\Vert f_{n}\right\Vert _{L_{s,0}^{2}\left( \mu ^{2}\right)
}^{2}\right) ^{2}\right]  \notag \\
&\rightarrow &3.  \label{lim3}
\end{eqnarray}

Now, due to (\ref{develop}),
\begin{eqnarray}
&&\mathbb{E}\left[ \left( F_{n}^{2}-2I_{2}^{\widehat{N}}\left( f_{n}\star
_{2}^{0}f_{n}\right) \right) ^{2}\right]  \label{sqd0} \\
&=&\left[ \left( 2\left\Vert f_{n}\right\Vert _{L_{s,0}^{2}\left( \mu
^{2}\right) }^{2}\right) ^{2}+\mathbb{E}\left( 16I_{2}^{\widehat{N}}\left(
f_{n}\star _{1}^{1}f_{n}\right) ^{2}\right) +\mathbb{E}\left( I_{4}^{%
\widehat{N}}\left( \widetilde{f_{n}\star _{0}^{0}f_{n}}\right) ^{2}\right) %
\right]  \label{sqd} \\
&&+\mathbb{E}\left( 16I_{3}^{\widehat{N}}\left( \widetilde{f_{n}\star
_{1}^{0}f_{n}}\right) ^{2}\right) +\mathbb{E}\left( 4I_{1}^{\widehat{N}%
}\left( f_{n}\star _{2}^{1}f_{n}\right) ^{2}\right) .  \notag
\end{eqnarray}

There are no cross terms because the multiple integrals have different
orders, and hence are orthogonal. The most complicated term in the square
bracket in (\ref{sqd}) is $\mathbb{E}\left( I_{4}^{\widehat{N}}\left(
\widetilde{f_{n}\star _{0}^{0}f_{n}}\right) ^{2}\right) $. Since we are
dealing with second order moments, the computations are as in the Gaussian
case. We therefore obtain, using e.g. formula (2) in \cite[p. 250]{PT}, that
\begin{eqnarray*}
\mathbb{E}\left( I_{4}^{\widehat{N}}\left( \widetilde{f_{n}\star
_{0}^{0}f_{n}}\right) ^{2}\right) &=&4!\left\Vert \widetilde{f_{n}\star
_{0}^{0}f_{n}}\right\Vert _{L^{2}\left( \mu ^{2}\right) }^{2} \\
&=&2\left( 2\left\Vert f_{n}\right\Vert _{L_{s,0}^{2}\left( \mu ^{2}\right)
}^{2}\right) ^{2}+16\left\Vert f_{n}\star _{1}^{1}f_{n}\right\Vert
_{L^{2}\left( \mu ^{2}\right) }^{2}.
\end{eqnarray*}%
As a consequence, (\ref{sqd0}) equals%
\begin{eqnarray}
&&\left( 2\left\Vert f_{n}\right\Vert _{L_{s,0}^{2}\left( \mu ^{2}\right)
}^{2}\right) ^{2}+\left[ 2\left( 2\left\Vert f_{n}\right\Vert
_{L_{s,0}^{2}\left( \mu ^{2}\right) }^{2}\right) ^{2}+16\left\Vert
f_{n}\star _{1}^{1}f_{n}\right\Vert _{L^{2}\left( \mu ^{2}\right) }^{2}%
\right]  \notag \\
&&+16\times 2\left\Vert f_{n}\star _{1}^{1}f_{n}\right\Vert _{L^{2}\left(
\mu ^{2}\right) }^{2}+16\times 3!\left\Vert f_{n}\star
_{1}^{0}f_{n}\right\Vert _{L^{2}\left( \mu ^{3}\right) }^{2}+4\left\Vert
f_{n}\star _{2}^{1}f_{n}\right\Vert _{L^{2}\left( \mu \right) }^{2}  \notag
\\
&=&3\left( 2\left\Vert f_{n}\right\Vert _{L_{s,0}^{2}\left( \mu ^{2}\right)
}^{2}\right) ^{2}+48\left\Vert f_{n}\star _{1}^{1}f_{n}\right\Vert
_{L^{2}\left( \mu ^{2}\right) }^{2}+96\left\Vert f_{n}\star
_{1}^{0}f_{n}\right\Vert _{L^{2}\left( \mu ^{3}\right) }^{2}+4\left\Vert
f_{n}\star _{2}^{1}f_{n}\right\Vert _{L^{2}\left( \mu \right) }^{2}.
\label{findev}
\end{eqnarray}%
Since (\ref{findev}) converges to 3, by (\ref{lim3}), and $2\left\Vert
f_{n}\right\Vert _{L_{s,0}^{2}\left( \mu ^{2}\right) }^{2}\rightarrow 1$ by
Assumption (\textbf{N}$_{1}$\textbf{-ii}), we conclude that $\left\Vert
f_{n}\star _{1}^{1}f_{n}\right\Vert _{L^{2}\left( \mu ^{2}\right)
}^{2}\rightarrow 0$ and $\left\Vert f_{n}\star _{2}^{1}f_{n}\right\Vert
_{L^{2}\left( \mu \right) }^{2}\rightarrow 0$, thus proving Point 2.

(\textit{Proof of point 3}) If $F_{n}\overset{\text{law}}{\rightarrow }%
N\left( 0,1\right) $ and $\left( F_{n}^{4}\right) $ is uniformly integrable,
then necessarily $\mathbb{E}\left( F_{n}^{4}\right) \rightarrow \mathbb{E}%
\left( N\left( 0,1\right) ^{4}\right) =3$, so that the proof is obtained by
combining Point 1 and Point 2 in the statement. $\blacksquare $

\bigskip

\textbf{Example -- }We now exhibit an elementary example of a sequence $%
f_{n}\in L_{s,0}^{2}\left( \mu ^{2}\right) $, $n\geq 1$, verifying
conditions (\ref{N-1}), (\ref{N0}), (\ref{N1}) and (\ref{G*}). To this end,
let $B_{j}$, $j\geq 1$, be a sequence of disjoint subsets of $Z$ such that $%
\mu \left( B_{j}\right) =1$, $j\geq 1$, and set
\begin{equation*}
B_{0,j}^{2}=\left\{ \left( x,y\right) \in B_{j}\times B_{j}:x\neq y\right\}
\text{, \ \ }j\geq 1;
\end{equation*}%
note that, since $\mu $ is non-atomic, $\mu ^{2}\left( B_{0,j}^{2}\right)
=\mu ^{2}\left( B_{j}\times B_{j}\right) =1$. For $n\geq 1$ and $\left(
x,y\right) \in Z^{2}$, we define $f_{n}\left( x,y\right) =\left( 2n\right)
^{-1/2}\sum_{j=1}^{n}\mathbf{1}_{B_{0,j}^{2}}\left( x,y\right) $. Of course,
$f_{n}\in L_{s,0}^{2}\left( \mu ^{2}\right) $ by definition, and we shall
prove that $\left( f_{n}\right) $ also satisfies (\ref{N-1}), (\ref{N0}), (%
\ref{N1}) and (\ref{G*}). Indeed, $\int_{Z}f_{n}\left( z,\cdot \right)
^{2}\mu \left( dz\right) $ $=$ $2n^{-1}\sum_{j=1}^{n}\mathbf{1}%
_{B_{j}}\left( \cdot \right) \in L^{2}\left( \mu \right) $ and $2\left\Vert
f_{n}\right\Vert ^{2}=1$ by definition, so that $\left( f_{n}\right) $
verifies (\ref{N-1}) and (\ref{N0}). On the other hand,
\begin{eqnarray*}
\int_{Z}\int_{Z}f_{n}\left( x,y\right) ^{4}\mu ^{2}\left( dx,dy\right) &=&%
\frac{1}{4n^{2}}\int_{Z}\int_{Z}\left( \sum_{j=1}^{n}\mathbf{1}%
_{B_{0,j}^{2}}\left( x,y\right) \right) ^{4}\mu ^{2}\left( dx,dy\right) \\
&=&\frac{1}{4n}\rightarrow 0\text{,}
\end{eqnarray*}%
and therefore (\ref{N1}) is verified. Finally,
\begin{equation*}
\int_{Z}\left( \int_{Z}f\left( z,x\right) ^{2}\mu \left( dz\right) \right)
^{2}\mu \left( dx\right) =\frac{1}{4n}\rightarrow 0
\end{equation*}%
and%
\begin{eqnarray*}
\int_{Z}\left( \int_{Z}f\left( x,z\right) f\left( y,z\right) \mu \left(
dz\right) \right) ^{2}\mu ^{2}\left( dx,dy\right) &=&\frac{1}{4n^{2}}%
\int_{Z}\int_{Z}\left( \sum_{j=1}^{n}\mathbf{1}_{B_{0,j}^{2}}\left(
x,y\right) \right) ^{2}\mu ^{2}\left( dx,dy\right) \\
&=&\frac{1}{4n}\rightarrow 0\text{,}
\end{eqnarray*}%
thus yielding that $\left( f_{n}\right) $ satisfies (\ref{G*}), by (\ref%
{norm1}) and (\ref{norm2}). Of course, since (due e.g. to (\ref{pproduct}))%
\begin{equation*}
I_{2}^{\widehat{N}}\left( f_{n}\right) =n^{-1/2}\sum_{j=1}^{n}2^{-1/2}\left(
\widehat{N}\left( B_{j}\right) ^{2}-\widehat{N}\left( B_{j}\right) -1\right)
\text{,}
\end{equation*}%
the central limit result $I_{2}^{\widehat{N}}\left( f_{n}\right) \overset{law%
}{\rightarrow }N\left( 0,1\right) $ can be verified directly, by using a
standard version of the Central Limit Theorem, as well as the fact that $%
\widehat{N}$ is independently scattered and the $B_{j}$'s are disjoint.

\section{Stable convergence of functionals of Gaussian processes}

We shall now use Theorem \ref{T : Main} to prove general sufficient
conditions, ensuring the stable convergence of functionals of Gaussian
processes towards mixtures of normal distributions. This extends part of the
results contained in \cite{NuPe} and \cite{PT}, and leads to quite general
criteria for the stable convergence of Skorohod integrals and multiple
Wiener-It\^{o} integrals. However, to keep the lenght of this paper within
bounds, we have deferred the discussion about multiple Wiener-It\^{o}
integrals, as well as some relations with Brownian martingales to a separate
paper, see \cite{PeTaq}.

\subsection{Preliminaries}

Consider a real separable Hilbert space $\mathfrak{H}$, as well as a
continuous resolution of the identity $\pi =\left\{ \pi _{t}:t\in \left[ 0,1%
\right] \right\} \in \mathcal{R}\left( \mathfrak{H}\right) $ (see Definition
D). Throughout this paragraph, $X$ $=X\left( \mathfrak{H}\right) $ $=\left\{
X\left( f\right) :f\in \mathfrak{H}\right\} $ stands for a centered Gaussian
family, defined on some probability space $\left( \Omega ,\mathcal{F},%
\mathbb{P}\right) $, indexed by the elements of $\mathfrak{H}$ and
satisfying the isomorphic condition (\ref{basicISO}). Note, that due to the
Gaussian nature of $X$, every vector as in (\ref{indincrements}) is composed
of independent random variables, and therefore, in this case, $\mathcal{R}%
\left( \mathfrak{H}\right) =\mathcal{R}_{X}\left( \mathfrak{H}\right) $.
When (\ref{basicISO}) is satisfied and $X\left( \mathfrak{H}\right) $ is a
Gaussian family, one usually says that $X\left( \mathfrak{H}\right) $ is an
\textit{isonormal Gaussian process}, or a \textit{Gaussian measure}, over $%
\mathfrak{H}$ (see e.g. \cite[Section 1]{Nualart} or \cite{Nualart2}). As
before, we write $L^{2}\left( \mathfrak{H},X\right) $ to indicate the
(Hilbert) space of $\mathfrak{H}$-valued and $\sigma \left( X\right) $%
-measurable random variables. The filtration $\mathcal{F}^{\pi }\left(
X\right) =\left\{ \mathcal{F}_{t}^{\pi }\left( X\right) :t\in \left[ 0,1%
\right] \right\} $ (which is complete by definition) is given by formula (%
\ref{resfiltration}).

\bigskip

In what follows, we shall apply to the Gaussian measure $X$ some standard
notions and results from Malliavin calculus (the reader is again referred to
\cite{Nualart} and \cite{Nualart2} for any unexplained notation or
definition). For instance, $D=D_{X}$ and $\delta =\delta _{X}$ stand,
respectively, for the usual Malliavin derivative and Skorohod integral with
respect to the Gaussian measure $X$ (the dependence on $X$ will be dropped,
when there is no risk of confusion); for $k\geq 1$, $\mathbb{D}_{X}^{1,2}$
is the space of differentiable functionals of $X$, endowed with the norm $%
\left\Vert \cdot \right\Vert _{1,2}$ (see \cite[Chapter 1]{Nualart} for a
definition of this norm); $dom\left( \delta _{X}\right) $ is the domain of
the operator $\delta _{X}$. Note that $D_{X}$ is an operator from $\mathbb{D}%
_{X}^{1,2}$ to $L^{2}\left( \mathfrak{H},X\right) $, and also that $%
dom\left( \delta _{X}\right) \subset L^{2}\left( \mathfrak{H},X\right) $.
For every $d\geq 1$, we define $\mathfrak{H}^{\otimes d}$ and $\mathfrak{H}%
^{\odot d}$ to be, respectively, the $d$th tensor product and the $d$th
\textit{symmetric} tensor product of $\mathfrak{H}$. For $d\geq 1$ we will
denote by $I_{d}^{X}$ the isometry between $\mathfrak{H}^{\odot d}$ equipped
with the norm $\sqrt{d!}\left\Vert \cdot \right\Vert _{\mathfrak{H}^{\otimes
d}}$ and the $d$th Wiener chaos of $X$.

\bigskip

The vector spaces $L_{\pi }^{2}\left( \mathfrak{H},X\right) $ and $\mathcal{E%
}_{\pi }\left( \mathfrak{H},X\right) $, composed respectively of adapted and
elementary adapted elements of $L^{2}\left( \mathfrak{H},X\right) $, are
once again defined as in Paragraph 3. We now want to link the above defined
operators $\delta _{X}$ and $D_{X}$ to the theory developed in the previous
sections. In particular, we shall use the facts that (i) for any $\pi \in
\mathcal{R}_{X}\left( \mathfrak{H}\right) $, $L_{\pi }^{2}\left( \mathfrak{H}%
,X\right) \subseteq dom\left( \delta _{X}\right) $, and (ii) for any $u\in
L_{\pi }^{2}\left( \mathfrak{H},X\right) $ the random variable $J_{X}^{\pi
}\left( u\right) $ can be regarded as a Skorohod integral. They are based on
the following (simple) result, proved for instance in \cite[Lemme 1]{Wu}.

\bigskip

\begin{proposition}
\label{P : SKO}Let the assumptions of this section prevail. Then, $L_{\pi
}^{2}\left( \mathfrak{H},X\right) \subseteq dom\left( \delta _{X}\right) $,
and for every $h_{1},h_{2}\in L_{\pi }^{2}\left( \mathfrak{H},X\right) $%
\begin{equation}
\mathbb{E}\left( \delta _{X}\left( h_{1}\right) \delta _{X}\left(
h_{2}\right) \right) =\left( h_{1},h_{2}\right) _{L_{\pi }^{2}\left(
\mathfrak{H},X\right) }.  \label{Wu-isometry}
\end{equation}%
Moreover, if $h\in \mathcal{E}_{\pi }\left( \mathfrak{H},X\right) $ has the
form $h=\sum_{i=1}^{n}h_{i}$, where $n\geq 1$, and $h_{i}\in \mathcal{E}%
_{\pi }\left( \mathfrak{H},X\right) $ is such that%
\begin{equation*}
h_{i}=\Phi _{i}\times \left( \pi _{t_{2}^{\left( i\right) }}-\pi
_{t_{1}^{\left( i\right) }}\right) f_{i}\text{, \ \ }f_{i}\in \mathfrak{H}%
\text{, \ \ }i=1,...,n,
\end{equation*}%
with $t_{2}^{\left( i\right) }>t_{1}^{\left( i\right) }$ and $\Phi _{i}$
square integrable and $\mathcal{F}_{t_{1}^{\left( i\right) }}^{\pi }\left(
X\right) $-measurable, then
\begin{equation}
\delta _{X}\left( h\right) =\sum_{i=1}^{n}\Phi _{i}\times \left[ X\left( \pi
_{t_{2}^{\left( i\right) }}f_{i}\right) -X\left( \pi _{t_{1}^{\left(
i\right) }}f_{i}\right) \right] .  \label{simpleSkorohod}
\end{equation}
\end{proposition}

\bigskip

Relation (\ref{Wu-isometry}) implies, in the terminology of \cite{Wu}, that $%
L_{\pi }^{2}\left( \mathfrak{H},X\right) $ is a closed subspace of the
\textit{isometric subset of }$dom\left( \delta _{X}\right) $, defined as the
class of those $h\in dom\left( \delta _{X}\right) $ s.t. $\mathbb{E}\left(
\delta _{X}\left( h\right) ^{2}\right) =\left\Vert h\right\Vert
_{L^{2}\left( \mathfrak{H},X\right) }^{2}$ (note that, in general, such an
isometric subset is not even a vector space; see \cite[p. 170]{Wu}).
Relation (\ref{simpleSkorohod}) applies to simple integrands $h$, but by
combining (\ref{Wu-isometry}), (\ref{simpleSkorohod}) and Proposition \ref{P
: ItoInt}, we deduce immediately that, for every $h\in L_{\pi }^{2}\left(
\mathfrak{H},X\right) $,%
\begin{equation}
\delta _{X}\left( h\right) =J_{X}^{\pi }\left( h\right) \text{, \ \ a.s.-}%
\mathbb{P}\text{.}  \label{coincidence}
\end{equation}%
where the random variable $J_{X}^{\pi }\left( h\right) $ is defined
according to Proposition \ref{P : ItoInt} and formula (\ref{elSkorohod}).
Observe that the definition of $J_{X}^{\pi }$ involves the resolution of the
identity $\pi $, whereas the definition of $\delta $ does not involve any
notion of resolution.

\bigskip

The next crucial result, which is partly a consequence of the continuity of $%
\pi $, is an abstract version of the \textit{Clark-Ocone formula} (see \cite%
{Nualart}): it is a direct corollary of \cite[Th\'{e}or\`{e}me 1, formula
(2.4) and Th\'{e}or\`{e}me 3]{Wu}, to which the reader is referred for a
detailed proof.

\bigskip

\begin{proposition}[Abstract Clark-Ocone formula; Wu, 1990]
\label{P : AbsCO}Under the above notation and assumptions (in particular, $%
\pi $ is a \textsf{continuous} resolution of the identity as in Definition
D), for every $F\in \mathbb{D}_{X}^{1,2}$,
\begin{equation}
F=\mathbb{E}\left( F\right) +\delta \left( proj\left\{ D_{X}F\mid L_{\pi
}^{2}\left( \mathfrak{H},X\right) \right\} \right) \text{,}  \label{ACO}
\end{equation}%
where $D_{X}F$ is the Malliavin derivative of $F$, and $proj\left\{ \cdot
\mid L_{\pi }^{2}\left( \mathfrak{H},X\right) \right\} $ is the orthogonal
projection operator on $L_{\pi }^{2}\left( \mathfrak{H},X\right) $.
\end{proposition}

\bigskip

\textbf{Remarks -- }(a) Note that the right-hand side of (\ref{ACO}) is well
defined since $D_{X}F\in L^{2}\left( \mathfrak{H},X\right) $ by definition,
and therefore
\begin{equation*}
proj\left\{ D_{X}F\mid L_{\pi }^{2}\left( \mathfrak{H},X\right) \right\} \in
L_{\pi }^{2}\left( \mathfrak{H},X\right) \subseteq dom\left( \delta
_{X}\right) ,
\end{equation*}%
where the last inclusion is stated in Proposition \ref{P : SKO}.

(b) Formula (\ref{ACO}) has been proved in \cite{Wu} in the context of
abstract Wiener spaces, but in the proof of (\ref{ACO}) the role of the
underlying probability space is immaterial. The extension to the framework
of isonormal Gaussian processes is therefore standard (see e.g. \cite[%
Section 1.1]{Nualart2}).

(c) Since $\mathbb{D}_{X}^{1,2}$ is dense in $L^{2}\left( \mathbb{P}\right) $
and $\delta _{X}\left( L_{\pi }^{2}\left( \mathfrak{H},X\right) \right) $ is
an isometry (due to relation (\ref{Wu-isometry})), the Clark-Ocone formula (%
\ref{ACO}) implies that every $F\in L^{2}\left( \mathbb{P},\sigma \left(
X\right) \right) $ admits a unique \textquotedblleft
predictable\textquotedblright\ representation of the form%
\begin{equation}
F=\mathbb{E}\left( F\right) +\delta _{X}\left( u\right) \text{, \ \ }u\in
L_{\pi }^{2}\left( \mathfrak{H},X\right) \text{;}  \label{Apredictable}
\end{equation}%
see also \cite[Remarque 2, p. 172]{Wu}.

(d) Since (\ref{coincidence}) holds, formula (\ref{ACO}) can be rewritten as
\begin{equation}
F=\mathbb{E}\left( F\right) +J_{X}^{\pi }\left( proj\left\{ D_{X}F\mid
L_{\pi }^{2}\left( \mathfrak{H},X\right) \right\} \right) .  \label{ACO2}
\end{equation}

\bigskip

Now consider, as before, an independent copy of $X$, noted $\widetilde{X}%
=\left\{ \widetilde{X}\left( f\right) :f\in \mathfrak{H}\right\} $, and, for
$h\in L_{\pi }^{2}\left( \mathfrak{H},X\right) $, define the random variable
$J_{\widetilde{X}}^{\pi }\left( h\right) $ according to Proposition \ref{P :
ItoInt} and (\ref{decoupSko}). The following result is an immediate
consequence of Proposition \ref{P : LK}, and characterizes $J_{\widetilde{X}%
}^{\pi }\left( h\right) $, $h\in L_{\pi }^{2}\left( \mathfrak{H},X\right) $,
as a conditionally Gaussian random variable.

\bigskip

\begin{proposition}
\label{P: QuadV}For every $h\in L_{\pi }^{2}\left( \mathfrak{H},X\right) $
and for every $\lambda \in \mathbb{%
\mathbb{R}
}$,
\begin{equation*}
\mathbb{E}\left[ \exp \left( i\lambda J_{\widetilde{X}}^{\pi }\left(
h\right) \right) \mid \sigma \left( X\right) \right] =\exp \left( -\frac{%
\lambda ^{2}}{2}\left\Vert h\right\Vert _{\mathfrak{H}}^{2}\right) .
\end{equation*}
\end{proposition}

\subsection{Stable convergence to a mixture of Gaussian distributions}

The following result, based on Theorem \ref{T : Main}, gives general
sufficient conditions for the stable convergence of Skorohod integrals to a
mixture of Gaussian distributions. In what follows, $\mathfrak{H}_{n}$, $%
n\geq 1$, is a sequence of real separable Hilbert spaces, and, for each $%
n\geq 1$, $X_{n}=X_{n}\left( \mathfrak{H}_{n}\right) =\left\{ X_{n}\left(
g\right) :g\in \mathfrak{H}_{n}\right\} $, is an isonormal Gaussian process
over $\mathfrak{H}_{n}$; for $n\geq 1$, $\widetilde{X}_{n}$ is an
independent copy of $X_{n}$ (note that $\widetilde{X}_{n}$ appears in the
proof of the next result, but not in the statement). Recall that $\mathcal{R}%
\left( \mathfrak{H}_{n}\right) $ is a class of resolutions of the identity $%
\pi $ (see Definition D), and that the Hilbert space $L_{\pi }^{2}\left(
\mathfrak{H}_{n},X_{n}\right) $ is defined after Relation (\ref{adaptation}).

\begin{theorem}
\label{T : Main 4}Suppose that the isonormal Gaussian processes $X_{n}\left(
\mathfrak{H}_{n}\right) $, $n\geq 1$, are defined on the probability space $%
\left( \Omega ,\mathcal{F},\mathbb{P}\right) $. Let, for $n\geq 1$, $\pi
^{\left( n\right) }\in \mathcal{R}\left( \mathfrak{H}_{n}\right) $ and $%
u_{n}\in L_{\pi ^{\left( n\right) }}^{2}\left( \mathfrak{H}_{n},X_{n}\right)
$. Suppose also that there exists a sequence $\left\{ t_{n}:n\geq 1\right\} $
$\subset \left[ 0,1\right] $ and $\sigma $-fields $\left\{ \mathcal{U}%
_{n}:n\geq 1\right\} $, such that
\begin{equation}
\left\Vert \pi _{t_{n}}^{\left( n\right) }u_{n}\right\Vert _{\mathfrak{H}%
_{n}}^{2}\overset{\mathbb{P}}{\rightarrow }0  \label{ASnegl}
\end{equation}%
and
\begin{equation}
\mathcal{U}_{n}\subseteq \mathcal{U}_{n+1}\cap \mathcal{F}_{t_{n}}^{\pi
^{\left( n\right) }}\left( X_{n}\right) .  \label{boxed2}
\end{equation}%
If%
\begin{equation}
\left\Vert u_{n}\right\Vert _{\mathfrak{H}}^{2}\overset{\mathbb{P}}{%
\rightarrow }Y\text{,}  \label{ASconv}
\end{equation}%
for some $Y\in L^{2}\left( \mathbb{P}\right) $ such that $Y\neq 0$, $Y\geq 0$
and $Y\in \mathcal{U}^{\ast }\triangleq \vee _{n}\mathcal{U}_{n}$, then, as $%
n\rightarrow +\infty $,
\begin{equation*}
\mathbb{E}\left[ \exp \left( i\lambda \delta _{X_{n}}\left( u_{n}\right)
\right) \mid \mathcal{F}_{t_{n}}^{\pi ^{\left( n\right) }}\left(
X_{n}\right) \right] \overset{\mathbb{P}}{\rightarrow }\exp \left( -\frac{%
\lambda ^{2}}{2}Y\right) ,\text{ \ \ }\forall \lambda \in
\mathbb{R}
\text{,}
\end{equation*}%
and%
\begin{equation*}
\delta _{X_{n}}\left( u_{n}\right) \rightarrow _{\left( s,\mathcal{U}^{\ast
}\right) }\mathbb{E}\mu \left( \cdot \right) \text{,}
\end{equation*}%
where $\mu \in \mathbf{M}$ verifies $\widehat{\mu }\left( \lambda \right)
=\exp \left( -\frac{\lambda ^{2}}{2}Y\right) $ (see (\ref{ftr}) for the
definition of $\widehat{\mu }$).
\end{theorem}

\begin{proof}
Since $\delta _{X_{n}}\left( u_{n}\right) =J_{X_{n}}^{\pi ^{\left( n\right)
}}\left( u_{n}\right) $ for every $n$, the result follows immediately from
Theorem \ref{T : Main} by observing that, due to Proposition \ref{P: QuadV},
\begin{equation*}
\mathbb{E}\left[ \exp \left( i\lambda J_{\widetilde{X}_{n}}^{\pi ^{\left(
n\right) }}\left( u_{n}\right) \right) \mid \sigma \left( X_{n}\right) %
\right] =\exp \left( -\frac{\lambda ^{2}}{2}\left\Vert u_{n}\right\Vert _{%
\mathfrak{H}}^{2}\right) \text{, }
\end{equation*}%
and therefore (\ref{ASconv}) that $\mathbb{E}\left[ \exp \left( i\lambda J_{%
\widetilde{X}_{n}}^{\pi ^{\left( n\right) }}\left( u_{n}\right) \right) \mid
\sigma \left( X_{n}\right) \right] \rightarrow \exp \left( -\lambda
^{2}Y/2\right) $ if, and only if, (\ref{ASconv}) is verified.
\end{proof}

\bigskip

By using the Clark-Ocone formula stated in Proposition \ref{P : AbsCO}, we
deduce immediately, from\ Theorem \ref{T : Main 4}, a useful criterion for
the stable convergence of (Malliavin) differentiable functionals.

\begin{corollary}
Let $\mathfrak{H}_{n}$, $X_{n}\left( \mathfrak{H}_{n}\right) $, $\pi
^{\left( n\right) }$, $t_{n}$ and $\mathcal{U}_{n}$, $n\geq 1$, satisfy the
assumptions of Theorem \ref{T : Main 4} (in particular, (\ref{boxed})
holds), and consider a sequence of random variables $\left\{ F_{n}:n\geq
1\right\} $, such that $\mathbb{E}\left( F_{n}\right) =0$ and $F_{n}\in
\mathbb{D}_{X_{n}}^{1,2}$ for every $n$. Then, a sufficient condition to
have that
\begin{equation*}
F_{n}\rightarrow _{\left( s,\mathcal{U}^{\ast }\right) }\mathbb{E}\mu \left(
\cdot \right)
\end{equation*}%
and
\begin{equation*}
\mathbb{E}\left[ \exp \left( i\lambda F_{n}\right) \mid \mathcal{F}%
_{t_{n}}^{\pi ^{\left( n\right) }}\left( X_{n}\right) \right] \overset{%
\mathbb{P}}{\rightarrow }\exp \left( -\frac{\lambda ^{2}}{2}Y\right) ,\text{
\ \ }\forall \lambda \in
\mathbb{R}
\text{,}
\end{equation*}%
where $\mathcal{U}^{\ast }\triangleq \vee _{n}\mathcal{U}_{n}$, $Y\geq 0$ is
s.t. $Y\in \mathcal{U}^{\ast }$ and $\widehat{\mu }\left( \lambda \right)
=\exp \left( -\frac{\lambda ^{2}}{2}Y\right) $, $\forall \lambda \in \mathbb{%
R}$, is
\begin{equation}
\left\Vert \pi _{t_{n}}^{\left( n\right) }proj\left\{ D_{X_{n}}F_{n}\mid
L_{\pi ^{\left( n\right) }}^{2}\left( \mathfrak{H}_{n},X_{n}\right) \right\}
\right\Vert _{\mathfrak{H}_{n}}^{2}\overset{\mathbb{P}}{\rightarrow }0\text{
\ \ and \ \ }\left\Vert proj\left\{ D_{X_{n}}F_{n}\mid L_{\pi ^{\left(
n\right) }}^{2}\left( \mathfrak{H}_{n},X_{n}\right) \right\} \right\Vert _{%
\mathfrak{H}_{n}}^{2}\overset{\mathbb{P}}{\underset{n\rightarrow +\infty }{%
\rightarrow }}Y\text{.}  \label{CVMal}
\end{equation}
\end{corollary}

\begin{proof}
Since, for every $n$, $F_{n}$ is a centered random variable in $\mathbb{D}%
_{X_{n}}^{1,2}$, the abstract Clark-Ocone formula ensures that $F_{n}=\delta
_{X_{n}}\left( proj\left\{ D_{X_{n}}F_{n}\mid L_{\pi ^{\left( n\right)
}}^{2}\left( \mathfrak{H}_{n},X_{n}\right) \right\} \right) $, the result
follows from Theorem\ \ref{T : Main 4}, by putting
\begin{equation*}
u_{n}=proj\left\{ D_{X_{n}}F_{n}\mid L_{\pi ^{\left( n\right) }}^{2}\left(
\mathfrak{H}_{n},X_{n}\right) \right\} .
\end{equation*}
\end{proof}

\subsection{Example: a \textquotedblleft switching\textquotedblright\
sequence of quadratic Brownian functionals}

We conclude the paper by providing a generalization, as well as a new proof,
of a result contained in \cite[Proposition 2.1]{PecYor}. Let $W_{t}$, $t\in %
\left[ 0,1\right] $, be a standard Brownian motion initialized at zero, and
define the time-reversed Brownian motion $W^{\ast }$ by%
\begin{equation*}
W_{t}^{\ast }=W_{1}-W_{1-t}\text{, \ \ }t\in \left[ 0,1\right] \text{.}
\end{equation*}%
Observe that
\begin{equation}
W_{0}^{\ast }=0\text{ \ \ and \ \ }W_{1}^{\ast }=W_{1}.  \label{oo}
\end{equation}

We are interested in the asymptotic behavior, for $n\rightarrow +\infty $,
of the \textquotedblleft switching\textquotedblright\ sequence
\begin{equation*}
A_{n}=\int_{0}^{1}t^{2n}\left[ \left( W_{1}^{\left( n\right) }\right)
^{2}-\left( W_{t}^{\left( n\right) }\right) ^{2}\right] dt\text{, \ \ }n\geq
1\text{,}
\end{equation*}%
where%
\begin{equation*}
W^{\left( n\right) }=\left\{
\begin{array}{ll}
W & n\text{ odd} \\
W^{\ast } & n\text{ even.}%
\end{array}%
\right.
\end{equation*}

In particular, we would like to determine the speed at which $A_{n}$
converges to zero as $n\rightarrow +\infty $, by establishing a stable
convergence result. We start by observing that the asymptotic study of $%
A_{n} $ can be reduced to that of a sequence of double stochastic integrals.
As a matter of fact, from the relation $\left( W_{t}^{\left( n\right)
}\right) ^{2}=t+2\int_{0}^{t}W_{s}^{\left( n\right) }dW_{s}^{\left( n\right)
}$ one gets
\begin{equation*}
A_{n}=\int_{0}^{1}t^{2n}\left[ 2\int_{0}^{1}W_{s}^{\left( n\right) }\mathbf{1%
}_{\left( t\leq s\right) }dW_{s}^{\left( n\right) }+1-t\right] dt,
\end{equation*}%
and it is easily deduced that
\begin{eqnarray*}
\sqrt{n}\left( 2n+1\right) A_{n} &=&2\sqrt{n}\int_{0}^{1}dW_{s}^{\left(
n\right) }W_{s}^{\left( n\right) }s^{2n+1}+\sqrt{n}\left( 2n+1\right)
\int_{0}^{1}\left( 1-t\right) t^{2n}dt \\
&=&2\sqrt{n}\int_{0}^{1}dW_{s}^{\left( n\right) }W_{s}^{\left( n\right)
}s^{2n+1}+o\left( 1\right) .
\end{eqnarray*}

Now define $\sigma \left( W\right) $ to be the $\sigma $-field generated by $%
W$ (or, equivalently, by $W^{\ast }$): we have the following

\begin{proposition}
\label{P : Switch}As $n\rightarrow +\infty $,%
\begin{equation}
\sqrt{n}\left( 2n+1\right) A_{n}\rightarrow _{\left( s,\sigma \left(
W\right) \right) }\mathbb{E}\mu _{1}\left( \cdot \right) \text{,}
\label{EX1}
\end{equation}%
where $\mu _{1}\left( \cdot \right) $ verifies, for $\lambda \in \mathbb{R}$,%
\begin{equation*}
\widehat{\mu }_{1}\left( \lambda \right) =\exp \left( -\frac{\lambda ^{2}}{2}%
W_{1}^{2}\right) \text{,}
\end{equation*}%
or, equivalently, for every $Z\in \sigma \left( W\right) $
\begin{equation}
\left( Z,\sqrt{n}\left( 2n+1\right) A_{n}\right) \overset{\text{law}}{%
\rightarrow }\left( Z,W_{1}\times N^{\prime }\right) ,  \label{exoCL}
\end{equation}%
where $N^{\prime }$ is a standard Gaussian random variable independent of $W$%
.
\end{proposition}

\bigskip

\textbf{Remark -- }In particular, if $W^{\left( n\right) }=W$ for every $n$
(no switching between $W$ and $W^{\ast }$), one gets the same convergence in
law (\ref{exoCL}). This last result was proved in \cite[Proposition 2.1]%
{PecYor} by completely different methods.

\bigskip

\textbf{Proof of Proposition \ref{P : Switch}. }The proof of (\ref{EX1}) is
based on Theorem \ref{T : Main 4}. First observe that the Gaussian family
\begin{equation}
X_{W}\left( h\right) =\int_{0}^{1}h\left( s\right) dW_{s}\text{, \ \ }h\in
L^{2}\left( \left[ 0,1\right] ,ds\right) \text{,}  \label{IsoG}
\end{equation}%
defines an isonormal Gaussian process over the Hilbert space $\mathfrak{H}%
_{1}\triangleq L^{2}\left( \left[ 0,1\right] ,ds\right) $; we shall write $%
X_{W}$ to indicate the isonormal Gaussian process given by (\ref{IsoG}). Now
define the following sequence of continuous resolutions of the identity on $%
\mathfrak{H}_{1}$: for every $n\geq 1$, every $t\in \left[ 0,1\right] $ and
every $h\in $ $\mathfrak{H}_{1}$%
\begin{equation*}
\pi _{t}^{\left( n\right) }h=\left\{
\begin{array}{ll}
\pi _{t}^{o}h\triangleq h\mathbf{1}_{\left[ 0,t\right] }\text{,} & n\text{
odd} \\
\pi _{t}^{e}h\triangleq h\mathbf{1}_{\left[ 1-t,1\right] }\text{,} & n\text{
even.}%
\end{array}%
\right. ,
\end{equation*}%
so that, for $t\in \left[ 0,1\right] $
\begin{equation}
\mathcal{F}_{t}^{\pi ^{\left( n\right) }}\left( X_{W}\right) =\left\{
\begin{array}{ll}
\mathcal{F}_{t}^{\pi ^{o}}\left( X_{W}\right) \triangleq \sigma \left\{
W_{u}:u\leq t\right\} \text{,} & n\text{ odd} \\
\mathcal{F}_{t}^{\pi ^{e}}\left( X_{W}\right) \triangleq \sigma \left\{
W_{1}-W_{s}:s\geq 1-t\right\} \text{,} & n\text{ even.}%
\end{array}%
\right. .  \label{filtr}
\end{equation}%
In this case, the class of adapted processes $L_{\pi ^{o}}^{2}\left(
\mathfrak{H}_{1},X_{W}\right) $ (resp. $L_{\pi ^{e}}^{2}\left( \mathfrak{H}%
_{1},X_{W}\right) $) is given by those elements of $L^{2}\left( \mathfrak{H}%
_{1},X_{W}\right) $ that are adapted to the filtration $\mathcal{F}_{\cdot
}^{\pi ^{o}}\left( X_{W}\right) $ (resp. $\mathcal{F}_{\cdot }^{\pi
^{e}}\left( X_{W}\right) $), as defined in (\ref{filtr}). Moreover, by
defining, for $n\geq 1$,%
\begin{equation*}
u_{n}\left( t\right) =2\sqrt{n}W_{t}^{\left( n\right) }t^{2n+1}\text{, \ \ }%
t\in \left[ 0,1\right] \text{,}
\end{equation*}%
we see immediately that $u_{n}\in L_{\pi ^{\left( n\right) }}^{2}\left(
\mathfrak{H}_{1},X_{W}\right) $ for every $n$, and hence%
\begin{equation*}
2\sqrt{n}\int_{0}^{1}W_{s}^{\left( n\right) }s^{2n+1}dW_{s}^{\left( n\right)
}=\int_{0}^{1}u_{n}\left( s\right) dW_{s}^{\left( n\right) }=\delta
_{W}\left( u_{n}\right) \text{,}
\end{equation*}%
where $\delta _{W}$ stands for a Skorohod integral with respect to $X_{W}$.
Indeed, if $n$ is even, $\int_{0}^{1}u_{n}\left( s\right) dW_{s}^{\left(
n\right) }$ is an It\^{o} integral with respect to $W$, and, when $n$ is
odd, $\int_{0}^{1}u_{n}\left( s\right) dW_{s}^{\left( n\right) }$ is again
an It\^{o} integral with respect to a time reversed Brownian motion (see
e.g. \cite[Section 4]{PTT} for more details). Now fix $\varepsilon \in
\left( 1/2,1\right) $, and set $t_{n}=\varepsilon ^{1/\sqrt{n}}$, so that, $%
t_{n}>1/2$ and $t_{n}\uparrow 1.$ Then,%
\begin{equation*}
\mathbb{E}\left[ \left\Vert \pi _{t_{n}}^{\left( n\right) }u_{n}\right\Vert
_{\mathfrak{H}_{1}}^{2}\right] =4n\int_{0}^{t_{n}}s^{4n+3}ds=\frac{4n}{4n+4}%
\varepsilon ^{\frac{4n+4}{\sqrt{n}}}\underset{n\rightarrow +\infty }{%
\longrightarrow }0,
\end{equation*}%
and%
\begin{eqnarray}
\left\Vert u_{n}\right\Vert _{\mathfrak{H}_{1}}^{2}
&=&4n\int_{0}^{1}ds\left( W_{s}^{\left( n\right) }\right) ^{2}s^{4n+2}
\notag \\
&=&\frac{4n}{4n+4}+\frac{8n}{4n+3}\int_{0}^{1}dW_{u}^{\left( n\right)
}W_{u}^{\left( n\right) }\left( 1-u^{4n+3}\right)  \notag \\
&=&o_{\mathbb{P}}\left( 1\right) +1+\frac{8n}{4n+3}\int_{0}^{1}dW_{u}^{%
\left( n\right) }W_{u}^{\left( n\right) }\overset{\mathbb{P}}{\rightarrow }%
W_{1}^{2}=\left( W_{1}^{\ast }\right) ^{2},  \label{CVex}
\end{eqnarray}%
by (\ref{oo}), where $o_{\mathbb{P}}\left( 1\right) $ stands for a sequence
converging to zero in probability (as $n\rightarrow +\infty $). We thus have
shown that \ relations (\ref{ASnegl}) and (\ref{ASconv}) of Theorem \ref{T :
Main 4} are satisfied. It remains to verify relation (\ref{boxed}), namely
to show that there exists a sequence of $\sigma $-fields $\left\{ \mathcal{U}%
_{n}:n\geq 1\right\} $ verifying $\mathcal{U}_{n}\subseteq \mathcal{U}%
_{n+1}\cap \mathcal{F}_{t_{n}}^{\pi ^{\left( n\right) }}\left( X_{W}\right) $
and $\vee _{n}\mathcal{U}_{n}=\sigma \left( W\right) $. The sequence
\begin{equation*}
\mathcal{U}_{n}\triangleq \sigma \left\{ W_{u}-W_{s}:1-t_{n}\leq s<u\leq
t_{n}\right\} \text{,}
\end{equation*}%
which is increasing and such that $\mathcal{U}_{n}\subseteq \mathcal{F}%
_{t_{n}}^{\pi ^{\left( n\right) }}\left( X_{W}\right) $ (see (\ref{filtr})),
verifies these properties. Therefore, Theorem \ref{T : Main 4} applies, and
we obtain the stable convergence result (\ref{EX1}). $\blacksquare $

\bigskip

\textbf{Remarks -- }(a) The sequence $A_{n}^{\prime }\triangleq \sqrt{n}%
\left( 2n+1\right) A_{n}$, although stably convergent and such that (\ref%
{CVex}) is verified, \textit{does not} admit a limit in probability. Indeed,
simple computations show that $A_{n}^{\prime }$ is not a Cauchy sequence in $%
L^{2}\left( \mathbb{P}\right) $ and therefore, since the $L^{2}$ and $L^{0}$
topologies coincide on any finite sum of Wiener chaoses (see e.g. \cite{Sch}%
), $A_{n}^{\prime }$ cannot converge in probability.

(b) Observe that, by using the notation introduced above (see e.g. \ref%
{exPROJ}), $rank\left( \pi ^{o}\right) =rank\left( \pi ^{e}\right) =1$%
.\bigskip

\bigskip

\textbf{Acknowledgments -- }G. Peccati thanks M.S. Taqqu for an invitation
to the Department of Mathematics and Statistics of the Boston University in
February 2006. G. Peccati is also indebted to S. Kwapien for a series of
inspiring discussions around the \textquotedblleft principle of
conditioning\textquotedblright , in the occasion of a visit to the Warsaw
University in October 2004.

\end{document}